\documentclass[reqno]{amsart}
\usepackage{amssymb,mathrsfs,amscd}

\makeatletter

\newcommand{\g}{\mathfrak g}
\newcommand{\id}{\operatorname{id}}
\newcommand{\Wh}{\operatorname{Wh}}
\newcommand{\gr}{\operatorname{gr}}
\newcommand{\ev}{\operatorname{ev}}

\newcommand{\s}{\mathfrak s}

\newcommand{\eps}{\varepsilon}
\newcommand{\col}{\operatorname{col}}
\newcommand{\row}{\operatorname{row}}

\newcommand{\C}{\mathbb C}

\newcommand{\Z}{\mathbb Z}
\newcommand{\N}{\mathbb N}
\newcommand{\V}{\mathbb V}

\newcommand{\Idem}{\operatorname{Idem}}
\newcommand{\Tab}{\operatorname{Tab}}
\newcommand{\Std}{\operatorname{Std}}
\newcommand{\Col}{\operatorname{Col}}

\newcommand{\End}{\operatorname{End}}
\newcommand{\Hom}{\operatorname{Hom}}
\renewcommand{\hom}{\operatorname{Hom}}
\newcommand{\op}{\operatorname{op}}

\newtheorem{thm}{Theorem}[section]
\newtheorem{cor}[thm]{Corollary}
\newtheorem{lem}[thm]{Lemma}

\theoremstyle{definition}

\theoremstyle{remark}
\newtheorem{rem}[thm]{Remark}

\numberwithin{equation}{section}

\newcommand{\bc}{{\text{\boldmath{$c$}}}}

\newcommand{\bp}{{\text{\boldmath{$p$}}}}
\newcommand{\bq}{{\text{\boldmath{$q$}}}}
\newcommand{\brr}{{\text{\boldmath{$r$}}}}
\newcommand{\bs}{{\text{\boldmath{$s$}}}}
\newcommand{\bt}{{\text{\boldmath{$t$}}}}
\newcommand{\bh}{{\text{\boldmath{$h$}}}}
\newcommand{\bi}{{\text{\boldmath{$i$}}}}
\newcommand{\bid}{{\text{\boldmath{$1$}}}}
\newcommand{\biz}{{\text{\boldmath{$0$}}}}
\newcommand{\bj}{{\text{\boldmath{$j$}}}}
\newcommand{\bk}{{\text{\boldmath{$k$}}}}
\newcommand{\bm}{{\text{\boldmath{$m$}}}}
\newcommand{\bmu}{{\text{\boldmath{$\mu$}}}}
\newcommand{\bu}{{\text{\boldmath{$u$}}}}
\newcommand{\bv}{{\text{\boldmath{$v$}}}}
\newdimen\Hoogte    \Hoogte=12pt
\newdimen\Breedte   \Breedte=12pt
\newdimen\Dikte     \Dikte=0.5pt
\newenvironment{Young}{\begingroup
       \def\vr{\vrule height0.8\Hoogte width\Dikte depth 0.2\Hoogte}
       \def\fbox##1{\vbox{\offinterlineskip
                    \hrule height\Dikte
                    \hbox to \Breedte{\vr\hfill##1\hfill\vr}
                    \hrule height\Dikte}}
       \vbox\bgroup \offinterlineskip \tabskip=-\Dikte \lineskip=-\Dikte
            \halign\bgroup &\fbox{##\unskip}\unskip  \crcr }
       {\egroup\egroup\endgroup}
\def\Diagram#1{\relax\ifmmode\vcenter{\,\begin{Young}#1\end{Young}\,}\else%
              $\vcenter{\,\begin{Young}#1\end{Young}\,}$\fi}

\begin{document}
\title[Schur-Weyl duality]{\boldmath Schur-Weyl duality for higher levels}
\author{Jonathan Brundan}
\address{Department of Mathematics\\University of Oregon\\Eugene OR 97403\\USA}
\email{brundan@uoregon.edu}
\author{Alexander Kleshchev}
\address{Department of Mathematics\\University of Oregon\\Eugene OR 97403\\USA}
\email{klesh@uoregon.edu}
\thanks{Research supported in part by NSF grant no. DMS-0654147.}

\begin{abstract}
We extend Schur-Weyl duality
to an arbitrary level $l \geq 1$, 
level one recovering the classical duality
between the symmetric and general linear groups.
In general,
the symmetric group is replaced by
the degenerate cyclotomic Hecke algebra over $\C$ parametrized by a 
dominant weight of level $l$ for the root system of type $A_\infty$.
As an application, we prove that the degenerate analogue of the 
quasi-hereditary cover of the 
cyclotomic Hecke algebra constructed by Dipper, James and Mathas is Morita equivalent to certain blocks of 
parabolic category $\mathcal{O}$ for the general linear Lie algebra.
\end{abstract}

\maketitle
\section{Introduction}

We will work over the ground field $\C$ (any algebraically closed field of
characteristic zero is fine too).
Recall the definition of the {\em degenerate affine Hecke algebra}
$H_d$ from \cite{D2}. This is the associative algebra 
equal as a vector space to 
the tensor product 
$\C[x_1,\dots,x_d] \otimes \C S_d$ of a polynomial algebra 
and the group algebra of the symmetric group $S_d$. 
We write $s_i$ for the basic transposition 
$(i\:i\!+\!1) \in S_d$, and identify $\C[x_1,\dots,x_d]$ and
$\C S_d$ with the subspaces $\C[x_1,\dots,x_d] \otimes 1$ and
$1 \otimes \C S_d$ of $H_d$, respectively. Multiplication is then defined 
so that
$\C[x_1,\dots,x_d]$ and $\C S_d$ are subalgebras of $H_d$,
$s_i x_j = x_j s_i$ if $j \neq i,i+1$, and
$$
s_i x_{i+1} = x_i s_i + 1.
$$
As well as being a subalgebra of $H_d$,
the group algebra $\C S_d$ is also a quotient, via the
homomorphism $H_d \twoheadrightarrow \C S_d$ mapping
$x_1 \mapsto 0$ and 
$s_i \mapsto s_i$ for each $i$.

The degenerate affine Hecke algebra arises when studying the
representation theory of the general linear Lie algebra 
$\mathfrak{g} = \mathfrak{gl}_N(\C)$, as follows.
Let $V$ be the natural $\mathfrak{g}$-module of column vectors,
let $M$ be any $\mathfrak{g}$-module, and consider the
tensor product $M \otimes V \otimes \cdots \otimes V$ ($d$ copies of $V$).
Let $e_{i,j} \in \mathfrak{g}$ 
denote the $ij$-matrix unit and
$\Omega = \sum_{i,j=1}^N e_{i,j} \otimes e_{j,i}$.
Then, by an observation of Arakawa and Suzuki \cite[$\S$2.2]{AS}, there is a well-defined right action of $H_d$ 
on $M \otimes V^{\otimes d}$, commuting with the left
action of $\mathfrak{g}$,
defined
so that $x_1$ acts as the endomorphism $\Omega \otimes 1^{\otimes d-1}$
and each
$s_i$ acts as $1^{\otimes i} 
\otimes \Omega \otimes 1^{\otimes (d-i-1)}$.

The goal of this article is to extend the classical 
Schur-Weyl duality to some remarkable finite dimensional quotients of
$H_d$: the {\em degenerate cyclotomic Hecke algebras}.
Let $\Lambda = \sum_{i \in \Z} m_i \Lambda_i$
be 
a dominant weight of level $l$ 
for the root system of type $A_\infty$,
so each $m_i$ is a non-negative integer
and $\sum_{i \in \Z} m_i = l$.
The degenerate cyclotomic Hecke algebra
parametrized by $\Lambda$ is the quotient
$H_d(\Lambda)$
of $H_d$ by the two-sided ideal generated by the
polynomial 
$\prod_{i \in \Z} (x_1 - i)^{m_i}$.
If $\Lambda = \Lambda_0$ (or any other weight of level $l=1$)
then $H_d(\Lambda) \cong \C S_d$ as above.

\vspace{2mm}

Fix now a partition $\lambda = (p_1 \leq \cdots \leq p_n)$ of $N$,
let $(q_1 \geq \cdots \geq q_l)$ be the transpose
partition, and take the dominant weight $\Lambda$ from the previous paragraph
to be the weight
$\Lambda = \sum_{i = 1}^{l} 
\Lambda_{q_i-n}$. Up to change of origin,
any $\Lambda$ arises in this way for suitable $\lambda$.
 We identify $\lambda$ with its Young diagram, numbering
rows by $1,2,\dots,n$ from top to bottom
in order of increasing length and columns 
by $1,2,\dots,l$ from left to right in order of decreasing length,
so that 
there are $p_i$ boxes in the $i$th row and $q_j$ boxes in the $j$th column.
Since these are not the standard conventions, we give an example.
Suppose $\lambda = (p_1,p_2,p_3) = (2,3,4)$, so 
$(q_1,q_2,q_3,q_4)=(3,3,2,1)$, 
$\Lambda = 2 \Lambda_0 + \Lambda_{-1}+\Lambda_{-2}$,
$N= 9$, $n=3$, and
the level $l = 4$.
The Young diagram is
$$
\Diagram{1&4\cr2&5&7\cr3&6&8&9\cr}
\begin{picture}(0,0)
\put(0,0){\makebox(7,-25){.}}
\end{picture}
$$
We always number the boxes of the diagram $1,2,\dots,N$ down columns starting
from the first column,
and write $\row(i)$ and 
$\col(i)$ for the row and column
numbers of the $i$th box. 
This identifies the boxes
with 
the standard basis $v_1,\dots, v_N$ of the natural $\mathfrak{g}$-module
$V$. 
Define
 $e \in \mathfrak{g}$ to be the nilpotent matrix of Jordan type $\lambda$
which maps the basis vector corresponding 
the $i$th box to the one immediately to its left, or to zero if there
is no such box;
in our example, $e = e_{1,4} + e_{2,5}+e_{5,7}+e_{3,6}+e_{6,8}+e_{8,9}$.
Finally, define a $\Z$-grading 
$\g = \bigoplus_{r \in \Z} \g_r$ 
on $\mathfrak{g}$ by
declaring that $e_{i,j}$ is of degree $(\col(j)-\col(i))$
for each $i,j=1,\dots,N$, and set
$\mathfrak{p} = \bigoplus_{r \geq 0} \mathfrak{g}_r,
\mathfrak{h} = \mathfrak{g}_0,
\mathfrak{m} = \bigoplus_{r < 0} \mathfrak{g}_r.
$
The element $e \in \mathfrak{g}_1$ is a Richardson element 
in the parabolic subalgebra $\mathfrak{p}$, so the
centralizer $\mathfrak{g}_e$ of $e$ in $\mathfrak{g}$
is a graded subalgebra of $\mathfrak{p}$.
Hence its universal enveloping algebra $U(\mathfrak{g}_e)$ is a graded
subalgebra of $U(\mathfrak{p})$.

Our first main result is best understood as a
filtered deformation of the following
extension of the classical Schur-Weyl duality
to the centralizer $\mathfrak{g}_e$.
Let $\C_l[x_1,\dots,x_d]$ be the level $l$ truncated polynomial algebra, that is, 
the quotient of
the polynomial algebra $\C[x_1,\dots,x_d]$
by the relations $x_1^l = \cdots = x_d^l = 0$.
Extend the usual right action of $S_d$ on $V^{\otimes d}$ to an action of 
the
twisted tensor product
$\C_l[x_1,\dots,x_d] \,{\scriptstyle{\rtimes\!\!\!\!\!\bigcirc}}\, \C S_d$ 
by defining the action of $x_i$ 
to be as the endomorphism
$1^{\otimes(i-1)}\otimes e \otimes 1^{\otimes (d-i)}$.
This commutes with the natural action of $\mathfrak{g}_e$, 
making $V^{\otimes d}$ into a $(U(\mathfrak{g}_e),
\C_l[x_1,\dots,x_d]
 \,{\scriptstyle{\rtimes\!\!\!\!\!\bigcirc}}\, \C S_d)$-bimodule.
We have defined a homomorphism $\phi_d$ and an antihomomorphism $\psi_d$:
$$
U(\mathfrak{g}_e) \stackrel{\phi_d}{\longrightarrow}
\End_{\C}(V^{\otimes d}) \stackrel{\psi_d}{\longleftarrow}
\C_l[x_1,\dots,x_d]
 \,{\scriptstyle{\rtimes\!\!\!\!\!\bigcirc}}\, \C S_d.
$$
Now a result of Vust \cite[$\S$6]{KP}
asserts that
these maps satisfy the double centralizer property,
i.e. the image of $\phi_d$ is the centralizer of the image of $\psi_d$
and vice versa.
This is a surprising consequence of the normality of the
closure of the conjugacy class containing $e$.

In our filtered deformation of this picture, 
we replace
$\C_l[x_1,\dots,x_d]
 \,{\scriptstyle{\rtimes\!\!\!\!\!\bigcirc}}\, \C S_d$
with the degenerate cyclotomic Hecke algebra $H_d(\Lambda)$
and $U(\mathfrak{g}_e)$ with the finite $W$-algebra $W(\lambda)$
associated to the partition $\lambda$.
Let us recall the definition of the latter algebra
following \cite{BKrep};
see also \cite{P, premet2, GG, BGK}.
Let $\eta:U(\mathfrak{p}) \rightarrow U(\mathfrak{p})$
be the algebra automorphism defined by
\begin{align*}
\eta(e_{i,j}) &= 
e_{i,j} + \delta_{i,j} (n-q_{\col(j)} - q_{\col(j)+1} - \cdots - q_l)
\end{align*}
for each $e_{i,j} \in \mathfrak{p}$.
Let $I_\chi$ be the kernel of the homomorphism $\chi:U(\mathfrak{m}) 
\rightarrow \C$ defined by
$x \mapsto (x,e)$ for all $x \in \mathfrak{m}$, where $(.,.)$ is the trace form
on $\mathfrak{g}$. 
Then, by the definition followed here, $W(\lambda)$
is the following subalgebra of $U(\mathfrak{p})$:
\begin{align*}
W(\lambda) 
&= \{u \in U(\mathfrak{p})\mid
[x, \eta(u)] \in U(\mathfrak{g}) I_\chi \text{ for all }x \in \mathfrak{m}\}.
\end{align*}
The grading on $U(\mathfrak{p})$ induces a filtration on
$W(\lambda)$ so that the associated graded algebra
$\gr W(\lambda)$ is naturally identified with a graded subalgebra
of $U(\mathfrak{p})$; this is the {\em loop} or {\em good filtration}
of the finite $W$-algebra, not the
more familiar
Kazhdan filtration.
A key point is that, by a result of Premet \cite[Proposition 2.1]{premet2},
the associated graded algebra $\gr W(\lambda)$
is equal to $U(\mathfrak{g}_e)$ as a graded 
subalgebra of $U(\mathfrak{p})$,
hence $W(\lambda)$ is indeed a filtered deformation of $U(\mathfrak{g}_e)$.

We still need to introduce a 
$(W(\lambda), H_d(\Lambda))$-bimodule
$V^{\otimes d}$
 which is a filtered deformation of
the $(U(\mathfrak{g}_e), \C_l[x_1,\dots,x_d]
 \,{\scriptstyle{\rtimes\!\!\!\!\!\bigcirc}}\, \C S_d)$-bimodule $V^{\otimes d}$ from above.
The left action of $W(\lambda)$ on $V^{\otimes d}$ is simply the restriction
of the natural action of $U(\mathfrak{p})$.
To define the right action of $H_d(\Lambda)$,
let 
$S_d$ act on the right 
by place permutation as usual.
Let $x_1$ act as the endomorphism
$$
\left(e + \sum_{j=1}^N (q_{\col(j)}-n) e_{j,j} \right)
\otimes 1^{\otimes (d-1)}
-
\sum_{k=2}^d \sum_{\substack{i,j=1 \\ \col(i) < \col(j)}}^N
e_{i,j} \otimes 1^{\otimes(k-2)} \otimes e_{j,i} \otimes 1^{\otimes (d-k)}.
$$
This extends uniquely to make 
$V^{\otimes d}$ into a 
 $(W(\lambda),H_d(\Lambda))$-bimodule.
To explain why,
let $Q_\chi$ denote the $\mathfrak{g}$-module
$U(\mathfrak{g}) / U(\mathfrak{g}) I_\chi$.
It is actually a $(U(\mathfrak{g}), W(\lambda))$-bimodule, the action of $u \in W(\lambda)$ arising by right
multiplication by $\eta(u)$ (which leaves $U(\mathfrak{g}) I_\chi$ invariant
by the above definition of $W(\lambda)$).
Let $\mathcal{C}(\lambda)$ be the category of
{\em generalized Whittaker modules}, that is, the $\mathfrak{g}$-modules
on which $(x - \chi(x))$ acts locally nilpotently for all
$x \in \mathfrak{m}$.
By
{\em Skryabin's theorem} \cite{Skry},
the functor 
$Q_\chi \otimes_{W(\lambda)} ?: W(\lambda)\text{-mod} \rightarrow \mathcal{C}(\lambda)$
is an equivalence of categories. 
The $W(\lambda)$-module $V^{\otimes d}$ 
corresponds under this equivalence
to the $\mathfrak{g}$-module
$M \otimes V^{\otimes d}$,
where
$M = Q_\chi \otimes_{W(\lambda)} \C$ and $\C$ is the
restriction of the trivial $U(\mathfrak{p})$-module.
The above formula for the action of $x_1$ 
on $V^{\otimes d}$ arises by
transporting its action
 on 
$M \otimes V^{\otimes d}$ from \cite[$\S$2.2]{AS}
through Skryabin's equivalence.

We have now defined a homomorphism $\Phi_d$ and an antihomomorphism
$\Psi_d$
$$ 
W(\lambda) \stackrel{\Phi_d}{\longrightarrow} \End_{\C}(V^{\otimes d})
\stackrel{\Psi_d}{\longleftarrow} H_d(\Lambda).
$$
Let $W_d(\lambda)$ denote the image of 
$\Phi_d$.
This finite dimensional algebra is
a natural
analogue of the classical {\em Schur algebra} for higher levels.
Let $H_d(\lambda)$ denote
the image of the homomorphism
$\Psi_d:H_d(\Lambda) \rightarrow \End_{\C}(V^{\otimes d})^{\op}$,
so that $V^{\otimes d}$ is also a $(W_d(\lambda), H_d(\lambda))$-bimodule.
Actually, if at least $d$ parts of $\lambda$ are equal to $l$,
then the map $\Psi_d$ is injective so $H_d(\lambda) = H_d(\Lambda)$.
In general $H_d(\lambda)$ is 
a sum of certain blocks of the finite dimensional algebra
$H_d(\Lambda)$.
In view of the fact that $W_d(\lambda)$ and $H_d(\lambda)$ are usually
not semisimple algebras, the following result was to us quite surprising.

\vspace{2mm}
{\noindent\bf Theorem A.}
{\em
The maps $\Phi_d$ and $\Psi_d$ satisfy the double centralizer property, i.e.
\begin{align*}
W_d(\lambda)=&\End_{H_d(\lambda)}(V^{\otimes d}),\\
&\End_{W_d(\lambda)}(V^{\otimes d})^{\op} = H_d(\lambda).
\end{align*}
Moreover, 
the functor
\begin{align*}
\Hom_{W_d(\lambda)}(V^{\otimes d}, ?):
&\,W_d(\lambda)\text{\rm-mod}
\rightarrow H_d(\lambda)\text{\rm-mod}
\end{align*}
is an equivalence of categories.
}
\vspace{2mm}

Our second main result is concerned with
the parabolic analogue of the BGG category
$\mathcal{O}$ for the Lie algebra $\mathfrak{g}$ relative to the
subalgebra $\mathfrak{p}$.
Let $P$ denote the
module $U(\mathfrak{g}) \otimes_{U(\mathfrak{p})} 
\C_{-\rho}$
induced from the one dimensional $\mathfrak{p}$-module $\C_{-\rho}$
on which each $e_{i,j} \in \mathfrak{p}$ acts as 
$\delta_{i,j}(q_1+q_2+\cdots+q_{\col(j)} - n)$.
This is an irreducible projective 
module in parabolic category $\mathcal{O}$.
Let $\mathcal{O}^d(\lambda)$ denote the 
Serre subcategory of parabolic category $\mathcal{O}$
generated by the module
$P \otimes V^{\otimes d}$.
We note that $\mathcal{O}^d(\lambda)$  is a sum of certain integral blocks of parabolic category
$\mathcal{O}$, and every integral block is equivalent to a block of
$\mathcal{O}^d(\lambda)$ for sufficiently large $d$.
Moreover, 
the module $P \otimes V^{\otimes d}$
is a {\em self-dual projective module} in $\mathcal{O}^d(\lambda)$, and every
self-dual projective indecomposable module in $\mathcal{O}^d(\lambda)$
is a summand of $P \otimes V^{\otimes d}$. In equivalent language, 
$P \otimes V^{\otimes d}$ is a {\em prinjective
generator} for $\mathcal{O}^d(\lambda)$ where by a {\em prinjective
module} we mean a module that is
both projective and injective.
Applying the construction from \cite[$\S$2.2]{AS} once more, we can view
$P \otimes V^{\otimes d}$ as a $(\mathfrak{g}, H_d)$-bimodule.
It turns out that the
 right action of $H_d$ on
$P \otimes V^{\otimes d}$ 
factors through the quotient $H_d(\lambda)$ of $H_d$ to make
$P \otimes V^{\otimes d}$ into a faithful right 
$H_d(\lambda)$-module,
i.e. $H_d(\lambda) \hookrightarrow \End_{\C}(P \otimes V^{\otimes d})^{\op}$.

\vspace{2mm}

{\noindent\bf Theorem B.} {\em
$\End_{\mathfrak{g}}(P \otimes V^{\otimes d})^{\operatorname{op}}
= H_d(\lambda).
$
}

\vspace{2mm}

The link between Theorems A and B is provided by the
{\em Whittaker functor} 
$$
\mathbb{V}: \mathcal{O}^d(\lambda)
\rightarrow W_d(\lambda)\text{-mod}
$$ 
introduced originally by
Kostant and Lynch \cite{Lynch} and studied recently in \cite[$\S$8.5]{BKrep}:
the $W_d(\lambda)$-module
$V^{\otimes d}$ from above is isomorphic to 
$\V(P \otimes V^{\otimes d})$.
We actually show
that
$\Hom_{\mathfrak{g}}(P \otimes V^{\otimes d}, ?)
\cong
\Hom_{W_d(\lambda)}(V^{\otimes d}, ?) \circ \V$,
i.e. the following diagram of functors commutes up to isomorphism:
$$
\begin{CD}
&\:\mathcal{O}^d(\lambda) &\\
\\\\
W_d(\lambda)\text{-mod}@>\sim > \Hom_{W_d(\lambda)}(V^{\otimes d}, ?)> &
H_d(\lambda)\text{-mod.}
\end{CD}
\begin{picture}(0,0)
\put(-143,-8){\makebox(0,0){$\swarrow$}}
\put(-135,12){\makebox(0,0){$\scriptstyle\mathbb{V}$}}
\put(-140,-5){\line(1,1){33}}
\put(-43.9,-8){\makebox(0,0){$\searrow$}}
\put(-47,-5){\line(-1,1){33}}
\put(-22,12){\makebox(0,0){$\scriptstyle\Hom_{\mathfrak{g}}(P \otimes V^{\otimes d}, ?)$}}
\end{picture}
$$
The categories $W_d(\lambda)\text{-mod}$ and $H_d(\lambda)\text{-mod}$ 
thus give two
different realizations of a
natural quotient of the category $\mathcal{O}^d(\lambda)$
in the general sense of \cite[$\S$III.1]{Gab}, 
the respective quotient functors being
the Whittaker functor $\mathbb{V}$ and the functor
$\Hom_{\mathfrak{g}}(P \otimes V^{\otimes d}, ?)$.

In many circumstances, the Whittaker functor 
turns out to be easier to 
work with than the functor $\Hom_{\mathfrak{g}}(P \otimes V^{\otimes d}, ?)$,
so this point of view
facilitates various other important computations
regarding the relationship between $\mathcal{O}^d(\lambda)$
and $H_d(\lambda)$-mod.
For example, 
we use it to identify
the images of arbitrary projective indecomposable modules in $\mathcal{O}^d(\lambda)$ with  the indecomposable summands
of the degenerate analogues of the
permutation modules introduced by Dipper, James and 
Mathas \cite{DJM}.
Also, we identify the images of parabolic Verma modules with {\em Specht modules}, thus recovering formulae for 
the latter's composition multiplicities directly from the Kazhdan-Lusztig
conjecture for $\mathfrak{g}$.
The degenerate analogue of Ariki's categorification theorem
from \cite{Ariki} follows
as an easy consequence of these results, as we will explain
in more detail in a subsequent article \cite{BKnew}.

Theorem B should be compared with Soergel's Endomorphismensatz from \cite{So};
the connection with Whittaker modules in that case is due to
Backelin \cite{Ba}.
As observed originally by Stroppel \cite[Theorem 10.1]{St}
(as an application of \cite[Theorem 2.10]{KSX}), 
there is also a version of Soergel's Struktursatz 
for our parabolic setup:
the functor $\Hom_{\mathfrak{g}}(P \otimes V^{\otimes d}, ?)$
is fully faithful on projective objects.
From this, we obtain 
the third main result of the article.

\vspace{2mm}

{\noindent\bf Theorem C.} {\em
If at least $d$ parts of $\lambda$ are equal to $l$,
there is an equivalence of categories between
$\mathcal{O}^d(\lambda)$ and the
category of finite dimensional modules over a degenerate analogue of the 
cyclotomic $q$-Schur algebra of Dipper, James and Mathas.
}

\vspace{2mm}

The proof of Theorem C explains how the category $\mathcal{O}^d(\lambda)$
itself can be 
reconstructed from the algebra $H_d(\lambda)$.
Unfortunately, this approach does not give a direct construction 
of the individual blocks of $\mathcal{O}^d(\lambda)$, 
rather, we obtain all the blocks simultaneously.
For this reason, 
it is important to find some more explicit information
about the individual
blocks of the degenerate cyclotomic Hecke algebras, especially their centers;
see \cite{Kh, St2} for further motivation for such questions.
As an application of the Schur-Weyl duality for higher levels
developed here,
the first steps in this direction have recently been made
in \cite{cyclo,springer, St3}: the centers of the blocks
of degenerate cyclotomic Hecke algebras (and of the parabolic category
$\mathcal{O}^d(\lambda)$) are now known explicitly
and are isomorphic to cohomology algebras of
Spaltenstein varieties.

\vspace{2mm}

Let us now explain the organization of the remainder of the article.
\begin{itemize}
\item
In section 2, we introduce some notation for elements of the centralizer
$\g_e$ and then review Vust's double centralizer property.
\item
In section 3, we define the algebras $W_d(\lambda)$ and $H_d(\lambda)$
and prove the filtered version of the double centralizer property, ultimately 
as a consequence of Vust's result at the level of
associated graded algebras.
\item
In section 4, we review some known results about parabolic category $\mathcal{O}$ and the subcategory $\mathcal{O}^d(\lambda)$. In particular, we reformulate 
Irving's classification  \cite{I} of self-dual projective indecomposable modules 
in terms of standard tableaux, and we introduce some natural 
projective modules in $\mathcal{O}^d(\lambda)$ built from divided powers.
\item
In section 5, we prove Theorems A and B, modulo one technical fact the proof of which is deferred to section 6. We also formulate a characterization of the
Whittaker functor $\V$ which is of independent interest, and we explain how to recover the degenerate analogue of a theorem of Dipper and Mathas \cite{DM} in our framework. 
\item
In section 6, we use
certain weight idempotents in
$W_d(\lambda)$ to relate the projective modules in
$\mathcal{O}^d(\lambda)$ built from divided powers to the permutation modules of Dipper, James and Mathas, completing the proof of Theorem C.
In similar fashion, we relate parabolic Verma modules to Specht modules. 
\item
Finally there is an appendix giving a short proof of the fact that the degenerate cyclotomic Hecke algebras are symmetric algebras.
\end{itemize}

\vspace{2mm}
\noindent{\em Acknowledgements.}
We thank Steve Donkin for drawing our attention to \cite{KP} and
Catharina Stroppel, Andrew Mathas,
Monica Vazirani and Weiqiang Wang for some helpful discussions.

\section{Vust's double centralizer property}

In this section we formulate Vust's double centralizer
theorem for the centralizer of the nilpotent matrix $e$.
This is the key technical ingredient needed to
prove the main double centralizer property 
formulated in Theorem A in the introduction.

\subsection{\boldmath The centralizer $\mathfrak{g}_e$}
We continue with the notation from the introduction.
So, $\mathfrak{g} = \mathfrak{gl}_N(\C)$ with natural module $V$,
and $\lambda = (0 \leq p_1 \leq \cdots \leq p_n)$ is a partition 
of $N$ with transpose partition
$(q_1 \geq \cdots \geq q_l \geq 0)$.
In \cite{BKrep}, we worked in slightly greater generality, replacing the Young diagram $\lambda$ with an arbitrary pyramid.
There seems to be little to be gained from that more general
setup (see \cite[Theorem 5.4]{MS} for a more precise statement), and many things later on become harder to prove,
so we will stick here to the left-justified case.
To help in translation, we 
define the shift matrix
$\sigma = (s_{i,j})_{1 \leq i,j \leq n}$ from \cite[$\S$3.1]{BKrep}
by $s_{i,j} = p_j - p_i$ if $i \leq j$, 
$s_{i,j} = 0$ if $i \geq j$.
Numbering the boxes of $\lambda$ as in the introduction, let 
$L(i)$ resp. $R(i)$ be the number of the box immediately to the left resp.
right of the $i$th box in the diagram $\lambda$, or 
$\varnothing$ if there is no such box.
Let
\begin{align*}
I &= \{1,2,\dots,N\},\\
J &= \{(i,j)\in I \times I\mid\col(i) \leq \col(j), R(j) = \varnothing\},\\
K &= \{(i,j,r)\mid 1 \leq i,j \leq n, s_{i,j} \leq r < p_j\}.
\end{align*}
The map $J \rightarrow K, (i,j) \mapsto (\row(i), \row(j),
\col(j)-\col(i))$ is a bijection.

The nilpotent matrix $e \in \mathfrak{g}$
maps $v_i \mapsto v_{L(i)}$, interpreting
$v_{\varnothing}$ as $0$.
Let $\mathfrak{g}_e$ be the  centralizer of $e$ in $\mathfrak{g}$ and
let $G_e$ be the corresponding algebraic group, i.e. the centralizer
of $e$ in the group $G = GL_N(\C)$.∂
It is well known, see e.g. \cite[Lemma 7.3]{BKshifted}, that
the Lie algebra $\mathfrak{g}_e$ has basis
$\{e_{i,j;r}\mid (i,j,r) \in K\}$ where
\begin{equation}\label{hunny}
e_{i,j;r} = \sum_{\substack{h,k \in I \\ \row(h) = i, \row(k) = j \\
\col(k) - \col(h) = r}} e_{h,k}.
\end{equation}
For example, if $\lambda$ is the partition whose Young diagram is 
displayed in the introduction and $n=3$, then
$$
e_{3,2;0} = e_{3,2}+e_{6,5}+e_{8,7},\qquad
e_{2,3;1} = e_{2,6}+e_{5,8}+e_{7,9},
\qquad
e_{1,3;3} = e_{1,9}.
$$
Note according to this definition that the notation
$e_{i,j;r}$ still makes sense even if $r \geq p_{j}$,
but then it is simply equal to zero.
We will often parametrize this basis instead by the
set
$J$: for $(i,j) \in J$, let
\begin{equation}\label{hunny2}
\xi_{i,j} 
= e_{i,j} + e_{L(i),L(j)} + \cdots + e_{L^k(i), L^k(j)}
\end{equation}
where $k = \col(i)-1$.
This is the same as the element $e_{\row(i), \row(j), \col(j)-\col(i)}$
from (\ref{hunny}), so in this alternate notation our basis for
$\mathfrak{g}_e$ is the set $\{\xi_{i,j}\mid (i,j) \in J\}$.

\subsection{\boldmath The graded Schur algebra associated to $e$}
Consider the
associative algebra $M_e$ of all $N \times N$ matrices
over $\C$ that commute with $e$.
Of course, $M_e$ is equal to $\mathfrak{g}_e$ as a vector space,
indeed, $\mathfrak{g}_e$ is the Lie algebra
equal to $M_e$ as a vector space with Lie bracket
being the commutator.
Multiplication in $M_e$ satisfies
\begin{equation*}
e_{i,j;r} e_{h,k;s} = \delta_{h,j} e_{i,k;r+s},
\end{equation*}
for $1 \leq i,j,h,k \leq n, r \geq s_{i,j}$ and $
s \geq s_{h,k}$.
Note also that the algebraic group $G_e$ is the
principal open subset of the
affine variety $M_e$ defined by the non-vanishing of determinant.
The coordinate algebra $\C[M_e]$ of $M_e$
is the polynomial algebra $\C[x_{i,j;r}\mid 1 \leq i,j \leq n,
s_{i,j} \leq r < p_j]$, where $x_{i,j;r}$ is the coordinate function
picking out the $e_{i,j;r}$-coefficient of an element of $M_e$
when expanded in terms of the above basis.
The natural monoid structure on $M_e$ induces a bialgebra
structure on $\C[M_e]$, with comultiplication
$\Delta:\C[M_e] \rightarrow \C[M_e] \otimes \C[M_e]$
and counit $\eps:\C[M_e] \rightarrow \C$ given 
on generators by
\begin{align*}
\Delta(x_{i,j;r}) &= \sum_{k=1}^n \sum_{\substack{s+t=r\\s\geq s_{i,k}\\ t \geq s_{k,j}}}
x_{i,k;s} \otimes x_{k,j;t},\\
\eps(x_{i,j;r}) &= \delta_{i,j} \delta_{r,0}
\end{align*}
for $1 \leq i,j \leq n$ and $s_{i,j} \leq r < p_{j}$.
The localization of $\C[M_e]$ at determinant is the 
coordinate algebra $\C[G_e]$
of the algebraic group
$G_e$, so it is a Hopf algebra.
We make $\C[M_e]$ into a graded algebra, 
$\C[M_e] = \bigoplus_{d \geq 0} \C[M_e]_d$,
by defining the degree
of each polynomial generator $x_{i,j;r}$ to be $1$.
Each $\C[M_e]_d$ is then a finite dimensional subcoalgebra of the
bialgebra $\C[M_e]$.
Hence the dual vector space
$\C[M_e]_d^*$ has the natural structure of a finite dimensional 
 algebra.

We pause to introduce some multi-index notation.
Let $I^d$ denote the set of all tuples $\bi = (i_1,\dots,i_d)$
with each $i_k \in I$. 
Let $J^d$ denote the set of all
pairs $(\bi,\bj)$ of tuples $\bi = (i_1,\dots,i_d)$
and $\bj = (j_1,\dots,j_d)$ with each $(i_k,j_k) \in J$.
Let $K^d$ denote the set of all triples
$(\bi,\bj,\brr)$ of tuples $\bi = (i_1,\dots,i_d), \bj = (j_1,\dots,j_d)$
and $\brr = (r_1,\dots,r_d)$ with each $(i_k,j_k,r_k) \in K$.
The symmetric group $S_d$ acts on the right on the set
$I^d$; we write
$\bi \sim \bj$ if the tuples 
$\bi$ and $\bj$ lie in the same $S_d$-orbit.
Similarly, $S_d$ acts diagonally on the right on the sets
$J^d$ and $K^d$. Choose sets of
orbit representatives $J^d / S_d$ and $K^d / S_d$.
For $\bi \in I^d$, let $\row(\bi) = (\row(i_1),\dots,\row(i_d))$
and define $\col(\bi)$ similarly. Then the map
\begin{equation}\label{bj}
J^d \rightarrow K^d, \qquad
(\bi,\bj) \mapsto (\row(\bi), \row(\bj), \col(\bj)-\col(\bi))
\end{equation}
is an $S_d$-equivariant bijection.
Finally, given $\bi \in I^d$ and $1 \leq j \leq d$, let 
\begin{equation}\label{ljd}
L_j(\bi) = \left\{
\begin{array}{ll}
(i_1,\dots,i_{j-1},L(i_j),i_{j+1},\dots,i_d)&\text{if $L(i_j) \neq \varnothing$,}\\
 \varnothing&\text{otherwise.}
\end{array}\right.
\end{equation}
Define $R_j(\bi)$ similarly.

Now we can write down an explicit basis for the algebra
$\C[M_e]_d^*$.
The monomials
$x_{\bi,\bj;\brr} = x_{i_1,j_1;r_1} \cdots x_{i_d,j_d;r_d}$
for $(\bi,\bj,\brr) \in K^d / S_d$
clearly form a basis for $\C[M_e]_d$.
Let
$\{\xi_{\bi,\bj;\brr}\mid (\bi,\bj,\brr) \in K^d / S_d\}$
denote the dual basis for
$\C[M_e]_d^*$.
Multiplication in the algebra
$\C[M_e]_d^*$ then satisfies
\begin{equation}\label{thealg}
\xi_{\bi,\bj;\brr} \xi_{\bh,\bk;\bs}
=
\sum_{(\bp,\bq,\bt) \in K^d / S_d}
a_{\bi,\bj,\bh,\bk,\bp,\bq;\brr,\bs,\bt}
\xi_{\bp,\bq;\bt},
\end{equation}
where
\begin{equation*}
a_{\bi,\bj,\bh,\bk,\bp,\bq;\brr,\bs,\bt}
=
\#
\left\{
(\bm,\bu,\bv)\:\Bigg|\:
\begin{array}{l}(\bp,\bm,\bu) \sim (\bi,\bj,\brr)\\
(\bm,\bq,\bv) \sim (\bh,\bk,\bs)\\
\bu+\bv = \bt\end{array}
\right\}.
\end{equation*}
All this is a straightforward generalization of \cite[ch.~2]{Green}.

Next we want to bring the universal enveloping algebra
$U(\mathfrak{g}_e)$ into the picture.
The augmentation ideal of the bialgebra $\C[M_e]$
is the kernel of $\eps:\C[M_e] \rightarrow \C$.
The algebra of distributions
\begin{equation*}
\operatorname{Dist}(M_e) = 
\left\{u \in \C[M_e]^*\mid u((\ker \eps)^m) = 0 \text{ for $m \gg 0$}\right\}
\end{equation*}
has a natural bialgebra structure dual to that on $\C[M_e]$ itself.
It is even a Hopf algebra because
there is an isomorphism
$U(\mathfrak{g}_e)
\stackrel{\sim}{\rightarrow} \operatorname{Dist}(M_e)$
induced by the Lie algebra homomorphism
$\mathfrak{g}_e
\rightarrow \operatorname{Dist}(M_e),
e_{i,j;r}\mapsto\eps \circ \frac{d}{dx}_{_{\!i,j;r}}$;
see \cite[I.7.10(1)]{Jantzen}.
We will always from now on identify
$U(\mathfrak{g}_e)$ with
$\operatorname{Dist}(M_e)$ in this way, i.e.
elements of $U(\mathfrak{g}_e)$
are functions on $\C[M_e]$.
Define an algebra homomorphism
\begin{equation*}
\pi_d: 
U(\mathfrak{g}_e)
\rightarrow \C[M_e]_d^*
\end{equation*}
by letting $\pi_d(u)(x) = u(x)$ for all
$u \in 
U(\mathfrak{g}_e)$ and
$x \in \C[M_e]_d$.

\begin{lem}\label{flu1}
$\pi_d$ is surjective.
\end{lem}

\begin{proof}
Suppose not. Then we can find $0 \neq x \in \C[M_e]_d$
such that $\pi_d(u)(x) = 0$ for all $u \in 
U(\mathfrak{g}_e)$.
Hence, for each $m \geq 1$, we have that
$u(x) = 0$ for all
$u \in \C[M_e]^*$ with $u((\ker \eps)^m) = 0$, i.e.
$x \in (\ker \eps)^m$.
But $\bigcap_{m \geq 1} (\ker \eps)^m = 0$.
\end{proof}

\subsection{The graded double centralizer property}
The natural $\mathfrak{g}$-module $V$ 
can be viewed as a right $\C[M_e]$-comodule with structure map 
\begin{equation}\label{xi}
V \rightarrow V \otimes \C[M_e],\qquad
v_j\mapsto
 \sum_{i \in I}
 v_i \otimes x_{\row(i), \row(j); \col(j)-\col(i)}
\end{equation}
where
the sum is over all $i$ satisfying $s_{i,j} \leq \col(j)-\col(i) < p_j$.
Since $\C[M_e]$ is a bialgebra
we get from this a
right $\C[M_e]$-comodule structure
on the tensor space $V^{\otimes d}$.
The image of the structure map of this comodule lies in 
$V^{\otimes d} \otimes \C[M_e]_d$, so
$V^{\otimes d}$ is actually a right $\C[M_e]_d$-comodule, hence a left
$\C[M_e]_d^*$-module.
Let
\begin{equation*}
\omega_d: \C[M_e]_d^* \rightarrow \End_{\C}(V^{\otimes d})
\end{equation*}
be the associated representation.
For $\bi,\bj \in I^d$, let
$v_\bi = v_{i_1}\otimes\cdots\otimes v_{i_d} \in V^{\otimes d}$
and $e_{\bi,\bj}
= e_{i_1,j_1}\otimes \cdots\otimes e_{i_d,j_d}
\in \End_\C(V^{\otimes d})$.
For $(\bi,\bj) \in J^d$, let
\begin{equation}\label{anew}
\xi_{\bi,\bj} = \sum_{(\bh,\bk) \sim (\bi,\bj)} \xi_{h_1, k_1} 
\otimes\cdots\otimes \xi_{h_d,k_d},
\end{equation}
recalling (\ref{hunny2}).
Using (\ref{xi}), one checks 
for $(\bi,\bj) \in J^d$ that
\begin{equation*}
\xi_{\bi,\bj}
=\omega_d(\xi_{\row(\bi),\row(\bj);\col(\bj)-\col(\bi)}).
\end{equation*}
Hence, recalling (\ref{bj}), $\omega_d$
maps the basis $\{\xi_{\bi,\bj;\brr}\mid (\bi,\bj,\brr) \in K^d / S_d\}$
of the algebra $\C[M_e]_d^*$ to the set
$\{\xi_{\bi,\bj}\mid (\bi,\bj) \in J^d / S_d\}$ in $\End_{\C}(V^{\otimes d})$.

\begin{lem}\label{flu2}
$\omega_d$ is injective.
\end{lem}

\begin{proof}
The vectors
$\{\xi_{\bi,\bj}\mid(\bi,\bj) \in J^d / S_d\}$ 
are linearly independent.
\end{proof}

Introduce the
twisted tensor product algebra 
$\C_l[x_1,\dots,x_d] \,{\scriptstyle{\rtimes\!\!\!\!\!\bigcirc}}\, \C S_d$
acting on the right on the space $V^{\otimes d}$
as in the introduction.

\begin{lem}\label{flu3}
Suppose we are given scalars $a_{\bi,\bj} \in \C$ 
such that the
endomorphism $\sum_{\bi,\bj \in I^d} a_{\bi,\bj} e_{\bi,\bj}$
of $V^{\otimes d}$ commutes with the action of
$\C_l[x_1,\dots,x_d] \,{\scriptstyle{\rtimes\!\!\!\!\!\bigcirc}}\, \C S_d$.
Then
 $a_{\bi,\bj}=a_{\bi\cdot w, \bj \cdot w}$
for all $w \in S_d$,
 $a_{\bi,\bj} = 0$  if
 $\col(i_k) > \col(j_k)$
for some $k$, and
$$
a_{\bi,\bj}=a_{R_k(\bi),R_k(\bj)}
$$
for all $k$
with
$R_k(\bj) \neq \varnothing$ (interpreting
the right hand side as $0$ in case $R_k(\bi) = \varnothing$).
\end{lem}

\begin{proof}
The image of $v_{\bj}$ under the given endomorphism
is $\sum_{\bi \in I^d} a_{\bi,\bj} v_{\bi}$. Applying $w \in S_d$, we get
$\sum_{\bi \in I^d} a_{\bi,\bj} 
v_{\bi \cdot w}$. This must equal the image of
$v_{\bj\cdot w}$ under our endomorphism, i.e.
$\sum_{\bi \in I^d} a_{\bi\cdot w,\bj\cdot w} v_{\bi\cdot w}$.
Equating coefficients gives that
$a_{\bi,\bj} = a_{\bi\cdot w, \bj \cdot w}$.
A similar argument using $x_k$ in place of $w$
gives that
$a_{R_k(\bi),\bj} = a_{\bi,L_k(\bj)}$,
interpreting the left resp. right hand side as zero if 
$R_k(\bi) = \varnothing$
resp. $L_k(\bj) = \varnothing$.
The remaining statements follow easily from this.
\end{proof}

Finally let $\phi_d:U(\mathfrak{g}_e) \rightarrow
\End_{\C}(V^{\otimes d})$ be the representation induced by the natural
action of the Lie algebra $\mathfrak{g}_e$
on the tensor space.
By the definitions,
$\phi_d$ is precisely the composite map
$\omega_d \circ \pi_d$, so Lemmas~\ref{flu1} and \ref{flu2} show that 
its image is isomorphic to $\C[M_e]_d^*$.
Also let $\psi_d:\C_l[x_1,\dots,x_d] \,{\scriptstyle{\rtimes\!\!\!\!\!\bigcirc}}\, \C S_d\rightarrow \End_{\C}(V^{\otimes d})^{\op}$ be the homomorphism arising from the
right action of $\C_l[x_1,\dots,x_d] \,{\scriptstyle{\rtimes\!\!\!\!\!\bigcirc}}\, \C S_d$ on $V^{\otimes d}$. Thus, we have defined maps
\begin{equation}\label{grgrpre}
U(\mathfrak{g}_e) \stackrel{\phi_d}{\longrightarrow}
\End_{\C}(V^{\otimes d}) \stackrel{\psi_d}{\longleftarrow}
\C_l[x_1,\dots,x_d]
 \,{\scriptstyle{\rtimes\!\!\!\!\!\bigcirc}}\, \C S_d.
\end{equation}
The $\Z$-grading on $\mathfrak{g}$ from the 
introduction extends to a grading
$U(\mathfrak{g}) = \bigoplus_{r \in \Z} U(\mathfrak{g})_r$
on its universal enveloping algebra, and $U(\mathfrak{g}_e)$
is a graded subalgebra.
Make $V$ into a graded module
by declaring
that each $v_i$ is of degree $(l- \col(i))$.
There are induced gradings
on $V^{\otimes d}$ and its endomorphism algebra
$\End_{\C}(V^{\otimes d})$,
so that the map $\phi_d$ is then a homomorphism of graded algebras.
Also define a grading on $\C_l[x_1,\dots,x_d]
 \,{\scriptstyle{\rtimes\!\!\!\!\!\bigcirc}}\, \C S_d$
by declaring that each $x_i$ is of degree $1$ and each $w \in S_d$ is of degree
$0$.
The map $\psi_d$ is then a homomorphism of graded algebras too.

\begin{thm}\label{grgr}
The maps $\phi_d$ and $\psi_d$ satisfy the double centralizer property, i.e.
\begin{align*}
\phi_d(U(\mathfrak{g}_e))=&\End_{\C_l[x_1,\dots,x_d] 
\,{\scriptscriptstyle{\rtimes\!\!\!\!\!\bigcirc}}\, \C S_d}(V^{\otimes d}),\\
&
\End_{\mathfrak{g}_e}(V^{\otimes d})^{\operatorname{op}} = \psi_d(\C_l[x_1,\dots,x_d] \,{\scriptstyle{\rtimes\!\!\!\!\!\bigcirc}}\, \C S_d).
\end{align*}
Moreover, if at least $d$ parts of $\lambda$ are equal to $l$,
then the map $\psi_d$ is injective.
\end{thm}

\begin{proof}
The second equality follows by a theorem of Vust; see \cite[$\S$6]{KP}.
For the first equality, note by
 Lemma~\ref{flu3} that
any element of
$\End_{\C_l[x_1,\dots,x_d] \,{\scriptscriptstyle{\rtimes\!\!\!\!\!\bigcirc}}\, \C S_d}(V^{\otimes d})$ is a linear combination of 
the elements $\xi_{\bi,\bj}$
from (\ref{anew}).
These belong to the image of $\phi_d$
by Lemma~\ref{flu1}.
Finally, assume that at least $d$ parts of $\lambda$ are equal to $l$.
Take any $\bi \in I^d$ such that $i_1,\dots,i_d$
are all distinct and $\col(i_1)=\cdots=\col(i_d) = l$.
Then, for $0 \leq r_1,\dots,r_d < l$ and $w \in S_d$, we have that
$v_\bi x_1^{r_1} \cdots x_d^{r_d} w
=
v_\bj$
where $\bj = L_1^{r_1} \circ \cdots \circ L_d^{r_d}(\bi) \cdot w$.
These vectors are obviously linearly independent, hence the vector
$v_\bi \in V^{\otimes d}$ generates a copy of the regular
$\C_l[x_1,\dots,x_d] \,{\scriptstyle{\rtimes\!\!\!\!\!\bigcirc}}\, \C S_d$ 
module.
This implies that $\psi_d$ is injective.
\end{proof}

\begin{rem}\label{drem}\rm
Theorem~\ref{grgr} is true over an arbitrary 
commutative ground ring $R$ (instead
of just over the field $\mathbb{C}$) 
providing one replaces the universal enveloping algebra
$U(\mathfrak{g}_e)$ everywhere with the algebra of distributions
of the group scheme $G_e$ over $R$.
To obtain this generalization, one first treats the case of algebraically
closed fields of positive characteristic,
using the arguments above together with a result of
Donkin \cite[Theorem 2.3c]{D} in place of \cite[$\S$6]{KP}.
From there, the result can be lifted to
$\mathbb{Z}$, hence to any other commutative ground ring.
\end{rem}

\section{Higher level Schur algebras}\label{hlsa}

In this section 
we develop the filtered deformation 
of Vust's double centralizer property.
This leads naturally to the 
definition of certain
finite dimensional quotients $W_d(\lambda)$
of finite $W$-algebras, which we call
higher level Schur algebras.

\subsection{\boldmath The finite $W$-algebra as a deformation of $U(\mathfrak{g}_e)$}
Recall that 
$$
\mathfrak{m} = \bigoplus_{r < 0}
\mathfrak{g}_r,
\quad
\mathfrak{h} = \mathfrak{g}_0,
\quad 
\mathfrak{p} = \bigoplus_{r \geq 0} \mathfrak{g}_r.
$$ 
Define the
finite $W$-algebra $W(\lambda)$ 
associated to the nilpotent matrix $e$
as in the introduction; see also \cite[$\S$3.2]{BKrep}.
As an algebra, $W(\lambda)$ is generated by 
certain elements 
$$
\left\{T_{i,j}^{(r+1)}\:\Big|\:(i,j,r) \in K\right\}
$$
of $U(\mathfrak{p})$ described explicitly  in \cite[$\S$3.7]{BKrep}.
The grading on $\mathfrak{p}$
extends to give a grading $U(\mathfrak{p}) = \bigoplus_{r \geq 0}
U(\mathfrak{p})_r$ on its universal enveloping algebra.
In general, $W(\lambda)$ is not a graded subalgebra, but the
grading on $U(\mathfrak{p})$ 
does at least induce a filtration
$\operatorname{F}_0 W(\lambda) \subseteq \operatorname{F}_1 W(\lambda) 
\subseteq \cdots$
on $W(\lambda)$ with
$$
\operatorname{F}_r W(\lambda) = W(\lambda) 
\cap \bigoplus_{s=0}^r U(\mathfrak{p})_s.
$$ 
The associated graded algebra $\gr W(\lambda)$
is canonically identified with a graded subalgebra of
$U(\mathfrak{p})$.
The first statement of the following lemma is a special case of 
general result
about finite $W$-algebras due
to Premet \cite[Proposition 2.1]{premet2}; see also \cite[Theorem 3.8]{BGK}.

\begin{lem}\label{ag}
We have that
$\gr W(\lambda) = U(\mathfrak{g}_e)$.
In fact, for $(i,j,r) \in K$,
we have that $T_{i,j}^{(r+1)}\in\mathrm{F}_r W(\lambda)$ and
$\gr_r T_{i,j}^{(r+1)} = (-1) e_{i,j;r}$, where $e_{i,j;r}$ is the element 
of $U(\mathfrak{g}_e)$
from (\ref{hunny}).
\end{lem}

\begin{proof}
The fact that $\gr_r T_{i,j}^{(r+1)} = (-1)^r e_{i,j;r}$
follows by a direct calculation using \cite[Lemma 3.4]{BKrep}.
By \cite[Lemma 3.6]{BKrep}, ordered monomials in the elements
$\{T_{i,j}^{(r+1)}\mid(i,j,r)\in K\}$
form a basis for $W(\lambda)$.
Combining this with the PBW theorem for $U(\mathfrak{g}_e)$,
we deduce that $\gr W(\lambda) = U(\mathfrak{g}_e)$.
\end{proof}

\subsection{Filtered deformation of tensor space}
For the remainder of the article, $\bc = (c_1,\dots,c_l)\in\C^l$ will denote
a fixed choice of {origin}; 
in the introduction we discussed only the 
most interesting situation  
$\bc = \biz$, i.e. when $c_1=\cdots=c_l = 0$.
Let $\eta_{\bc}: U(\mathfrak{p}) \rightarrow
U(\mathfrak{p})$ denote the automorphism
mapping $e_{i,j} \mapsto e_{i,j} + \delta_{i,j} c_{\col(i)}$
for each $e_{i,j} \in \mathfrak{p}$.
Let $V^{\otimes d}_{\bc}$ denote the graded $U(\mathfrak{p})$-module
equal as a graded vector space to $V^{\otimes d}$, but with action 
obtained by twisting the natural action by the automorphism
$\eta_{\bc}$, i.e.
$u \cdot v = \eta_{\bc}(u) v$.
Also let $\C_{\bc} = \C 1_\bc$ denote the one dimensional $\mathfrak{p}$-module
on which each $e_{i,j} \in \mathfrak{p}$ acts
by multiplication by $\delta_{i,j} c_{\col(i)}$; in particular,
$\C_\biz$ is the trivial module $\C$.
Obviously, we can identify
\begin{equation*}
\C_{\bc} \otimes V^{\otimes d} = V^{\otimes d}_{\bc}
\end{equation*}
so that $1_\bc \otimes v =v$ for every $v \in V^{\otimes d}$.
By restriction, we can view 
$V^{\otimes d}_{\bc}$ also as a graded $\mathfrak{g}_e$-module.
The resulting representation 
\begin{equation}\label{grepre2}
\phi_{d,\bc}: U(\mathfrak{g}_e)
\rightarrow \End_{\C}(V^{\otimes d}_{\bc})
\end{equation}
is just the composition $\phi_d \circ \eta_{\bc}$
where $\phi_d$ is the map from (\ref{grgrpre}).
Because the automorphism $\eta_{\bc}$ of $U(\mathfrak{p})$
leaves the subalgebra $U(\mathfrak{g}_e)$ invariant, the image of
$\phi_{d,\bc}$ is the same as the image of $\phi_d$ itself.
So the statement of Theorem~\ref{grgr} remains true if we
replace $\phi_d$ with $\phi_{d,\bc}$.

Of course, we can also consider $\C_{\bc}$ and $V^{\otimes d}_{\bc}$ as
$W(\lambda)$-modules by restriction.
Let
\begin{equation}\label{cdd1}
\Phi_{d,\bc}: W(\lambda) \rightarrow \End_{\C}(V^{\otimes d}_{\bc})
\end{equation}
be the associated representation.
View the graded algebra $\End_{\C}(V^{\otimes d}_\bc)$ 
as a filtered algebra with degree $r$ filtered piece defined
from
$$
\operatorname{F}_r \End_{\C}(V^{\otimes d}_\bc)
= \bigoplus_{s \leq r} \End_{\C}(V^{\otimes d}_\bc)_s.
$$
The homomorphism $\Phi_{d,\bc}$ is then a homomorphism of filtered
algebras. Moreover,
identifying the associated graded algebra 
$\gr \End_{\C}(V^{\otimes d}_{\bc})$ 
simply with $\End_{\C}(V^{\otimes d}_\bc)$ in the obvious way,
Lemma~\ref{ag} shows that the associated graded map 
$$\gr \Phi_{d,\bc}:\gr W(\lambda) \rightarrow
\End_{\C}(V^{\otimes d}_{\bc})$$
coincides with the map
$\phi_{d,\bc}$.
We stress that the automorphism
$\eta_{\bc}$ of $U(\mathfrak{p})$
does not in general leave the subalgebra $W(\lambda)$
invariant (unless $c_1=\cdots=c_l$), so unlike $\phi_{d,\bc}$
the image
of the map $\Phi_{d,\bc}$ definitely does depend on the
choice of $\bc$.

\subsection{Skryabin's equivalence of categories}
Recall the $(U(\mathfrak{g}), W(\lambda))$-bimodule
$Q_\chi = U(\mathfrak{g}) / U(\mathfrak{g}) I_\chi$
and the automorphism $\eta$
from the introduction.
We write $1_\chi$ for the coset of $1$ in $Q_\chi$.
Let $\mathcal{C}(\lambda)$ denote the category of all 
$\mathfrak{g}$-modules on which $(x - \chi(x))$ acts locally nilpotently
for all $x \in \mathfrak{m}$. Then, by Skryabin's theorem \cite{Skry},
the functor
$$
Q_\chi \otimes_{W(\lambda)} ?: W(\lambda)\text{-mod} \rightarrow \mathcal{C}(\lambda)$$
is an equivalence of categories;
see also \cite[$\S$8.1]{BKrep}.
Conversely,
given $M \in \mathcal{C}(\lambda)$, the subspace
\begin{equation}
\Wh(M) = \{v \in M\mid xv = \chi(x)v\text{ for all }x \in \mathfrak{m}\}
\end{equation}
is naturally a $W(\lambda)$-module with action
$u \cdot v = \eta(u) v$ for $u \in W(\lambda), v \in \Wh(M)$.
This defines a functor
$\Wh: \mathcal{C}(\lambda) \rightarrow W(\lambda)\text{-mod}$
which gives
an inverse equivalence to $Q_\chi \otimes_{W(\lambda)} ?$.

Given any $M \in \mathcal{C}(\lambda)$ and any finite dimensional
$\mathfrak{g}$-module $X$, the tensor product $M\otimes X$
belongs to $\mathcal{C}(\lambda)$.
This defines an exact functor 
$$
? \otimes X: \mathcal{C}(\lambda) \rightarrow \mathcal{C}(\lambda).
$$
Transporting through Skryabin's equivalence of categories, we obtain
a corresponding exact functor
\begin{equation*}
? \circledast X: W(\lambda)\text{-mod} \rightarrow
W(\lambda)\text{-mod},
\quad
M \mapsto \Wh((Q_\chi \otimes_{W(\lambda)} M) \otimes X);
\end{equation*}
see \cite[$\S$8.2]{BKrep}.

Given another finite dimensional $\mathfrak{g}$-module $Y$
and any $W(\lambda)$-module $M$,
there is a natural associativity isomorphism
\cite[(8.8)]{BKrep}:
\begin{equation}\label{a}
a_{M,X,Y}: (M \circledast X) \circledast Y \stackrel{\sim}{\longrightarrow} M \circledast (X \otimes Y).
\end{equation}

Suppose that $M$ is actually a $\mathfrak{p}$-module, and
view $M$ and $M \otimes X$ as $W(\lambda)$-modules by restricting
from $U(\mathfrak{p})$ to $W(\lambda)$.
Consider the map
\begin{equation}\label{muhat}
\hat\mu_{M,X}:(Q_\chi \otimes_{W(\lambda)} M) \otimes X
\rightarrow M \otimes X
\end{equation}
such that
$(\eta(u) 1_\chi \otimes m) \otimes x
\mapsto um \otimes x$
for $u \in U(\mathfrak{p})$, $m \in M$ and $x \in X$.
The restriction of this map to the subspace
$\Wh((Q_\chi \otimes_{W(\lambda)} M) \otimes X)$
obviously defines a {homomorphism}
of $W(\lambda)$-modules
\begin{equation}\label{mu}
\mu_{M, X}: 
M \circledast X \stackrel{\sim}{\longrightarrow}
M \otimes X.
\end{equation}
Much less obviously,
this homomorphism is actually
an {\em isomorphism};
see \cite[Corollary 8.2]{BKrep}.

We are going to apply these things to 
the $d$th power $(? \circledast V)^d$ of the functor
$? \circledast V$; we will often denote this simply by 
$(? \circledast V^{\circledast d})$.
To start with, take the module $M$ to be
the one dimensional $\mathfrak{p}$-module $\C_{\bc}$.
Composing $d$ isomorphisms built from  (\ref{a}) and (\ref{mu})
in any order that makes sense, we obtain an isomorphism
\begin{equation}\label{mud}
\mu_d: \C_\bc \circledast V^{\circledast d}
\stackrel{\sim}{\longrightarrow}
\C_\bc \otimes V^{\otimes d} = V_\bc^{\otimes d}
\end{equation}
of $W(\lambda)$-modules.
The associativity ``pentagons'' from \cite[(8.9)--(8.10)]{BKrep}
show that this definition is independent of the precise 
order chosen.
For example, 
$\mu_0:\C_\bc \rightarrow \C_\bc$ is the identity map and
$\mu_d$ is the composite
\begin{equation}\label{m1}
\begin{CD}
(\C_\bc \circledast V^{\circledast (d-1)}) \circledast V
&@>\mu_{d-1} \circledast\id_{V}>>
&V_\bc^{\otimes (d-1)} \circledast V
&@>\mu_{V_\bc^{\otimes (d-1)},V}>>
&V_\bc^{\otimes d}
\end{CD}
\end{equation}
for $d \geq 1$.
Alternatively, 
$\mu_d$ is the composite
\begin{equation}\label{m2}
\begin{CD}
\C_\bc \circledast V^{\circledast d}
&@>a_d>>
&\C_\bc \circledast V^{\otimes d}
&@>\mu_{\C_\bc, V^{\otimes d}}>>
&V_\bc^{\otimes d}
\end{CD}
\end{equation}
where $a_0:\C_\bc \rightarrow \C_\bc$ is the identity map and
 $a_d$ is the composite
$$
\begin{CD}
(\C_\bc \circledast V^{\circledast (d-1)}) \circledast V
&@>a_{d-1} \circledast \id_V>>
&(\C_\bc \circledast V^{\otimes (d-1)}) \circledast V
&@>a_{\C_\bc, V^{\otimes(d-1)}, V}>>
&\C_\bc \circledast V^{\otimes d}
\end{CD}
$$
for $d \geq 1$.
The following technical lemma
will be needed to compute the isomorphism
$\mu_{\C_\bc, V^{\otimes d}}$ appearing in (\ref{m2})
explicitly later on.
Recall $I^d$ denotes the set of multi-indexes
$\bi = (i_1,\dots,i_d) \in I^d$
and $v_{\bi} = v_{i_1}\otimes\cdots\otimes v_{i_d}$.

\begin{lem}\label{xij}
For all $\bi,\bj \in I^d$,
there exist elements
$x_{\bi,\bj} \in U(\mathfrak{p})$ such that
\begin{itemize}
\item[(i)]
$[e_{i,j}, \eta(x_{\bi,\bj})] + \sum_{\bk} \eta(x_{\bk,\bj}) 
\in U(\mathfrak{g}) I_\chi$
for each $e_{i,j} \in \mathfrak{m}$,
where the sum is over all $\bk \in I^d$ obtained from 
$\bi$ by replacing an entry equal to $i$ by $j$;
\item[(ii)]
$x_{\bi,\bj}$ acts on the module $\C_{\bc}$ as the scalar $\delta_{\bi,\bj}$.
\end{itemize}
The inverse image
of $v_{\bj} \in V_\bc^{\otimes d}$ under the isomorphism
$\mu_{\C_\bc, V^{\otimes d}}$ 
is equal to
$$
\sum_{\bi \in I^d} (\eta(x_{\bi,\bj}) 1_\chi \otimes 1_\bc) \otimes v_{\bi}
\in \C_\bc \circledast V^{\otimes d} \subseteq (Q_\chi \otimes_{W(\lambda)} \C_\bc) \otimes V^{\otimes d}
$$
for any such choice of elements $x_{\bi,\bj}$.
\end{lem}

\begin{proof}
The existence of such elements
follows from \cite[Theorem 8.1]{BKrep};
one needs to
choose
the projection $\mathrm{p}$ there
in a similar way to \cite[(8.4)]{BKrep}
so that all $u \in U(\mathfrak{p})$ with
$\mathrm{p}(u \cdot 1_\chi) = 0$ act as zero on $\C_{\bc}$.
The second statement then follows
from the explicit description of the map
$\mu_{\C_\bc, V^{\otimes d}}$ given in \cite[Corollary 8.2]{BKrep}.
\end{proof}

\subsection{Action of the degenerate affine Hecke algebra}
Now we bring the degenerate affine Hecke algebra $H_d$ into the picture
following the approach of Chuang and Rouquier
\cite[$\S$7.4]{CR}.
Consider the functor $? \otimes V$
on $\mathfrak{g}$-modules.
Define an endomorphism
$x:? \otimes V \rightarrow ? \otimes V$ by letting
\begin{equation}\label{xd}
x_M:M \otimes V \rightarrow M \otimes V
\end{equation}
denote the endomorphism
arising from left multiplication by
$\Omega = \sum_{i,j=1}^N e_{i,j} \otimes e_{j,i}$,
for each $\mathfrak{g}$-module $M$.
Define an endomorphism
$s:(? \otimes V) \otimes V \rightarrow (? \otimes V) \otimes V$
by letting
\begin{equation}\label{sd}
s_M:M \otimes V \otimes V \rightarrow M \otimes V \otimes V
\end{equation}
denote the endomorphism 
arising from left multiplication by $1 \otimes \Omega$, i.e.
$m \otimes v \otimes v' \mapsto
m \otimes v' \otimes v$.
More generally, for any $d \geq 1$, 
there are endomorphisms
$x_1,\dots,x_d$, $s_1,\dots,s_{d-1}$
of the functor $(? \otimes V^{\otimes d})= (? \otimes V)^d$
defined by
\begin{equation}\label{sk1}
x_i = \bid^{d-i} x \bid^{i-1},
\qquad
s_j = \bid^{d-j-1} s \bid^{j-1}.
\end{equation}
By \cite[Lemma 7.21]{CR}, these endomorphisms
induce  a well-defined right action of 
the degenerate affine Hecke
algebra $H_d$ on 
$M \otimes V^{\otimes d}$ for any $\mathfrak{g}$-module $M$. This is precisely the action
from \cite[$\S$2.2]{AS} described at the beginning of the introduction.
Transporting $x$ and $s$ through Skryabin's equivalence of categories,
we get endomorphisms
also denoted $x$ and $s$ of the functors
$? \circledast V$ and $(? \circledast V) \circledast V$;
see \cite[(8.13)-(8.14)]{BKrep}.
Once again, we set
\begin{equation}\label{sk2}
x_i = \bid^{d-i} x \bid^{i-1},
\qquad
s_j = \bid^{d-j-1} s \bid^{j-1},
\end{equation}
to get endomorphisms $x_1,\dots,x_d$, $s_1,\dots,s_{d-1}$ of 
the functor $(? \circledast V^{\circledast d})=(? \circledast V)^d$
for all $d \geq 1$.
By \cite[(8.18)--(8.23)]{BKrep},
these induce a well-defined right action of the
degenerate affine Hecke algebra $H_d$ on
$M \circledast V^{\circledast d}$,
for any $W(\lambda)$-module $M$.

Returning to the special case
$M = \C_{\bc}$ once again, we 
have made $\C_\bc \circledast V^{\circledast d}$
into a $(W(\lambda), H_d)$-bimodule.
Using the isomorphisms
$\C_\bc \circledast V^{\circledast d} \cong
\C_\bc \circledast V^{\otimes d}
\cong V_\bc^{\otimes d}$ from (\ref{m2}),
we lift the $H_d$-module structure on
$\C_\bc \circledast V^{\circledast d}$
first to 
$\C_\bc \circledast V^{\otimes d}$
then to $V_\bc^{\otimes d}$,
to make each of these spaces into 
$(W(\lambda), H_d)$-bimodules too.
Let us describe the resulting actions of the generators
of $H_d$ on
$\C_\bc \circledast V^{\otimes d}$
and on $V_\bc^{\otimes d}$ more explicitly.
For $\C_\bc \circledast V^{\otimes d}$, 
this is already explained in
\cite[$\S$8.3]{BKrep}:
the conclusion is that
$x_i$ and $s_j$ act on
$$
\C_\bc \circledast V^{\otimes d} \subseteq
(Q_\chi \otimes_{W(\lambda)} \C_\bc) \otimes V^{\otimes d}
$$
in the same way as the endomorphisms defined by multiplication by
$\sum_{h=1}^i \Omega^{[h,i+1]}$
and $\Omega^{[j+1,j+2]}$, respectively.
Here, $\Omega^{[r,s]}$ denotes $\sum_{i,j=1}^N 
1^{\otimes(r-1)} \otimes
e_{i,j} \otimes 1^{(s-r-1)} \otimes e_{j,i} \otimes 1^{(d+1-s)}$,
 treating $Q_\chi \otimes_{W(\lambda)} \C_\bc$
as the first tensor position and the copies of $V$
as positions $2,3,\dots,(d+1)$.
To describe the action of $H_d$ on
$V_\bc^{\otimes d}$, each $s_i$ acts by permuting the $i$th
and $(i+1)$th tensor positions, as follows from the naturality
of 
$\mu_{\C_\bc, V^{\otimes d}}: \C_\bc \circledast V^{\otimes d}
\rightarrow V_\bc^{\otimes d}$.
The action of each $x_j$ is more subtle, and is described 
by the following lemma; the formula for the action of $x_1$
on $V^{\otimes d}$ from the introduction is a special case of this.

\begin{lem}\label{action}
For $\bi \in I^d$
 and $1 \leq j \leq d$, 
the action of $x_j$ on $V_{\bc}^{\otimes d}$ satisfies
$$
v_{\bi} x_j
= 
v_{L_j(\bi)} 
+ (c_{\col(i_j)}+q_{\col(i_j)}-n) v_{\bi}
+
\sum_{\substack{1 \leq k < j \\ \col(i_k) \geq \col(i_j)}}
v_{\bi\cdot(k\:j)}
-
\sum_{\substack{j < k \leq d \\ \col(i_k) < \col(i_j)}}
v_{\bi\cdot(j\:k)}.
$$
(Recalling (\ref{ljd}),  the term $v_{L_j(\bi)}$ should be 
interpreted as $0$ in
case $L_j(\bi) = \varnothing$.)
\end{lem}

\begin{proof}
Note 
that $x_{j+1} = s_j x_{j} s_j + s_j$.
Using this and induction on $j$, one reduces just to checking the given formula
in the case $j=1$. 
Let $\mu$ denote the map $\hat\mu_{\C_c,V^{\otimes d}}$
from (\ref{muhat}).
Choosing elements $x_{\bi,\bj} \in U(\mathfrak{p})$ as in
Lemma~\ref{xij}, we have that
\begin{align*}
v_{\bi} x_1 &= 
\mu\bigg(
\Omega^{[1,2]} \sum_{\bj \in I^d} (\eta(x_{\bj,\bi}) 1_\chi \otimes 1_\bc) \otimes v_{\bj}
\bigg)
=
\sum_{\bj,\bk} \mu((e_{j_1,k_1} \eta(x_{\bj,\bi}) 1_\chi \otimes 1_\bc) 
\otimes v_{\bk})
\end{align*}
where the last sum is over $\bj,\bk \in I^d$
with $j_2=k_2,\cdots,j_d = k_d$.
Suppose first that $\col(j_1) \leq \col(k_1)$.
By Lemma~\ref{xij}(ii),
$\mu((e_{j_1,k_1} \eta(x_{\bj,\bi}) 1_\chi \otimes 1_\bc) 
\otimes v_{\bk}) = 0$ unless $\bi=\bj=\bk$, in which case
$$
\mu((e_{j_1,k_1} \eta(x_{\bj,\bi}) 1_\chi \otimes 1) 
\otimes v_{\bk})
=(c_{\col(i_1)}+q_{\col(i_1)}+\cdots+q_l - n) v_{\bi}.
$$
Now consider the terms with $\col(j_1) > \col(k_1)$.
By Lemma~\ref{xij}(i), we have that
$e_{j_1,k_1}
 \eta(x_{\bj,\bi}) 1_\chi
=
\eta(x_{\bj,\bi}) e_{j_1,k_1} 1_\chi - \sum_{\bh} \eta(x_{\bh,\bi}) 1_\chi$
where the sum is over all $\bh \in I^d$ obtained from $\bj$ by replacing 
an entry equal to $j_1$ by $k_1$.
Since $\chi(e_{j_1,k_1})$ is zero unless $k_1 = L(j_1)$ when it is simply
equal to $1$, and using Lemma \ref{xij}(ii),
it follows that these terms contribute
$v_{L_1(\bi)}
-
\sum_{\bj,\bk} v_{\bk}$
where the sum is over all $\bj,\bk \in I^d$ such that
$\col(j_1) > \col(k_1)$,
$j_2 = k_2,\dots,j_d = k_d$, 
and $\bi$ is obtained from $\bj$ by replacing an entry equal to $j_1$ by $k_1$.
This simplifies further to give
$$
v_{L_1(\bi)}
- (q_{\col(i_1)+1} + \cdots + q_l)v_{\bi}
- \sum_{\substack{1 < k \leq d \\ \col(i_k) < \col(i_1)}} v_{\bi \cdot (1 \, k)}
$$
which added to our earlier term give the conclusion.
\end{proof}

\subsection{Degenerate cyclotomic Hecke algebras}
Now we pass from the degenerate affine Hecke algebra
to its cyclotomic quotients.

\begin{lem}\label{minimal}
For $d \geq 1$, the minimal polynomial of the
endomorphism
of $V_\bc^{\otimes d}$ defined by the action of
$x_1$ is $\prod_{i=1}^l (x - (c_i+q_i-n))$.
\end{lem}

\begin{proof}
In the special case $d=1$,
this is clear from Lemma~\ref{action}.
By the definition of $x_1$, it
follows that for any $d \geq 1$,
the minimal polynomial certainly divides
$\prod_{i=1}^l (x- (c_i+q_i-n))$.
Finally, by Lemma~\ref{action}, 
the endomorphism of $V_\bc^{\otimes d}$ defined by 
$x_1$ 
belongs to $\operatorname{F}_1 \End_{\C}(V_\bc^{\otimes d})$,
and the associated graded endomorphism of $V_\bc^{\otimes d}$
is equal to $e \otimes 1^{\otimes(d-1)}$.
This clearly has minimal polynomial exactly equal to 
$x^l$, since $l$ is the size of the
largest Jordan block of $e$.
This implies that the minimal polynomial 
of the filtered endomorphism
cannot be of degree smaller than $l$, completing the proof.
\end{proof}

Let $\Lambda = \sum_{i=1}^l \Lambda_{c_i+q_i-n}$, an element
of the free abelian group 
generated by the symbols
$\{\Lambda_a\mid a \in \C\}$. 
The corresponding degenerate cyclotomic Hecke algebra
is the quotient
$H_d(\Lambda)$ of $H_d$ by the two-sided ideal 
generated
by $\prod_{i=1}^l (x_1-(c_i+q_i-n))$.
In view of Lemma~\ref{minimal}, the right action of $H_d$ on 
$V_\bc^{\otimes d}$ factors through the quotient $H_d(\Lambda)$, and we obtain
a homomorphism
\begin{equation}\label{psid}
\Psi_d: H_d(\Lambda) \rightarrow \End_{\C}(V_\bc^{\otimes d})^{\op}.
\end{equation}
Define a filtration 
$\operatorname{F}_0 H_d(\Lambda) \subseteq \operatorname{F}_1 H_d(\Lambda)
\subseteq \cdots$ by declaring that
$\operatorname{F}_r H_d(\Lambda)$ is the span of all
$x_1^{i_1} \cdots x_d^{i_d} w$ for $i_1,\dots,i_d \geq 0$ and $w \in S_d$
with $i_1+\cdots+i_d \leq r$.
Recalling the graded algebra
$\C_l[x_1,\dots,x_d]
 \,{\scriptstyle{\rtimes\!\!\!\!\!\bigcirc}}\, \C S_d$
from the previous section,
the relations imply that
there is a well-defined surjective homomorphism of graded algebras
$$
\zeta_d:
\C_l[x_1,\dots,x_d]
 \,{\scriptstyle{\rtimes\!\!\!\!\!\bigcirc}}\, \C S_d
\twoheadrightarrow
\gr H_d(\Lambda)
$$
such that $x_i \mapsto \gr_1 x_i$ for each $i$ and
$s_j \mapsto \gr_0 s_j$ for each $j$.
Moreover, by Lemma~\ref{action}, 
the map $\Psi_d$ from (\ref{psid}) 
is filtered and 
\begin{equation}\label{bag}
(\gr \Psi_d) \circ \zeta_d = \psi_d.
\end{equation}

\begin{lem}\label{ze}
The map $\zeta_d:
\C_l[x_1,\dots,x_d]
 \,{\scriptstyle{\rtimes\!\!\!\!\!\bigcirc}}\, \C S_d
\rightarrow
\gr H_d(\Lambda)$  is an isomorphism of graded algebras.
If in addition
at least $d$ parts of $\lambda$ are equal to $l$
then the map $\gr \Psi_d$, hence also the map $\Psi_d$ itself,
is injective.
\end{lem}

\begin{proof}
We can add extra parts equal to $l$ to the partition $\lambda$
if necessary to assume
that at least $d$ parts of $\lambda$ are equal to $l$ without affecting
$\Lambda$ or the map $\zeta_d$.
Under this assumption, the last sentence of Theorem~\ref{grgr}
asserts that the map $\psi_d$
is injective. 
Hence by (\ref{bag}) 
the maps $\zeta_d$ and $\gr \Psi_d$ are injective too.
\end{proof}

\subsection{The filtered double centralizer property}
To prove the main result of the section we just need one more general lemma.

\begin{lem}\label{thetrick}
Let 
$\Phi:B \rightarrow A$ and $\Psi:C \rightarrow A$ be homomorphisms of
filtered algebras such that
$\Phi(B) \subseteq Z_A(\Psi(C))$, where
$Z_A(\Psi(C))$ denotes the centralizer of $\Psi(C)$ in $A$.
View the subalgebras $\Phi(B)$, $\Psi(C)$ and $Z_A(\Psi(C))$
of $A$ as filtered algebras with filtrations 
induced by the one on $A$, so that the associated graded
algebras are naturally subalgebras of $\gr A$.
Then
$$
(\gr \Phi)(\gr B) \subseteq \gr \Phi(B)
\subseteq \gr Z_A(\Psi(C))
\subseteq Z_{\gr A} (\gr \Psi(C))
\subseteq Z_{\gr A}( (\gr \Psi)(\gr C)).
$$
\end{lem}

\begin{proof}
Exercise.
\end{proof}

We are going to apply this to the maps from
(\ref{cdd1}) and (\ref{psid}):
\begin{equation}\label{co2}
W(\lambda) \stackrel{\Phi_{d,\bc}}{\longrightarrow}
\End_{\C}(V_\bc^{\otimes d})
\stackrel{\Psi_d}{\longleftarrow}
H_d(\Lambda).
\end{equation}
As we have already said,
Lemma~\ref{ag} 
means that the associated
graded map $\gr \Phi_{d,\bc}$ is identified
with the map $\phi_{d,\bc}$ from (\ref{grepre2}).
Also Lemma~\ref{ze} and (\ref{bag}) identify
 $\gr \Psi_d$ with $\psi_d$.
Hence, taking
$A = \End_{\C}(V_\bc^{\otimes d})$,
Theorem~\ref{grgr} establishes that
$(\gr \Phi)(\gr B) = Z_{\gr A}((\gr \Psi)(\gr C))$
if $(B,C, \Phi,\Psi) = 
(W(\lambda), H_d(\Lambda)^{\operatorname{op}}, \Phi_{d,\bc}, \Psi_d)$
or if 
$(B,C,\Phi,\Psi) = (H_d(\Lambda)^{\operatorname{op}}, W(\lambda), \Psi_d, \Phi_{d,\bc})$.
The following theorem now follows by Lemma~\ref{thetrick}.

\begin{thm}\label{dcthm}
The maps $\Phi_{d,\bc}$ and $\Psi_d$ satisfy the double centralizer
property, i.e.
\begin{align*}
\Phi_{d,\bc}(W(\lambda))=&\End_{H_d(\Lambda)}(V_\bc^{\otimes d}),\\
&\End_{W(\lambda)}(V_\bc^{\otimes d})^{\operatorname{op}} 
= \Psi_d(H_d(\Lambda)).
\end{align*}
Moreover, at the level of associated graded algebras,
$\gr \Phi_{d,\bc}(W(\lambda))
= \phi_{d,\bc}(U(\mathfrak{g}_e))$
and
$\gr \Psi_d(H_d(\Lambda))
= \psi_d(\C_l[x_1,\dots,x_d]
 \,{\scriptstyle{\rtimes\!\!\!\!\!\bigcirc}}\, \C S_d)$.
\end{thm}

\subsection{Basis theorem for higher level Schur algebras}
Define
 the {\em higher level Schur algebra}
$W_{d}(\lambda,\bc)$
to be the image of the homomorphism $\Phi_{d,\bc} : W(\lambda)
\rightarrow \Hom_{\C}(V_\bc^{\otimes d})$.
Also let $H_{d}(\lambda,\bc)$ denote the image
of 
$\Psi_d:H_d(\Lambda) \rightarrow \Hom_{\C}(V_\bc^{\otimes d})^{\operatorname{op}}$.
In the most important case $\bc = \biz$,
we denote
$W_{d}(\lambda,\bc)$ and $H_{d}(\lambda,\bc)$
simply by $W_d(\lambda)$ and $H_d(\lambda)$.
Both $W_{d}(\lambda,\bc)$ and $H_{d}(\lambda,\bc)$
are filtered algebras with filtrations induced by the filtration on
$\End_{\C}(V_\bc^{\otimes d})$.
The last part of Theorem~\ref{dcthm}
shows that
\begin{align*}
\gr W_d(\lambda,\bc) &= \phi_{d,\bc}(U(\mathfrak{g}_e))
=\phi_d(U(\mathfrak{g}_e)),\\
\gr H_d(\lambda,\bc) &= \psi_d(\C_l[x_1,\dots,x_d] \,{\scriptstyle{\rtimes\!\!\!\!\!\bigcirc}}\, \C S_d).
\end{align*}
In particular,
recalling Lemmas~\ref{flu1} and \ref{flu2},
the algebra $W_d(\lambda,\bc)$ is a filtered deformation of the algebra
$\C[M_e]_d^*$, so
\begin{equation}
\dim W_d(\lambda,\bc) = \binom{\dim \mathfrak{g}_e + d-1}{d}.
\end{equation}
We want to construct an explicit basis
for $W_d(\lambda,\bc)$ lifting the basis 
$\{\xi_{\bi,\bj}\mid (\bi,\bj) \in J^d / S_d\}$
for the associated graded algebra from (\ref{anew}).

\begin{lem}\label{rec}
Suppose we are given scalars $a_{\bi,\bj} \in \C$ 
such that the endomorphism $\sum_{\bi,\bj \in I^d} a_{\bi,\bj} e_{\bi,\bj}$ of $V_\bc^{\otimes d}$ commutes with the action of
$H_d(\Lambda)$.
Then
 $a_{\bi,\bj}=a_{\bi\cdot w, \bj \cdot w}$
for all $w \in S_d$,
 $a_{\bi,\bj} = 0$  if
 $\col(i_k) > \col(j_k)$
for some $k$, and
\begin{multline*}
a_{\bi,\bj} = a_{R_k(\bi), R_k(\bj)}
+ (c_{\col(i_k)}+q_{\col(i_k)} - c_{\col(j_k)+1}-q_{\col(j_k)+1}) a_{\bi,R_k(\bj)}\\
+ \sum_{\substack{h \neq k \\ \col(i_h) \leq \col(i_k)\\ \col(j_h) \leq \col(j_k)}} a_{\bi \cdot (h\,k), R_k(\bj)}
- \sum_{\substack{h \neq k \\ \col(i_h) > \col(i_k)\\ \col(j_h) > \col(j_k)}} a_{\bi \cdot (h\,k), R_k(\bj)}
\end{multline*}
for all $k$ with
$R_k(\bj) \neq \varnothing$ (interpreting
the first term on the right hand side as $0$ in case $R_k(\bi) = \varnothing$).
\end{lem}

\begin{proof}
The fact that $a_{\bi,\bj} = a_{\bi\cdot w, \bj \cdot w}$ for all
$w \in S_d$ is proved  as in Lemma~\ref{flu3}.
By Theorem~\ref{dcthm} we know that 
all endomorphisms commuting with $H_d(\Lambda)$ belong to $W_d(\lambda,\bc)$.
By the definition of the latter algebra,
they are therefore contained in
the image of the algebra $U(\mathfrak{p})$ under its representation
on $V_\bc^{\otimes d}$. This proves that $a_{\bi,\bj} = 0$
if $\col(i_k) > \col(j_k)$ for some $k$. For the final formula, 
a calculation like in the proof of Lemma~\ref{flu3}, using
Lemma~\ref{action} to compute the action of $x_k$, gives that
\begin{multline*}
a_{R_k(\bi), \bj} + (c_{\col(i_k)}+q_{\col(i_k)}-n)a_{\bi,\bj} + \!\!\!\!\sum_{\substack{1 \leq h < k \\ \col(i_h) \leq \col(i_k)}} a_{\bi\cdot(h\,k), \bj} - \!\!\!\!\sum_{\substack{k < h \leq d \\ \col(i_h) > \col(i_k)}}
a_{\bi\cdot(h\,k), \bj}\\
=
a_{\bi,L_k(\bj)} + (c_{\col(j_k)}+q_{\col(j_k)}-n) a_{\bi,\bj}
+ \!\!\!\sum_{\substack{1 \leq h < k \\ \col(j_h) \geq \col(j_k)}} a_{\bi,\bj\cdot (h\,k)} - \!\!\!\sum_{\substack{k < h \leq d \\ \col(j_h) < \col(j_k)}} a_{\bi,\bj\cdot (h\,k)},
\end{multline*}
interpreting the first terms on either side as zero if
$R_k(\bi)$ or $L_k(\bj)$ is $\varnothing$.
Now replace $\bj$ by $R_k(\bj)$ and simplify.
\end{proof}

\begin{thm}\label{basis}
For each $(\bi,\bj) \in J^d$, there exists a unique element
$$
\Xi_{\bi,\bj} = \sum_{\bh,\bk \in I^d} a_{\bh,\bk}  e_{\bh,\bk}
\in  W_d(\lambda,\bc)
$$
such that the coefficients $a_{\bh,\bk}$ satisfy $a_{\bi,\bj} = 1$ and
$a_{\bh,\bk} = 0$ for all $(\bh,\bk) \in J^d$ with
$(\bh,\bk) \not\sim (\bi,\bj)$.
Moreover:
\begin{itemize}
\item[\rm(i)] 
The elements $\{\Xi_{\bi,\bj}\mid (\bi,\bj) \in J^d / S_d\}$
form a basis for $W_d(\lambda,\bc)$.
\item[\rm(ii)]
Letting $r = \col(j_1)-\col(i_1)+\cdots+\col(j_d)-\col(i_d)$,
$\Xi_{\bi,\bj}$ belongs to $\mathrm{F}_r W_d(\lambda)$ and
$\gr_r \Xi_{\bi,\bj} = \xi_{\bi,\bj}$.
\end{itemize}
\end{thm}

\begin{proof}
Lemma~\ref{rec} and Theorem~\ref{dcthm}
show that 
any element $\sum_{\bi,\bj \in I^d} a_{\bi,\bj}  e_{\bi,\bj}$
of $W_d(\lambda,\bc)$ is 
uniquely determined just by the values of the coefficients $a_{\bi,\bj}$
for $(\bi,\bj) \in J^d / S_d$.
Since $\dim W_d(\lambda,\bc) = |J^d / S_d|$, this
implies that for each $(\bi,\bj) \in J^d$ there exists a unique
element $\Xi_{\bi,\bj} \in W_d(\lambda,\bc)$
with coefficients
$a_{\bh,\bk}$ satisfying $a_{\bi,\bj} = 1$ and
$a_{\bh,\bk} = 0$ for all $(\bh,\bk) \in J^d$ with
$(\bh,\bk) \not\sim (\bi,\bj)$.
Repeating this argument with $W_d(\lambda,\bc)$ replaced by $\operatorname{F}_r
W_d(\lambda,\bc)$ implies that $\Xi_{\bi,\bj}$ actually lies in
$\operatorname{F}_r W_d(\lambda,\bc)$, where
$r = \col(j_1)-\col(i_1)+\cdots+\col(j_d)-\col(i_d)$.
Hence the coefficients $a_{\bh,\bk}$ of $\Xi_{\bi,\bj}$ are zero
unless $\col(k_1)-\col(h_1)+\cdots+\col(k_d)-\col(h_d) \leq r$.
Moreover, by Lemma~\ref{rec}, $a_{\bh,\bk} = 0$ if
$\col(h_j) > \col(k_j)$ for some $j$, and
$a_{\bh,\bk} = a_{\bh \cdot w, \bk \cdot w}$ for all $w \in S_d$.
Now the recurrence relation in Lemma~\ref{rec} gives a recipe to compute 
all remaining coefficients from the known 
coefficients $a_{\bh,\bk}$ for $(\bh,\bk) \in J^d$, proceeding by
downward induction on $\col(k_1)+\cdots+\col(k_d)$.
This establishes in particular that 
the coefficients of top degree $r$, i.e. the $a_{\bh,\bk}$'s
with $\col(k_1)-\col(h_1) + \cdots + \col(k_d)-\col(h_d)=r$,
satisfy the same recurrence as the coefficients of $\xi_{\bi,\bj}$
from Lemma~\ref{flu3}. Hence $\gr_r \Xi_{\bi,\bj} = \xi_{\bi,\bj}$.
Finally, the fact that the elements $\{\Xi_{\bi,\bj}\mid (\bi,\bj) \in J^d / S_d\}$
form a basis for $W_d(\lambda,\bc)$ follows because the
elements $\{\xi_{\bi,\bj}\mid (\bi,\bj) \in J^d / S_d\}$ form a basis
for $\gr W_d(\lambda,\bc)$.
\end{proof}

\begin{rem}\label{Rform}\rm
Assume throughout this remark that $\bc = \biz$.
In that case, it is clear from the proof of Theorem~\ref{basis}
that all the coefficients of $\Xi_{\bi,\bj}$ are actually integers.
Using this and Remark~\ref{drem}, we can extend
some aspects of Theorem~\ref{dcthm}
to an arbitrary commutative ground ring $R$.
Let $H_d(\Lambda)_{\Z}$ be the subring of $H_d(\Lambda)$
generated by $x_1,\dots,x_d$ and $S_d$.
Set $H_d(\Lambda)_R = R \otimes_{\Z} H_d(\Lambda)_{\Z}$.
Let $V_{\Z}^{\otimes d}$ be the $\Z$-submodule
of $V^{\otimes d}$ spanned by
$\{v_{\bi}\mid \bi \in I^d\}$ 
and $V_R^{\otimes d} = R \otimes_{\Z} V_\Z^{\otimes d}$.
By Lemma~\ref{action},
the action of $H_d(\Lambda)$ restricts to a well-defined
action of $H_d(\Lambda)_{\Z}$ on $V^{\otimes d}_{\Z}$,
hence we also get a right action of
$H_d(\Lambda)_R$ on $V^{\otimes d}_R$ by
extending scalars. Let $H_d(\lambda)_R$ denote the image of
the resulting homomorphism
$H_d(\Lambda)_R \rightarrow \End_R(V^{\otimes d}_R)^{\op}$.
Let $W_d(\lambda)_{R} = \End_{H_d(\Lambda)_R}(V_{R}^{\otimes d})$.
It is then the case that
$$
W_d(\lambda)_R = R \otimes_{\Z} W_d(\lambda)_\Z,
\qquad
H_d(\lambda)_R = R \otimes_{\Z} H_d(\lambda)_\Z.
$$
Moreover, the double centralizer property
$\End_{W_d(\lambda)_R}(V_R^{\otimes d})^{\operatorname{op}}
=H_d(\lambda)_R$
still holds over $R$.
These statements are proved by an argument similar
to the above, but using
Remark~\ref{drem}
in place of Theorem~\ref{grgr}.
The first step is to 
observe that the elements
$\{\Xi_{\bi,\bj}\mid (\bi,\bj) \in J^d / S_d\}$
from Theorem~\ref{basis} make sense over any ground ring $R$,
and give a basis for $W_d(\lambda)_{R}$ as a free $R$-module.
\end{rem}

\section{\boldmath Parabolic category $\mathcal{O}$}

A basic theme from now on is that results
about the representation theory of $H_d(\Lambda)$ 
should be deduced from known results about 
category $\mathcal{O}$ for the Lie algebra
$\mathfrak{g}$ relative to the parabolic subalgebra $\mathfrak{p}$.
In this section we give a brief review of the latter
theory, taking our notation from \cite[ch.~4]{BKrep}.
We will assume for the remainder of the article
that the
origin
$\bc = (c_1,\dots,c_l) \in \C^l$ satisfies the condition
\begin{equation}\label{cond}
c_i - c_j \in \Z
\quad\Rightarrow\quad
c_i = c_j.
\end{equation}
There is no loss in generality in making this assumption,
because for every choice 
of $a_1,\dots,a_l \in \C$ it is possible to write
$\Lambda = \Lambda_{a_1}+\cdots+\Lambda_{a_l}$
as $\sum_{i=1}^l \Lambda_{c_i+q_i-n}$
for suitable $\lambda$ and $\bc$ satisfying (\ref{cond}).

\subsection{Combinatorics of tableaux}
By a {\em $\lambda$-tableau} we mean a filling of the boxes
of the diagram $\lambda$ by complex numbers.
We adopt the following notations.
\begin{itemize}
\item $\Tab(\lambda)$ denotes the set of all 
$\lambda$-tableaux.
\item $\Tab^d(\lambda)$ is the subset
of $\Tab(\lambda)$ consisting of the tableaux 
whose entries are non-negative integers summing to $d$.
\item
$\gamma(A) = (a_1,\dots,a_N) \in \C^N$
is the {\em column reading} of a $\lambda$-tableau $A$,
that is, the tuple 
defined so that $a_i$ is the entry in the $i$th box of $A$.
\item $\theta(A)$ denotes the {\em content} of $A$,
that is,
the element $\sum_{i=1}^N \gamma_{a_i}$ of the free abelian group
on generators $\{\gamma_a\mid a \in \C\}$.
\item $A_0$ denotes the {\em ground-state tableau}, that is,
the $\lambda$-tableau having all entries on its $i$th row equal
to $(1-i)$ (recall we number rows $1,\dots,n$ from top to bottom).
\item $A_{\bc}$ denotes
the tableau obtained from $A_0$ by adding
$c_i$ to all of the entries in the $i$th column  for each
$i=1,\dots,l$.
\item $\Col(\lambda)$ denotes the set of all {column
strict $\lambda$-tableaux}, where
a tableau is called {\em column strict}
if in every column its entries are strictly increasing
from bottom to top, working always with the partial order on $\C$
defined by $a \leq b$ if $(b-a) \in \N$.
\item
$\Col_{\bc}(\lambda)$ denotes the set of all $A \in \Col(\lambda)$
with the property that 
every entry in the $i$th column of $A$ belongs to $c_i  + \Z$ 
for each $i=1,\dots,l$.
\item
$\Std_{\bc}(\lambda)$ denotes the set of all 
{standard tableaux}
in $\Col_{\bc}(\lambda)$, where 
a column strict tableau $A$ is
{\em standard} if in every row its entries are non-decreasing from left to 
right; this means that if $x < y$ are entries on the same row of $A$ then
$x$ appears to the left of $y$.
\item $\Col_{\bc}^d(\lambda)$ 
denotes the set of all $A \in \Col_{\bc}(\lambda)$ 
such that
the entrywise difference
$(A - A_{\bc})$ lies in $\Tab^d(\lambda)$.
\item $\Std_\bc^d(\lambda)$
denotes
$\Std_\bc(\lambda) \cap \Col_{\bc}^d(\lambda)$.
\end{itemize}
Finally let $\geq$ denote the Bruhat ordering on the set $\Col(\lambda)$,
defined as in \cite[$\S$2]{qla} (where it is denoted $\geq'$).
Roughly speaking, this is generated by the basic move of swapping an entry in 
some column
with a smaller entry in some column further to the right,
then reordering columns in increasing order.
For example,
$$
\Diagram{4&6\cr 3&4&5\cr 1&1&2\cr}
\:>\: \Diagram{4&5\cr 3&4&6\cr 1&1&2\cr}
\:>\: \Diagram{4&5\cr 3&2&6\cr 1&1&4\cr}
\:>\: \Diagram{3&5\cr 2&4&6\cr 1&1&4\cr}\:\:.
$$
\vspace{0.5mm}
We note by \cite[Lemma 1]{qla}
that the restriction of this ordering to 
the subset $\Std_\bc(\lambda)$ coincides with the restriction of the
Bruhat ordering on  row standard tableaux from \cite[(4.1)]{BKrep}.
The first combinatorial lemma is the reason the 
assumption (\ref{cond}) is necessary.

\begin{lem}\label{silly}
Suppose that $A \in \Col_\bc^d(\lambda)$ and $B \in \Col_\bc(\lambda)$ satisfy
$\theta(A) = \theta(B)$. Then $B \in \Col_\bc^d(\lambda)$ too.
\end{lem}

\begin{proof}
Suppose that $B \notin \Col_{\bc}^d(\lambda)$.
Then, for some $i$ and $j$,
the entry in the $i$th row and $j$th column of $B$ is $< (c_j+1-i)$.
Hence as $B$ is column strict, the entry in the $n$th row and $j$th
column of $B$ is $< (c_j+1-n)$.
As $\theta(A) = \theta(B)$, this implies that some entry of $A$
is $< (c_j+1-n)$ too.
This contradicts the assumption that $A \in \Col_{\bc}^d(\lambda)$.
\end{proof}

\subsection{\boldmath Parabolic category $\mathcal O$}
Now let $\mathcal{O}(\lambda)$ denote the category of all
finitely generated $\mathfrak{g}$-modules
that are locally finite dimensional over $\mathfrak{p}$
and semisimple over $\mathfrak{h}$; recall here that $\mathfrak{h}$
is the standard Levi subalgebra
$\mathfrak{gl}_{q_1}(\C) \oplus \cdots \oplus
\mathfrak{gl}_{q_l}(\C)$ of 
$\mathfrak{g}$.
Modules in $\mathcal{O}(\lambda)$ are automatically
diagonalizable with respect to the standard 
Cartan subalgebra consisting of diagonal matrices in $\mathfrak{g}$.
We say that a vector $v$ is of {\em weight}
$(a_1,\dots,a_N) \in \C^N$ if $e_{i,i} v = a_i v$
for each $i=1,\dots,N$.
There is the usual dominance ordering on the set
$\C^N$ of weights, positive roots coming from
the standard Borel subalgebra
of upper triangular matrices in $\mathfrak{g}$.
Also let $?^\#$ denote the duality 
on $\mathcal{O}(\lambda)$ mapping a module to the direct
sum of the duals of its weight spaces with 
$\mathfrak{g}$-action defined via the antiautomorphism $x \mapsto x^t$
(matrix transposition).

To a column strict tableau
$A \in \Col(\lambda)$ we associate the following
modules:
\begin{itemize}
\item $N(A)$ is the usual parabolic Verma module in $\mathcal O(\lambda)$
of highest weight
$(a_1,a_2+1,\dots,a_N+N-1)$ where
$\gamma(A) = (a_1,\dots,a_N)$ is the column reading of $A$.
\item $K(A)$ is the unique irreducible quotient of $N(A)$.
\item $P(A)$ is the projective cover of $K(A)$
in $\mathcal O(\lambda)$.
\item $V(A)$ denotes the finite dimensional irreducible
$\mathfrak{h}$-module of highest weight
$\gamma(A - A_0)$, viewed as a $\mathfrak{p}$-module
via the natural projection $\mathfrak{p}\twoheadrightarrow\mathfrak{h}$.
\item $V(A_0)$ is the trivial $\mathfrak{p}$-module $\C$.
\item $V(A_{\bc})$ is the one dimensional $\mathfrak{p}$-module $\C_{\bc}$
on which $e_{i,j}$ acts as $\delta_{i,j}c_{\col(i)}$
as in the previous section.
\end{itemize}
We point out that
\begin{equation}\label{nadef}
N(A) = U(\mathfrak{g}) \otimes_{U(\mathfrak{p})} (\C_{-\rho} \otimes V(A))
\end{equation}
where $\C_{-\rho}$ is the one dimensional $\mathfrak{p}$-module from the introduction; cf. \cite[(4.38)]{BKrep}.
Moreover,
$\{K(A)\mid A \in \Col(\lambda)\}$ is a full set of 
irreducible modules in $\mathcal{O}(\lambda)$.

A standard fact is that $K(A)$ and $K(B)$ have the same central
character if and only if $\theta(A) = \theta(B)$; indeed,
the composition multiplicity $[N(A):K(B)]$ is zero unless
$A \geq B$ in the Bruhat ordering. 
Recall also that BGG reciprocity in this setting
asserts that 
$P(A)$ has a parabolic Verma flag, and the multiplicity
of $N(B)$ in any parabolic Verma flag of $P(A)$ is given by
\begin{equation}\label{bgg}
(P(A):N(B)) = 
\dim \Hom_{\mathfrak{g}}(P(A), N(B)^\#) = [N(B):K(A)].
\end{equation}
One consequence of this is that
\begin{equation}\label{Ext}
\operatorname{Ext}^1_{\mathcal O(\lambda)}(N(A), N(B)) = 0
\end{equation}
for all $A, B \in \Col(\lambda)$ with $A \not< B$.

From now on, we restrict our attention to the
Serre subcategory $\mathcal O_{\bc}(\lambda)$ of $\mathcal O(\lambda)$
generated by the irreducible modules
$\{K(A)\mid  A \in \Col_{\bc}(\lambda)\}$.
This is a sum of certain blocks of $\mathcal O(\lambda)$, and
conversely every block of $\mathcal O(\lambda)$ is a block of
$\mathcal O_{\bc}(\lambda)$ for some $\bc$ chosen as in (\ref{cond}).
For $d \geq 0$, let $\mathcal O_\bc^d(\lambda)$ denote the Serre subcategory
of $\mathcal O_\bc(\lambda)$ generated by
the irreducible modules
$\{K(A)\mid 
A \in \Col_\bc^d(\lambda)\}$. 
Letting  $\bc - r\bid$ denote the tuple $(c_1-r,\dots,c_l - r)$
for each $r \in \N$,
it is obvious that
$\Col_\bc^d(\lambda) \subseteq \Col_{\bc-r\bid}^{d+rN}(\lambda)$ and
$\Col_\bc(\lambda) = \Col_{\bc - r \bid}(\lambda)$.
Hence, 
\begin{equation}\label{expafter}
\mathcal O_\bc^d(\lambda) 
\subseteq \mathcal O_{\bc - r \bid}^{d+rN}(\lambda),
\qquad
\mathcal O_\bc(\lambda) = \mathcal O_{\bc - r \bid}(\lambda).
\end{equation}
This is useful because, by Lemma~\ref{silly}, 
$\mathcal O_\bc^d(\lambda)$ is a sum of blocks
of $\mathcal O_\bc(\lambda)$, and conversely
every block of $\mathcal O_\bc(\lambda)$ is a block of
$\mathcal O_{\bc - r \bid}^d(\lambda)$ for some
$r,d \geq 0$.
In this way, questions about $\mathcal O(\lambda)$
can usually be reduced to questions about the subcategories
$\mathcal O_\bc^d(\lambda)$.
In particular, since $\Col_{\bc}^0(\lambda) = \{A_\bc\}$,
we get that 
$K(A_\bc)$ is the unique (up to isomorphism) irreducible module in 
its block, hence $P(A_\bc) = N(A_\bc) = K(A_\bc)$.
This module plays a special role; we denote it simply by 
$P_\bc$. 

\begin{lem}\label{allirr}
Let $A \in \Col(\lambda)$.
Then $K(A)$ is a composition factor of $P_\bc \otimes V^{\otimes d}$
if and only if $A \in \Col_\bc^d(\lambda)$.
\end{lem}

\begin{proof}
Recall that $\C_\bc \otimes V^{\otimes d} = V_{\bc}^{\otimes d}$.
By the classical Schur-Weyl duality,
$V_\bc^{\otimes d}$ decomposes as an $\mathfrak{h}$-module
as a direct sum of copies of the modules $V(A)$ for $A \in \Col_\bc^d(\lambda)$,
and each such $V(A)$ appears at least once.
Hence, $V_\bc^{\otimes d}$ has a filtration as a $\mathfrak{p}$-module
with sections of the form
$V(A)$ for $A \in \Col_\bc^d(\lambda)$, and each such
$V(A)$ appears at least once.
By the tensor identity,
$$
P_\bc \otimes V^{\otimes d}
\cong U(\mathfrak{g}) \otimes_{U(\mathfrak{p})} ( \C_{-\rho}\otimes 
V_\bc^{\otimes d}).
$$ 
Recalling (\ref{nadef}), these two statements imply that
$P_\bc \otimes V^{\otimes d}$ has a parabolic Verma flag with sections
of the form $N(A)$ for $A \in \Col_\bc^d(\lambda)$, and each such $N(A)$
appears at least once. This shows that 
every $K(A)$ for $A \in \Col_\bc^d(\lambda)$
is a composition factor of $P_\bc \otimes V^{\otimes d}$.
Conversely, if $K(A)$ is a composition factor of $P_\bc \otimes V^{\otimes d}$
for some $A \in \Col(\lambda)$,
it is also a composition factor of $N(B)$ for some $A \leq 
B \in \Col_\bc^d(\lambda)$.
Hence by Lemma~\ref{silly} 
we see that $A \in \Col_\bc^d(\lambda)$ too.
\end{proof}

\subsection{Special projective functors}
We also need the special projective functors
$f_i$ and $e_i$ on $\mathcal O(\lambda)$.
We refer the reader to \cite[$\S$4.4]{BKrep}
for a more detailed discussion of the signficance of these
functors; see also \cite{BKtf} and \cite[$\S$7.4]{CR}.
For a module $M$, we have by definition that
$f_i M$ resp. $e_i M$ is the generalized
$i$-eigenspace resp. generalized $-(N+i)$-eigenspace
of $\Omega$ acting on $M \otimes V$ resp. $M \otimes V^*$,
where $V^*$ is the dual of the natural $\mathfrak{g}$-module $V$.
In particular,
we have that
$$
M \otimes V = \bigoplus_{i \in \C} f_i M,
\qquad\qquad
M \otimes V^* = \bigoplus_{i \in \C} e_i M.
$$
An important point is that the functors $f_i$ and $e_i$ are both left
and right adjoint to each other, 
so they send projective
modules in $\mathcal O(\lambda)$ to projective modules. 
The following lemma, which is a well known
consequence of the tensor identity,
shows that $f_i$ and $e_i$ map objects in $\mathcal O_\bc(\lambda)$
to objects in $\mathcal O_{\bc}(\lambda)$, and 
$f_i$ maps objects in $\mathcal O_{\bc}^d(\lambda)$ to
objects in $\mathcal O_{\bc}^{d+1}(\lambda)$.

\begin{lem}\label{or}
For $A \in \Col_\bc(\lambda)$ and $i \in \C$, we have that
\begin{itemize}
\item[(i)] 
$f_i N(A)$ has a multiplicity-free parabolic Verma flag with sections
of the form $N(B)$, one for each $B \in \Col_\bc(\lambda)$ obtained from 
$A$ by replacing an entry equal to $i$ with $(i+1)$;
\item[(ii)] 
$e_i N(A)$ has a multiplicity-free parabolic Verma flag with sections
of the form $N(B)$, one for each $B \in \Col_\bc(\lambda)$ obtained from 
$A$ by replacing an entry equal to $(i+1)$ with $i$.
\end{itemize}
\end{lem}

\begin{cor}\label{v}
For any $A \in \Col_\bc(\lambda)$, there exist $i_1,\dots,i_r \in \C$ and
an irreducible parabolic Verma module
$N(B)$ for some
$B \in \Std_\bc(\lambda)$ such that $N(A)$ is isomorphic to a submodule of
$e_{i_1} \cdots e_{i_r} N(B)$.
\end{cor}

\begin{proof}
Let $B$ be the tableau obtained from $A$ by adding $k_j$ to all
entries in the $j$th column of $A$ for each $j=2,\dots,l$, 
where $k_2,\dots, k_l \geq 0$
are chosen so that the entries in column $j$ of $B$
are $\not\leq$ the entries in columns $1,\dots,j-1$ for each $j$.
It is then
automatic that $B \in \Std_\bc(\lambda)$, and $N(B)$ is irreducible
since $B$ is minimal in the Bruhat ordering.
Say the entries in the $j$th column of $A$ are
$a_{j,1} < \cdots < a_{j,q_j}$.
Apply the functors
\begin{multline*}
(e_{a_{j,q_j}} e_{a_{j,q_j}+1}\cdots e_{a_{j,q_j}+k_j-1})
\cdots 
(e_{a_{j,2}} e_{a_{j,2}+1}\cdots e_{a_{j,2}+k_j-1})\\(e_{a_{j,1}} e_{a_{j,1}+1}\cdots e_{a_{j,1}+k_j-1})
\end{multline*}
for $j=2$, then $j=3$, \dots, then $j=l$
to the parabolic Verma module $N(B)$.
Now use Lemma~\ref{or}, recalling by (\ref{Ext}) 
that parabolic Verma flags can be arranged in any order
refining the Bruhat order (largest at the bottom), 
to show that the result involves $N(A)$ as a submodule.
\end{proof}

The effect of the projective functors $f_i$ and $e_i$ on irreducible modules
is reflected by a crystal structure 
$(\Col(\lambda), \tilde e_i, \tilde f_i, \eps_i, \varphi_i, \theta)$
on the set $\Col(\lambda)$ 
in the general sense of Kashiwara \cite{Ka}.
This particular crystal 
is defined combinatorially in \cite[$\S$4.3]{BKrep};
there we only described the situation for $i \in \Z$ but it extends in
obvious fashion to arbitrary $i \in \C$.
In particular, \cite[Theorem 4.5]{BKrep}
asserts for $A \in \Col(\lambda)$
that $f_i  K(A)$ is non-zero if and only if
$\varphi_i(A) \neq 0$, in which case $f_i K(A)$ has irreducible
socle and cosocle isomorphic to $K(\tilde f_i A)$.
There is a similar statement about $e_i$.
The following lemma follows from the explicit
combinatorial definition of this crystal structure.

\begin{lem}\label{ccpt}
The subset $\Std_\bc(\lambda)$ of $\Col_\bc(\lambda)$ is 
the connected component of the crystal 
$(\Col(\lambda), \tilde e_i, \tilde f_i, \eps_i, \varphi_i, \theta)$
generated by the ground-state tableau $A_\bc$.
Moreover, for $d > 0$, 
every $A \in \Std_\bc^d(\lambda)$ is equal to
$\tilde f_i B$ for some $B \in \Std_\bc^{d-1}(\lambda)$ and $i \in \C$.
\end{lem}

\begin{cor}\label{q}
If $A \in \Std_\bc^d(\lambda)$ then 
$K(A)$ is a quotient of $P_\bc \otimes V^{\otimes d}$.
\end{cor}

\begin{proof}
Proceed by induction on $d$, the case $d=0$ being obvious.
For $d > 0$, Lemma~\ref{ccpt} shows that
$A = \tilde f_i B$ for some $B \in \Std_\bc^{d-1}(\lambda)$.
By induction, $K(B)$ is a quotient of
$P_\bc \otimes V^{\otimes (d-1)}$.
So $f_i K(B)$ is a quotient of $f_i (P_\bc \otimes V^{\otimes (d-1)})$,
which is itself a summand of $P_\bc \otimes V^{\otimes d}$.
Since $K(A)$ is a quotient of $f_i K(B)$
by \cite[Theorem 4.5]{BKrep}, this completes the proof.
\end{proof}

Let $[\mathcal O_\bc(\lambda)]$ denote the Grothendieck group of the
category $\mathcal O_\bc(\lambda)$.
It has (at least) two natural bases $\{[N(A)]\mid A \in \Col_\bc(\lambda)\}$
and $\{[K(A)]\mid A \in \Col_\bc(\lambda)\}$ corresponding to the
parabolic Verma modules and the irreducible modules.
The exact functors $f_i$ and $e_i$ induce operators
on $[\mathcal O_\bc(\lambda)]$ too.
The combinatorics of all these things is well 
understood by the Kazhdan-Lusztig
conjecture proved in \cite{BB, BrK}; 
see \cite[Theorem 4.5]{BKrep} for the application to this specific situation.
With that in mind, the next lemma
is a consequence
of \cite[Theorem 26]{qla}.
For a representation theoretic explanation, we note that $\{K(A)\mid A \in \Std_\bc(\lambda)\}$ are exactly the 
irreducible modules in $\mathcal{O}_{\bc}(\lambda)$ with maximal Gelfand-Kirillov dimension; see Irving's proof of Theorem~\ref{i} below.

\begin{lem}\label{layer}
The subgroup of the Grothendieck group $[\mathcal O_\bc(\lambda)]$
generated by the classes $\{[K(A)]\mid A \in \Col_\bc(\lambda) \setminus \Std_\bc(\lambda)\}$
is stable under the actions of $f_i$ and $e_i$. 
\end{lem}

\subsection{Self-dual projective modules}
Next, we record an important result originating in work of Irving \cite{I}.
Since our combinatorial setup is unconventional, we include the
proof, though this is essentially a translation of Irving and 
Shelton's proof from \cite[$\S$A.3]{I} into our language.
We learned this argument
from \cite[Theorem 5.1]{MS} which treats the case of a 
regular block.
In the statement of the theorem,
we have stuck with Irving's original terminology of
self-dual projective indecomposable modules;
note though that these modules are also the indecomposable 
prinjective
modules (as in \cite{KSX} or \cite{MS}) and the indecomposable projective tilting modules (as in \cite{St2}).

\begin{thm}\label{i} For $A \in \Col_\bc(\lambda)$, the following are equivalent:
\begin{itemize}
\item[(i)] $A \in \Std_\bc(\lambda)$;
\item[(ii)] $P(A)$ is self-dual, i.e. $P(A)^\# \cong P(A)$;
\item[(iii)] $K(A)$ is isomorphic to a submodule of
$N(B)$ for some $B \in \Col_\bc(\lambda)$.
\end{itemize}
\end{thm}

\begin{proof}
(i)$\Rightarrow$(ii).
Replacing $\bc$ by $\bc - r\bid = (c_1-r, \dots, c_l-r)$
for some $r \geq 0$ if necessary,
we may assume that $A \in \Col_{\bc}^d(\lambda)$
for some $d \geq 0$.
Note that $P_\bc \otimes V^{\otimes d}$ is a self-dual projective module.
By Corollary~\ref{q}, $P(A)$ is a summand of $P_\bc \otimes V^{\otimes d}$.
Hence $P(A)$ is itself self-dual; see \cite[Lemma 4.2]{I}.

(ii)$\Rightarrow$(iii).
If $P(A)$ is self-dual, then $K(A)$ is a submodule
of it. But $P(A)$ has a parabolic Verma flag, hence $K(A)$ is also 
a submodule of some $N(B)$.

(iii)$\Rightarrow$(i).
If $K(A)$ is a submodule of some parabolic Verma module, then
by Corollary~\ref{v} it is also a submodule of
$e_{i_1} \cdots e_{i_r} K(B)$ for some $B \in \Std_\bc(\lambda)$
and $i_1,\dots,i_r \in \C$. Hence by adjointness, the irreducible
module $K(B)$
appears in the cosocle of
$f_{i_r} \cdots f_{i_1} K(A)$.
This implies that $A \in \Std_\bc(\lambda)$ by Lemma~\ref{layer}.
\end{proof}

\begin{cor}\label{only}
If $K(A)$ is a quotient of $P_\bc \otimes V^{\otimes d}$
then $A \in \Std_\bc^d(\lambda)$.
\end{cor}

\begin{proof}
By Lemma~\ref{allirr}, we may assume that $A \in \Col_\bc(\lambda)$ and need
to prove that
$A \in \Std_\bc(\lambda)$.
Since $K(A)$ is a quotient of $P_\bc \otimes V^{\otimes d}$, 
$P(A)$ is a summand of it.
Since $P_\bc \otimes V^{\otimes d}$ is a self-dual projective module,
this implies that $P(A)$ is too.
Hence $A \in \Std_\bc(\lambda)$ by Theorem~\ref{i}.
\end{proof}

\subsection{Divided power modules}
Finally we introduce certain {\em divided power modules} in  $\mathcal{O}_\bc^d(\lambda)$.
For $1\leq j\leq l$, let $V_j$ denote the 
$\mathfrak{p}$-submodule of $V$ spanned by 
$\{v_i\mid \col(i) \leq j\}$.
Given a tableau $A \in \Tab^d(\lambda)$ with $\gamma(A) = (a_1,\dots,a_N)$, set
\begin{equation}\label{after}
Z^A (V) = Z^{a_1} \!\left(V_{\col(1)}\right) \otimes \cdots \otimes Z^{a_N} \!\left(V_{\col(N)}\right),
\end{equation}
where we write $Z^d (E)$ for the $d$th {\em divided power}
of a vector space $E$, 
i.e. the subspace of $E^{\otimes d}$ consisting of all symmetric tensors.
Clearly this is a $\mathfrak{p}$-submodule of the tensor space
$V^{\otimes d}$.
Let $Z^A_\bc(V)$ denote the $\mathfrak{p}$-submodule of $V^{\otimes d}_\bc$ 
obtained from $Z^A(V)$ by twisting the action by
the automorphism $\eta_\bc$.
By the tensor identity, the induced $\mathfrak{g}$-module
\begin{equation}\label{zac}
Z(A,\bc) = U(\mathfrak{g}) \otimes_{U(\mathfrak{p})} (\C_{-\rho} \otimes Z_\bc^A (V))
\end{equation}
is a submodule of $P_\bc \otimes V^{\otimes d}$, so it belongs to the
category $\mathcal O^d_\bc(\lambda)$.
As usual, in the most interesting case $\bc = \biz$, we denote this simply
by $Z(A)$.

To formulate next lemma, we need certain generalized Kostka numbers.
Take $A \in \Tab^d(\lambda)$ and $B \in \Col_\bc^d(\lambda)$.
Define an $l$-tuple of partitions 
$(\mu^{(1)}, \dots,\mu^{(l)})$ so that the parts of 
$\mu^{(j)}$ are the entries in the $j$th column of $(B - A_\bc)$. 
Letting $\gamma(A) = (a_1,\dots,a_N)$, define $K_{B,A}$ to be the number of 
tuples $(T^{(1)},\dots,T^{(l)})$ such that
\begin{itemize}
\item[(a)] for each $j=1,\dots,l$, $T^{(j)}$ is a standard $\mu^{(j)}$-tableau
with (not necessarily distinct) 
entries chosen from the set $\{i\in I\mid \col(i) \geq j\}$;
\item[(b)] for each $i \in I$, the total number of entries equal to 
$i$ in all of $T^{(1)}, \dots, T^{(l)}$ is equal to $a_i$.
\end{itemize}

\begin{lem}\label{s}
For any $A \in \Tab^d(\lambda)$, there is an isomorphism of
$\mathfrak{h}$-modules
$$
Z_\bc^A (V)
\cong \bigoplus_{B \in \Col_\bc^d(\lambda)} V(B)^{\oplus K_{B,A}}.
$$
Hence, 
$Z(A,\bc)$ has a parabolic Verma flag in which the
parabolic Verma module $N(B)$ appears with multiplicity
$K_{B,A}$ for each $B \in \Col_\bc^d(\lambda)$,
arranged in any order refining the Bruhat ordering on $\Col_\bc^d(\lambda)$
(most dominant at the bottom).
\end{lem}

\begin{proof}
The first statement is a consequence of the Littlewood-Richardson rule.
Hence, as a $\mathfrak{p}$-module,
$Z_\bc^A(V)$ has a filtration with sections
of the form $V(B)$ for $B \in \Col_\bc^d(\lambda)$,
each $V(B)$ appearing with multiplicity $K_{B,A}$.
Now apply the exact functor $U(\mathfrak{g}) \otimes_{U(\mathfrak{p})} (\C_{-\rho} \otimes ?)$
and use the definition (\ref{nadef})
and (\ref{Ext}).
\end{proof}

The goal in the remainder of the 
section is to prove by induction on level that
$Z(A,\bc)$ is a projective module in $\mathcal O^d_\bc(\lambda)$.
Let us set up some convenient notation for the proof.
Suppose $N = N' + N''$, let $\mathfrak{g}' = \mathfrak{gl}_{N'}(\C)$,
$\mathfrak{g}'' = \mathfrak{gl}_{N''}(\C)$, and identify $\mathfrak{g}'
\oplus \mathfrak{g}''$ with a standard Levi subalgebra of $\mathfrak{g}$
in the usual way. 
Given a $\mathfrak{g}'$-module $M'$ and a $\mathfrak{g}'$-module $M''$,
we will write $M' \boxtimes M''$ for the $\mathfrak{g}'\oplus \mathfrak{g}''$-module obtained from their (outer) tensor product.
Let $\mathfrak{q}$ be the standard (upper triangular) parabolic subalgebra of $\mathfrak{g}$ with Levi subalgebra $\mathfrak{g}'\oplus \mathfrak{g}''$ and 
nilradical $\mathfrak{r}$.
Let 
$\beta$ be the weight
$$
\beta = 
(\underbrace{0,\dots,0}_{N'\text{ times}},\underbrace{N',\dots,N'}_{N''\text{ times}})
$$
and 
$\C_\beta$ be the corresponding 
one dimensional $\mathfrak{g}'\oplus \mathfrak{g}''$-module.
Define a functor
\begin{align*}
R&:\mathfrak{g}' \oplus \mathfrak{g}''\text{-mod}
\rightarrow \mathfrak{g}\text{-mod},
&M &\mapsto U(\mathfrak{g})\otimes_{U(\mathfrak{q})} (\C_\beta \otimes M),\\
\intertext{where $\C_\beta\otimes M$ is viewed as a $\mathfrak{q}$-module via the natural projection $\mathfrak{q}\twoheadrightarrow \mathfrak{g}' \oplus \mathfrak{g}''$.
Thus, $R$ is the usual Harish-Chandra induction functor, but shifted by $\beta$ as that makes subsequent
notation tidier. Note the exact functor $R$ has a right adjoint}
T&:\mathfrak{g}\text{-mod}\rightarrow \mathfrak{g}'\oplus \mathfrak{g}''\text{-mod},
&M &\mapsto \C_{-\beta}\otimes M^{\mathfrak{r}}
\end{align*}
arising by taking $\mathfrak{r}$-invariants.

Now assume in addition that $l = l' + l''$ and that $\lambda'$
resp. $\lambda''$ is the partition of $N'$ resp. $N''$ 
whose Young diagram consists of the $l'$ leftmost resp. $l''$ rightmost columns of the Young diagram of $\lambda$. Let $\bc' = (c_1,\dots,c_{l'})$ and $\bc'' = (c_{l'+1},\dots,c_{l})$.
Given any $A \in \Col_\bc(\lambda)$ we let $A'$ resp. $A''$ 
denote the element of $\Col_{\bc'}(\lambda')$ resp. $\Col_{\bc''}(\lambda'')$ consisting of the 
$l'$ leftmost resp. $l''$ rightmost columns of the tableau $A$.
For any $d \geq 0$, this defines a bijection
$$
\Col_\bc^d(\lambda) \longrightarrow \dot{\bigcup_{d'+d''=d}} \Col_{\bc'}^{d'}(\lambda')
\times \Col_{\bc''}^{d''}(\lambda''),
\qquad
A \mapsto (A',A'').
$$
Consider parabolic category $\mathcal O$
for $\mathfrak{g}' \oplus \mathfrak{g}''$ relative to the standard parabolic
subalgebra having Levi subalgebra $\mathfrak{h}$.
Let $\mathcal{O}_{\bc'}^{d'}(\lambda') \boxtimes \mathcal O_{\bc''}^{d''}(\lambda'')$ denotes the Serre subcategory 
generated by the modules $\{K(A') \boxtimes K(A'')\mid A' \in \Col^{d'}_{\bc'}(\lambda'),
A'' \in \Col^{d''}_{\bc''}(\lambda'')\}$.
By Lemma~\ref{silly}, this is a sum of blocks of the full
parabolic category $\mathcal O$.
Moreover, the parabolic Verma modules in $\mathcal{O}_{\bc'}^{d'}(\lambda')
\boxtimes \mathcal{O}_{\bc''}^{d''}(\lambda'')$
are the modules $N(A') \boxtimes N(A'')$, and 
the projective cover
of $K(A') \boxtimes K(A'')$
in $\mathcal O_{\bc'}^{d'}(\lambda') \boxtimes \mathcal O_{\bc''}^{d''}(\lambda'')$
is $P(A')\boxtimes P(A'')$,
for each $A' \in \Col_{\bc'}^{d'}(\lambda'),
A'' \in \Col_{\bc''}^{d''}(\lambda'')$.
The following well known lemma implies in particular that the functors $R$ and $T$
restrict to give well-defined functors
\begin{equation}\label{RT}
\bigoplus_{d'+d''=d}^{\phantom{H}}
\mathcal O_{\bc'}^{d'}(\lambda') \boxtimes \mathcal O_{\bc''}^{d''}(\lambda'')
\:\:\phantom{\longrightarrow} \:\:
\mathcal O_{\bc}^d(\lambda).
\begin{picture}(0,0)
  \put(-45,4){\makebox(0,0){$\stackrel{\stackrel{\,\scriptstyle{R}_{\phantom{,}}}{\displaystyle\longrightarrow}}{\stackrel{\displaystyle\longleftarrow}{\scriptstyle T}}$}}
\end{picture}
\end{equation}

\begin{lem}\label{stad}
Let $A, B \in \Col_\bc^d(\lambda)$ be tableaux such that
$A',B' \in \Col_{\bc'}^{d'}(\lambda')$ and $A'',B'' \in \Col_{\bc''}^{d''}(\lambda'')$ for some $d'+d''=d$.
Then
\begin{itemize}
\item[(i)] $R(N(A') \boxtimes N(A'')) \cong N(A)$;
\item[(ii)] $T(N(A)^\#) \cong N(A')^\# \boxtimes N(A'')^\#$;
\item[(iii)] $T(K(A)) \cong K(A') \boxtimes K(A'')$;
\item[(iv)] $[N(A):K(B)] = [N(A'):K(B')][N(A''):K(B'')]$.
\end{itemize}
\end{lem}

\begin{cor}\label{pi}
Let $A \in \Col_\bc^d(\lambda)$ such that
$A' \in \Col_{\bc'}^{d'}(\lambda')$ and
$A'' \in \Col_{\bc''}^{d''}(\lambda'')$ for some
$d'+d''=d$.
\begin{itemize}
\item[(i)] If $A$ is minimal 
in the sense that every
$A \geq B \in \Col_\bc^d(\lambda)$ also has
$B' \in \Col_{\bc'}^{d'}(\lambda')$ and
$B'' \in \Col_{\bc''}^{d''}(\lambda'')$ then
$R(K(A') \boxtimes K(A'')) \cong K(A)$.
\item[(ii)] If $A$ is maximal in the sense that every
$A \leq B \in \Col_\bc^d(\lambda)$ also has
$B' \in \Col_{\bc'}^{d'}(\lambda')$ and
$B'' \in \Col_{\bc''}^{d''}(\lambda'')$ then
$R(P(A') \boxtimes P(A'')) \cong P(A)$.
\end{itemize}
\end{cor}

\begin{proof}
(i)
By Lemma~\ref{stad}(iv) and the minimality assumption on $A$, 
we have for every $B \leq A$ that
$[N(A):K(B)] = [N(A'):K(B')][N(A''):K(B'')]$. Hence in the Grothendieck group, we can write
$$
[N(A') \boxtimes N(A'')] = \sum_{B \leq A} [N(A):K(B)] [K(B') \boxtimes K(B'')].
$$
The functor $R$ is exact, so we get by Lemma~\ref{stad}(i) that
$$
[N(A)] = \sum_{B \leq A} [N(A):K(B)] [R(K(B') \boxtimes K(B''))]
= \sum_{B \leq A} [N(A):K(B)] [K(B)].
$$
Since $R(K(B') \boxtimes K(B''))$ certainly contains at least one composition factor isomorphic to $K(B)$, the last equality here 
implies that for all $B \leq A$ we must have that
$R(K(B') \boxtimes K(B'')) \cong K(B)$.

(ii) By adjointness of $R$ and $T$
and Lemma~\ref{stad}(iii), 
$R(P(A')\boxtimes P(A''))$ has irreducible cosocle isomorphic
to $K(A)$.
Hence there is a surjective homomorphism
$P(A) \twoheadrightarrow R(P(A') \boxtimes P(A''))$.
To prove that this is an isomorphism, observe 
by Lemma~\ref{stad}(i) that
$R(P(A') \boxtimes P(A''))$ has a parabolic Verma flag, so it 
suffices to check that
$$
(P(A):N(B)) = (R(P(A') \boxtimes P(A'')):
N(B))
$$
for all $B \in \Col^d_\bc(\lambda)$.
By (\ref{bgg}), the 
left hand side is equal to $[N(B):K(A)]$, which is zero unless
$A \leq B$.
Hence by the maximality assumption, the left hand side is zero
unless $B' \in \Col_{\bc'}^{d'}(\lambda')$
and $B'' \in \Col_{\bc''}^{d''}(\lambda'')$, in which case 
by Lemma~\ref{stad}(iv) it is equal to
$[N(B'):K(A')][N(B''):K(A'')]$.
On the other hand, by Lemma~\ref{stad}(ii), we have that
\begin{align*}
(R(P(A')\boxtimes P(A'')):&N(B)) =
\dim \Hom_{\mathfrak{g}}(R(P(A') \boxtimes P(A'')), N(B)^\#)\\
&=
\dim \Hom_{\mathfrak{g}'\oplus \mathfrak{g}''}
(P(A')\boxtimes P(A''), N(B')^\#\boxtimes N(B'')^\#)\\
&=
[N(B'):K(A')] [N(B''):K(A'')].
\end{align*}
This completes the proof.
\end{proof}

\begin{cor}\label{equivcat}
Assume that no entry of $\bc' = (c_1,\dots,c_{l'})$ lies in the same coset
of $\C$ modulo $\Z$ as an entry of $\bc'' = (c_{l'+1},\dots,c_l)$.
Then, the category $\mathcal O_{\bc}^d(\lambda)$ decomposes as
$$
\mathcal O_{\bc}^d(\lambda) = \bigoplus_{d'+d'' = d}
\mathcal O_{\bc}^{d',d''}(\lambda),
$$
where $\mathcal O_{\bc}^{d',d''}(\lambda)$
denote the Serre subcategory of $\mathcal O_{\bc}^d(\lambda)$
generated by the irreducible modules
$\{K(A)\mid A \in \Col^d_{\bc}(\lambda)\text{ such that }
A' \in \Col^{d'}_{\bc'}(\lambda'), A'' \in \Col^{d''}_{\bc''}(\lambda'')\}$.
Moreover, 
the functors $R$ and $T$ from (\ref{RT}) are mutually inverse
equivalences of categories. 
\end{cor}

\begin{proof}
Fix $d'+d''=d$.
Given $A \in \Col_{\bc}^d(\lambda)$
with $A' \in \Col_{\bc'}^{d'}(\lambda')$ and
$A'' \in \Col_{\bc''}^{d''}(\lambda'')$,
the assumption on $\bc$ ensures that any other
$B \in \Col_{\bc}^d(\lambda)$ with
$\theta(A) = \theta(B)$
also has
$B' \in \Col_{\bc'}^{d'}(\lambda')$ and
$B'' \in \Col_{\bc''}^{d''}(\lambda'')$.
Hence, $\mathcal O_{\bc}^{d',d''}(\lambda)$ is a sum of blocks
of $\mathcal O_{\bc}^{d}(\lambda)$, proving the first statement.
Moreover, by Lemma~\ref{stad}, the functors $R$ and $T$ restrict to well-defined functors
$$\mathcal O_{\bc'}^{d'}(\lambda') \boxtimes \mathcal O_{\bc''}^{d''}(\lambda'')
\:\:\phantom{\longrightarrow} \:\:
\mathcal O_{\bc}^{d',d''}(\lambda).
\begin{picture}(0,18)
  \put(-59.5,4.5){\makebox(0,0){$\stackrel{\stackrel{\,\scriptstyle{R}_{\phantom{,}}}{\displaystyle\longrightarrow}}{\stackrel{\displaystyle\longleftarrow}{\scriptstyle T}}$}}
\end{picture}
$$
\vspace{0.1mm}

\noindent
To prove the second statement, it suffices to show that these
restricted functors are mutually inverse equivalences of categories.

We claim to start with that the functor $T$ is exact.
Since it is automatically left exact, we just 
need to check given any epimorphism
$M \twoheadrightarrow Q$ in $\mathcal O_{\bc}^{d',d''}(\lambda)$
that every $\mathfrak{r}$-invariant in $Q$ lifts to an
$\mathfrak{r}$-invariant in $M$.
The proof of this easily reduces to the case that
$Q$ is irreducible. In that case, by Lemma~\ref{stad}(iii), it is enough 
to show that the 
highest weight vector $v \in Q$ lifts to an $\mathfrak{r}$-invariant in $M$.
To see that, let $\bar v \in M$ be any weight vector that is a pre-image of $v$, then observe
by weight considerations that $\bar v$ is automatically 
an $\mathfrak{r}$-invariant.

Now we can show that the unit 
and the counit
of any adjunction between $R$ and $T$ are both isomorphisms,
to complete the proof. Since both functors are exact, it suffices to 
check that the unit and counit give isomorphisms on every irreducible module,
which follows easily from Lemma~\ref{stad}(iii) and Corollary~\ref{pi}(i).
\end{proof}

Now we can prove the desired projectivity of the modules
$Z(A,\bc)$ from (\ref{zac}).

\begin{thm}\label{pgen}
Let $A \in \Tab^d(\lambda)$ and define $B$ to be the unique tableau in 
$\Col_\bc^d(\lambda)$ such that $(B - A_\bc)$ and $A$ have the same 
entries in each column.
Then, 
$$
Z(A,\bc) \cong
\bigoplus_{C \in \Col^d_\bc(\lambda)} n_{B,C} P(C)
$$
for $n_{B,C} \in \N$ with $n_{B,B} = 1$ and $n_{B,C} = 0$ unless
$\gamma(C) \geq \gamma(B)$ in the dominance ordering.
In particular, $\bigoplus_{A\in \Tab^d(\lambda)} Z(A,\bc)$ is a projective generator for the category $\mathcal O^d_\bc(\lambda)$.
\end{thm}

\begin{proof}
By Lemma~\ref{s}, $Z(A,\bc)$ has a parabolic Verma flag
involving $N(B)$ with multiplicity one and various other $N(C)$'s for
$C \in \Col^d_\bc(\lambda)$ with $\gamma(C) \geq \gamma(B)$.
This reduces the proof of the theorem simply to showing that 
$Z(A,\bc)$ is a projective module in
$\mathcal O^d_\bc(\lambda)$.
For that, we proceed by induction on the level $l$, the case $l=1$
being immediate.

For the induction step, 
take any $A \in \Tab^d(\lambda)$.
Let $l' = l-1, l'' = 1$ and define $\lambda',\lambda'',
\mathfrak{g}'$ and $\mathfrak{g}''$ as above.
Let $\bc' = (c_1,\dots,c_{l-1})$ and write
$\C_{c_l+q_l-n}$ for the one dimensional $\mathfrak{g}''$-module
on which each diagonal matrix unit acts as the scalar $c_l+q_l-n$.
Let $A'$ be the tableau obtained from $A$ by removing the rightmost column.
By the tensor identity, we have that
$$
Z(A,\bc)
\cong 
R(Z(A', \bc') \boxtimes \C_{c_l+q_l-n})
\otimes 
Z^{a_{N-q_l+1}}(V) \otimes\cdots\otimes Z^{a_N}(V).
$$
By induction, $Z(A',\bc')$ is a direct sum of modules of the form
$P(B')$ for various $B' \in \Col_{\bc'}^{d'}(\lambda')$ and 
$0 \leq d' \leq d$.
Therefore, since tensoring with the finite dimensional 
$\mathfrak{g}$-module 
$Z^{a_{N-q_l+1}}(V) \otimes\cdots\otimes Z^{a_N}(V)$
sends projectives to projectives, it suffices to show that
$R(P(B') \boxtimes \C_{c_l+q_l-n})$
is a projective module in $\mathcal O_\bc^{d'}(\lambda)$ for
each $B' \in \Col_{\bc'}^{d'}(\lambda')$ and $d' \geq 0$.
This follows by Corollary~\ref{pi}(ii) applied to the (maximal) $\lambda$-tableau
obtained from $B'$ by adding one more column with entries
$c_l+q_l-n,\dots,c_l+1-n$ from top to bottom.
\end{proof}

\section{Whittaker functor}

Now we return to the representations
of the finite dimensional quotients $W_d(\lambda,\bc)$ 
of $W(\lambda)$ from
section \ref{hlsa}, still assuming that $\bc$ satisfies (\ref{cond}).
We are going to study the
{\em Whittaker functor} from \cite[$\S$8.5]{BKrep},
which will be viewed here
as a functor
from the parabolic category
$\mathcal O(\lambda)$ to modules over the finite $W$-algebra
$W(\lambda)$. The definition 
of Whittaker functor in a more general setting
originates in work of Kostant \cite{Kostant}
and Lynch \cite{Lynch}.

\subsection{\boldmath Polynomial and rational representations of $W(\lambda)$}
For $d \geq 0$, 
let $\mathcal R_\bc^d(\lambda)$ denote the category
of all finite dimensional left 
$W_d(\lambda,\bc)$-modules, viewed as a full subcategory of
the category of all left $W(\lambda)$-modules.
For $r \in \Z$, let $\det^r$ denote the $r$th power of the 
one dimensional determinant representation of $\mathfrak{g}$.
Recalling the functor $? \circledast X$ from $\S$3, 
let $\mathcal R_\bc(\lambda)$ denote the category of all
$W(\lambda)$-modules $M$
with the property that
$$
M \circledast {\det}^r
\cong M_1\oplus\cdots\oplus M_k
$$
for some $r,k,d_1,\dots,d_k \geq 0$ and 
finite dimensional $W_{d_i}(\lambda,\bc)$-modules $M_i$.
Clearly, $\mathcal R_{\bc}(\lambda)$ is an abelian category
closed under taking submodules, quotients and
finite direct sums.
In particular, taking $\bc = \biz$, we call the categories $\mathcal R^d_\biz(\lambda)$ and $\mathcal R_\biz(\lambda)$
the 
categories of {\em polynomial representations of degree $d$} 
and {\em rational representations} of $W(\lambda)$, respectively.

\begin{lem}\label{lm1}
A $W(\lambda)$-module belongs to the category
$\mathcal R_\bc^d(\lambda)$ if and only if it is isomorphic to a subquotient 
of a direct sum of
finitely many copies of $V_\bc^{\otimes d}$.
Hence, 
\begin{itemize}
\item[(i)]
the functor $(? \circledast V)$ maps objects in $\mathcal R_\bc^d(\lambda)$
to objects in $\mathcal R_\bc^{d+1}(\lambda)$;
\item[(ii)]
the functor $(? \circledast E)$ maps objects in 
$\mathcal R_\bc(\lambda)$ to objects in $\mathcal R_\bc(\lambda)$,
for any $\mathfrak{g}$-module arising from 
a finite dimensional rational representation $E$ of $G$. 
\end{itemize}
\end{lem}

\begin{proof}
The first statement is a consequence of the fact that
$V_\bc^{\otimes d}$ is a faithful $W_d(\lambda,\bc)$-module.
To deduce (i), recall
that
$V_\bc^{\otimes d} \circledast V
\cong
V_\bc^{\otimes (d+1)}$.
Finally (ii)  follows from (i)
and associativity, 
because every irreducible rational 
representation of $G$ is a summand of a tensor product of some copies of
$V$ and $\det^{-1}$.
\end{proof}

\begin{cor}\label{ind}
For $r\in\N$, $\mathcal R_\bc^d(\lambda)
\subseteq \mathcal R_{\bc - r\bid}^{d+rN}(\lambda)$
and
$\mathcal R_\bc(\lambda) = \mathcal R_{\bc-r\bid}(\lambda)$.
\end{cor}

\begin{proof}
Observe that
$V^{\otimes d}_\bc
\cong V^{\otimes d}_{\bc -r\bid} \otimes \det^r$.
The latter is isomorphic
to a submodule of $V_{\bc - r\bid}^{\otimes (d+rN)}$, since
$\det^r$ is isomorphic to a submodule of $V^{\otimes rN}$.
\end{proof}

Next we wish to classify the irreducible modules in the categories
$\mathcal R_\bc^d(\lambda)$ and $\mathcal R_{\bc}(\lambda)$.
By  \cite[Theorem 7.9]{BKrep},
there is a finite dimensional irreducible $W(\lambda)$-module
denoted $L(A)$ for each $A \in \Std_\bc(\lambda)$.
Letting $a_{i,1}, \dots, a_{i,p_i}$ denote the entries on the $i$th row
of $A$, 
the irreducible module $L(A)$ is characterized 
uniquely up to isomorphism by the property that
there is a vector $0 \neq v_+ \in L(A)$ such that
\begin{itemize}
\item[(a)] $T_{i,j}^{(r)} v_+ = 0$ for all $1 \leq i < j \leq n$
and $r > p_j-p_i$;
\item[(b)] $T_{i,i}^{(r)} v_+ = e_{r}(a_{i,1}+i-1,\dots,a_{i,p_i}+i-1)$
for all $1 \leq i \leq n$ and $r > 0$,
where $e_{r}(x_1,\dots,x_m)$ denotes the $r$th elementary symmetric
function in $x_1,\dots,x_m$.
\end{itemize}
Such a vector $v_+$ is called a {\em highest weight vector} of type $A$,
and a module generated by a highest weight vector is a 
{\em highest weight module}.
By \cite[Lemma 6.13]{BKrep}, two irreducible modules
$L(A)$ and $L(B)$ for $A, B \in \Std_\bc(\lambda)$
have the same central character if and only if
$\theta(A) = \theta(B)$.

Recall also the finite dimensional $U(\mathfrak{p})$-modules
$V(A)$ for each column strict $A \in \Col_\bc(\lambda)$ defined in the paragraph
preceeding (\ref{nadef}). Often we view $V(A)$ 
as a $W(\lambda)$-module by restriction.
The $W(\lambda)$-modules $\{V(A)\mid  A \in \Col_\bc(\lambda)\}$
are called {\em standard modules}; see
\cite[(7.1)]{BKrep}. 
We warn the reader that the terminology ``standard modules''
being used here has nothing to do with notion of a standard module in a 
highest weight category.
If $A$ is a standard tableau
rather than merely being column strict,
then \cite[Theorem 7.13]{BKrep} shows that
$V(A)$ is a highest weight module generated by a highest weight
vector of type $A$; 
hence it has a unique irreducible quotient, which is isomorphic to $L(A)$.
By \cite[Lemma 6.9]{BKrep}, for any
$A \in \Col_\bc(\lambda)$, every element $z \in Z(W(\lambda))$
acts on the standard module
$V(A)$ as a scalar, i.e. standard modules have central characters.
Again, it is the case that two standard modules
$V(A)$ and $V(B)$ for $A, B \in \Col_\bc(\lambda)$ have the same central character if and only if $\theta(A) = \theta(B)$.

\begin{lem}\label{ic}
For any $A \in \Col_\bc^d(\lambda)$,
$V(A)$ belongs to the category
$\mathcal R_\bc^d(\lambda)$.
The modules
$\{L(A)\mid A \in \Std_\bc^d(\lambda)\}$
form a complete set of irreducible modules in 
$\mathcal R_{\bc}^d(\lambda)$.
\end{lem}

\begin{proof}
As noted in the proof of Lemma~\ref{allirr}, 
$V_\bc^{\otimes d}$ has a filtration as a $\mathfrak{p}$-module,
hence also as a $W(\lambda)$-module,
with sections of the form 
$V(A)$ for $A \in \Col_\bc^d(\lambda)$, and each such $V(A)$ appears
at least once.
Hence by Lemma~\ref{lm1},
each standard module $V(A)$ 
for $A \in \Col_\bc^d(\lambda)$
and each irreducible module
$L(A)$ for $A \in \Std_\bc^d(\lambda)$
belongs to $\mathcal R_\bc^d(\lambda)$.
Conversely, we have already observed that every composition factor
of $N(A)$ for $A \in \Col_\bc^d(\lambda)$ is of the form $K(B)$
for $B \in \Col_\bc^d(\lambda)$.
By \cite[Lemma 8.20]{BKrep} and \cite[Corollary 8.24]{BKrep}, it follows that
every  composition factor of $V(A)$ for $A \in \Col_\bc^d(\lambda)$
is of the form $L(B)$ for $B \in \Std_\bc^d(\lambda)$.
Since $V_\bc^{\otimes d}$ has a filtration with sections of
the form $V(A)$ for $A \in \Col_\bc^d(\lambda)$, we deduce that all its
composition factors are indeed of the form $L(B)$
for $B \in \Std_\bc^d(\lambda)$.
So these are all of the irreducible objects in
$\mathcal R_\bc^d(\lambda)$.
\end{proof}

\begin{cor}\label{cl2}
For any $A \in \Col_\bc(\lambda)$,
$V(A)$ belongs to the category $\mathcal R_\bc(\lambda)$.
The modules
$\{L(A)\mid A \in \Std_\bc(\lambda)\}$
form a complete set of irreducible modules in
$\mathcal R_\bc(\lambda)$.
\end{cor}

\begin{proof}
If $A \in \Col_\bc(\lambda)$, 
then $A \in \Col_{\bc - r \bid}^d(\lambda)$ for
some $r,d \geq 0$.
Hence $V(A)$ belongs to the category
$\mathcal R_\bc(\lambda) = \mathcal R_{\bc - r \bid}(\lambda)$
by Lemma~\ref{ic} and Corollary~\ref{ind}.
It remains to show that every irreducible module
$L \in \mathcal R_\bc(\lambda)$ is isomorphic to
$L(A)$ for some $A \in \Std_\bc(\lambda)$.
By Lemma~\ref{ic} and the definition of the category,
$L \circledast \det^r \cong L(A)$ for some
$r,d \geq 0$ and $A \in \Std_\bc^d(\lambda)$.
Now observe that
$L \cong L(A) \circledast \det^{-r} \cong L(B)$ where $B$ is obtained
from $A$ by subtracting $r$ from all of its entries.
\end{proof}

\subsection{Axiomatic approach to the Whittaker functor}
Assume from now on that
we are given an exact functor
$\V: \mathcal O_\bc(\lambda)
\rightarrow W(\lambda)\text{-mod}$
and an isomorphism of functors
$\nu:
\V \circ (? \otimes V)
\rightarrow (? \circledast V) \circ \V$
such that
\begin{itemize}
\item[(V1)] for each $A \in \Col_\bc(\lambda)$ we have that
$\V(K(A)) \cong L(A)$ if $A$ is standard
  and $\V(K(A)) = 0$ otherwise;
\item[(V2)] 
the following diagrams commute
for all $M \in \mathcal{O}_\bc(\lambda)$:
$$
\quad\qquad
\begin{CD}
\V(M \otimes V) &@>\nu_{M}>>& \V(M) \circledast V\phantom{,}\\
@V\V(x_M)VV&&@VVx_{\V(M)}V\\
\V(M \otimes V)&@>>\nu_{M}>&\V(M) \circledast V,
\end{CD}
$$$$
\quad\qquad
\begin{CD}
\V(M \otimes V \otimes V)&@>\nu_{M \otimes V}>>&\V(M\otimes V) \circledast V
&@>\nu_{M} \circledast \id_V>>&(\V(M) \circledast V) \circledast V\phantom{.}\\
@V\V(s_M)VV&&&&&&@VVs_{\V(M)}V\\
\V(M \otimes V \otimes V)&@>>\nu_{M \otimes V}>&\V(M\otimes V) \circledast V
&@>>\nu_{M} \circledast \id_V>&(\V(M) \circledast V) \circledast V.
\end{CD}
$$
\end{itemize}
Here $x$ and $s$ are the endomorphisms of $(? \otimes V)$
and $(? \otimes V) \otimes V$ defined in
(\ref{xd})--(\ref{sd}), and their analogues
for $(? \circledast V)$ and $(? \circledast V) \circledast V$.
The existence of such a pair $(\V, \nu)$
is established in \cite[$\S$8.5]{BKrep}:
for $\V$ one can take the 
{\em Whittaker functor} from \cite[Lemma 8.20]{BKrep}
and for $\nu$
one can take the isomorphism
from \cite[Lemma 8.18]{BKrep}.
For these choices, 
property (V1) is \cite[Theorem 8.21]{BKrep}, and property
(V2) is checked in the proof of \cite[Lemma 8.19]{BKrep}.
In fact, as we will explain precisely 
in Theorem~\ref{unicity} below, the above properties
determine $(\V, \nu)$ uniquely up to isomorphism.

Property (V1) implies in particular 
that $\V(P_\bc)$ is isomorphic to the $W(\lambda)$-module
$\C_\bc$. We fix such an isomorphism
$i_0:\V(P_\bc) \stackrel{\sim}{\rightarrow} \C_\bc$.
Recursively define an isomorphism
$i_d: \V(P_\bc \otimes V^{\otimes d}) \stackrel{\sim}{\rightarrow}
\C_\bc \circledast V^{\circledast d}$
for each $d \geq 1$
by setting $i_d = (i_{d-1} \circledast \id_V)
\circ \nu_{P_\bc \otimes V^{\otimes (d-1)}}$.
Composing with the isomorphism $\mu_d$ from (\ref{mud}), we get an isomorphism
\begin{equation}\label{bimod}
j_d: \V(P_\bc \otimes V^{\otimes d}) \stackrel{\sim}{\longrightarrow}
V_\bc^{\otimes d},\qquad j_d = \mu_d \circ i_d
\end{equation}
of $W(\lambda)$-modules.
Recall that 
$V_\bc^{\otimes d}$ 
is a $(W(\lambda), H_d)$-bimodule;
the right action
of $H_d$ was defined via (\ref{sk2}) and (\ref{mud}). 
Using (\ref{sk1}), we also make
 $P_\bc \otimes V^{\otimes d}$ into a $(U(\mathfrak{g}), H_d)$-bimodule.
By functoriality, 
we get an induced 
$(W(\lambda), H_d)$-bimodule structure 
on $\V(P_\bc \otimes V^{\otimes d})$  too.
The following lemma is a consequence of (V2)
and the definitions.

\begin{lem}\label{jd}
$j_d: \V(P_\bc \otimes V^{\otimes d}) \rightarrow V_\bc^{\otimes d}$
is an isomorphism of $(W(\lambda), H_d)$-bimodules.
\end{lem}

\subsection{Proof of Theorems A and B}
To prove Theorems A and B from the introduction, 
there is just one more important ingredient,
provided by the next lemma.
We postpone the proof of this until Theorem~\ref{defer} in the next section, preferring 
first to describe the consequences.

\begin{lem}\label{pt}
$\C_\bc$ is a projective module in
$\mathcal R_\bc(\lambda)$.
\end{lem}

We let
$Q(A) = \V(P(A))$ for each $A \in \Std_\bc(\lambda)$.

\begin{lem} \label{pim}
The category $\mathcal R_\bc(\lambda)$ has enough projectives.
More precisely,
for $A \in \Std_\bc(\lambda)$,
 $Q(A)$ is the projective cover of
$L(A)$ in the category $\mathcal R_\bc(\lambda)$.
Moreover, if $A \in \Std_\bc^d(\lambda)$ for some $d \geq 0$, then
$Q(A)$ belongs to
$\mathcal R_\bc^d(\lambda)$,
so it is also the projective cover of $L(A)$
in the category $\mathcal R_\bc^d(\lambda)$.
\end{lem}

\begin{proof}
We show by induction on $d$ that
$Q(A)$ belongs to $\mathcal R_\bc^d(\lambda)$
and that it is the projective cover of $L(A)$
in $\mathcal R_\bc(\lambda)$,
for every $A \in \Std_\bc^d(\lambda)$.
In view of Corollary~\ref{ind} and the observation that 
every $A \in \Std_\bc(\lambda)$ lies in
$\Std_{\bc - r \bid}^d(\lambda)$ for some $r,d \geq 0$,
this is good enough to prove the lemma.
For the base case $d = 0$, we have that $Q(A) \cong \C_\bc$
so are done by Lemma~\ref{pt}.

Now suppose that $d > 0$ and we are given 
$A \in \Std_\bc^d(\lambda)$.
By Lemma~\ref{ccpt} we can write $A$ 
as $\tilde f_i B$ for some $i \in \C$ and
$B \in \Std_\bc^{d-1}(\lambda)$. 
By the induction hypothesis, $Q(B)$ belongs
to $\mathcal R_\bc^{d-1}(\lambda)$
and is the projective cover of $L(B)$ in $\mathcal R_\bc(\lambda)$.
By \cite[Theorem 4.5]{BKrep}, $K(A)$ is a quotient of
$K(B) \otimes V$
so, 
since $P(B) \otimes V$ is projective,
$P(A)$ is a summand of $P(B) \otimes V$.
Applying $\V$, we get that $Q(A)$ is a summand of
$\V(P(B) \otimes V) \cong Q(B) \circledast V$.
Applying Lemma~\ref{lm1}(i), this shows to start with that
$Q(A)$ belongs to $\mathcal R_\bc^d(\lambda)$.
Recalling that 
$V^*$ is the dual $\mathfrak{g}$-module to $V$,
there is a canonical adjunction between the functors $(? \circledast V)$ 
and $(? \circledast V^*)$; see \cite[(8.11)--(8.12)]{BKrep}.
Since $(? \circledast V^*)$ is exact
and maps objects in $\mathcal R_\bc(\lambda)$ to objects in
$\mathcal R_\bc(\lambda)$ by Lemma~\ref{lm1}(ii),
we deduce that $(? \circledast V)$ maps projective
objects in $\mathcal R_\bc(\lambda)$ to 
projective objects in $\mathcal R_\bc(\lambda)$.
Hence, 
$Q(A)$ is a projective object in $\mathcal R_\bc(\lambda)$.

It just remains to show that
$\dim \Hom_{W(\lambda)}(Q(A), L(D)) = \delta_{A, D}$
for all $D \in \Std_\bc(\lambda)$.
Since $P(B) \otimes V$ is a self-dual projective module,
Theorem~\ref{i} implies that
$P(B) \otimes V \cong \bigoplus_{C \in \Std_\bc(\lambda)} P(C)^{\oplus m_C}$
for some integers $m_C \geq 0$ with $m_A \neq 0$.
Applying $\V$, we get that
$Q(B) \circledast V \cong \bigoplus_{C \in \Std_\bc(\lambda)} Q(C)^{\oplus m_C}$
too. We now proceed to show for each $C, D \in \Std_\bc(\lambda)$
with $m_C \neq 0$ that 
$\dim \Hom_{W(\lambda)}(Q(C), L(D)) = \delta_{C,D}$.
By (V1) and exactness, we have that
$\dim \Hom_{W(\lambda)}(Q(C), L(D)) \geq \delta_{C,D}$.
Hence, $\dim \Hom_{W(\lambda)}(Q(B) \circledast V, L(D))
\geq m_D$, and moreover equality holds if and only if
$\dim \Hom_{W(\lambda)}(Q(C), L(D)) = \delta_{C,D}$
whenever $m_C \neq 0$.
Now we calculate using the induction hypothesis:
\begin{align*}
\dim \Hom_{W(\lambda)}(Q(B) \circledast V, L(D))
&=
\dim \Hom_{W(\lambda)}(Q(B), L(D) \circledast V^*)\\
&=
[L(D)\circledast V^*:L(B)]=
[K(D) \otimes V^*: K(B)]\\
&=
\dim \Hom_{\mathfrak{g}}(P(B), K(D) \otimes V^*)\\
&= \dim \Hom_{\mathfrak{g}}(P(B) \otimes V, K(D))
= m_D.
\end{align*}
This completes the proof.
\end{proof}

\begin{cor}\label{sum}
 The category
$\mathcal R_\bc^d(\lambda)$
is a sum of blocks of $\mathcal R_\bc(\lambda)$.
Conversely, every block of $\mathcal R_\bc(\lambda)$ is a block
of $\mathcal R_{\bc - r \bid}^d(\lambda)$ for some $r,d \geq 0$.
\end{cor}

\begin{proof}
Recalling Lemma~\ref{silly},
the first statement follows from 
Lemma~\ref{pim}, the classification of irreducible modules
and central character considerations.
The second statement follows using Corollary~\ref{ind} too,
because every $A \in \Std_\bc(\lambda)$ belongs to 
$\Std_{\bc - r \bid}^d(\lambda)$ for some $r,d \geq 0$.
\end{proof}

\begin{cor}\label{tt}
The functor $\V$ maps 
objects in $\mathcal O_\bc(\lambda)$
to objects in $\mathcal R_\bc(\lambda)$ and
objects in $\mathcal O_\bc^d(\lambda)$
to objects in $\mathcal R_\bc^d(\lambda)$, i.e. to $W_d(\lambda,\bc)$-modules.
\end{cor}

\begin{proof}
We prove the first statement, the second being similar.
It is enough to show that
$\V(P(A))$ belongs to $\mathcal R_\bc(\lambda)$ for each $A \in 
\Col_\bc(\lambda)$.
We already checked this for $A \in \Std_\bc(\lambda)$ in Lemma~\ref{pim}.
In general, let $I$ be the injective hull of $P(A)$
in $\mathcal O_\bc(\lambda)$.
Since $P(A)$ has a parabolic Verma flag, 
Theorem~\ref{i} implies that $I$ is a direct sum of 
$P(B)$'s for $B \in \Std_\bc(\lambda)$.
Hence $\V (I)$ belongs to $\mathcal R_\bc(\lambda)$.
By exactness of $\V$ this implies that $\V (P(A))$ does too.
\end{proof}

Recall from Theorem~\ref{dcthm}
that the algebra $H_d(\lambda,\bc)$, i.e.
the image of the map (\ref{psid}),
coincides with the endomorphism algebra
$\End_{W(\lambda)}(V_\bc^{\otimes d})^{\op}$.
Taking $\bc= \biz$, the following theorem 
completes the proof of Theorem A
from the introduction.

\begin{thm}\label{cta}
For each $d \geq 0$,
$V_\bc^{\otimes d}$
is a projective generator for $\mathcal R_\bc^d(\lambda)$. Hence, 
the functor
$$
\Hom_{W(\lambda)}(V_\bc^{\otimes d}, ?):
W_d(\lambda,\bc)\text{-\rm{mod}}
\rightarrow
H_d(\lambda,\bc)\text{-\rm{mod}}
$$
is an equivalence of categories.
\end{thm}

\begin{proof}
By Corollaries~\ref{q}, \ref{only} and Theorem~\ref{i},
$P_\bc \otimes V^{\otimes d}$
is isomorphic to a direct sum of the modules $P(A)$
for $A \in \Std_\bc^d(\lambda)$, with each of them 
appearing at least once.
Since $\V(P_\bc \otimes V^{\otimes d}) \cong
V_\bc^{\otimes d}$, it follows easily from Lemma~\ref{pim}
and the classification of irreducible modules
that $V_\bc^{\otimes d}$ is a projective generator.
The second statement is a well known consequence.
\end{proof}

For Theorem B, we just need to combine Theorem~\ref{dcthm}
instead with the bimodule isomorphism
$j_d$ from (\ref{bimod}), recalling Lemma~\ref{jd}.

\begin{lem}\label{snow}
There is an isomorphism
$$\gamma:\Hom_{\mathfrak{g}}(P_\bc \otimes V^{\otimes d}, ?)
\stackrel{\sim}{\longrightarrow} 
\Hom_{W(\lambda)}(V_\bc^{\otimes d}, ?) \circ \V$$
of functors from $\mathcal O_\bc^d(\lambda)$ to
$H_d\text{-\rm{mod}}$.
\end{lem}

\begin{proof}
For $M \in \mathcal O_\bc^d(\lambda)$, 
define
$$
\gamma_M: 
\Hom_{\mathfrak{g}}(P_\bc \otimes V^{\otimes d}, M)
\rightarrow \Hom_{W(\lambda)}(V_\bc^{\otimes d}, \V(M))
$$
to be the map
$f \mapsto \V(f) \circ j_d^{-1}$.
This is an $H_d$-module homomorphism by Lemma~\ref{jd},
so defines a natural transformation $\gamma$ between
the functors
$\Hom_{\mathfrak{g}}(P_\bc \otimes V^{\otimes d}, ?)$
and $\Hom_{W(\lambda)}(V_\bc^{\otimes d}, ?) \circ \V$.
These are both exact functors
from $\mathcal O_\bc^d(\lambda)$ to $H_d\text{-mod}$,
so to see that $\gamma$ is actually an isomorphism,
it suffices to check that it gives an isomorphism
on each irreducible module in $\mathcal O^d_\bc(\lambda)$,
i.e. we need to show that
$\gamma_{K(B)}$ is an isomorphism
for each $B \in \Col^d_\bc(\lambda)$.
This follows from Lemma~\ref{jd} if we can show that the map
\begin{align*}
\hom_{\mathfrak{g}}(P_\bc \otimes V^{\otimes d},
K(B))
\rightarrow \hom_{W(\lambda)}(\mathbb{V}(P_{\bc} \otimes V^{\otimes d}),
\mathbb{V}(K(B))),
\quad
&f \mapsto \mathbb{V}(f)\\\intertext{is an isomorphism for each $B \in \Col^d_\bc(\lambda)$.
By Corollary~\ref{only}, the indecomposable summands of
$P_\bc \otimes V^{\otimes d}$ are of the form $P(A)$ for
$A \in \Std_\bc^d(\lambda)$. Therefore it is enough to
show that the map}
\hom_{\mathfrak{g}}(P(A),
K(B))
\rightarrow \hom_{W(\lambda)}(\mathbb{V}(P(A)),
\mathbb{V}(K(B))),
\quad
&f \mapsto \mathbb{V}(f)
\end{align*}
is an isomorphism
for each $A \in \Std_\bc^d(\lambda)$ and
$B \in \Col_\bc^d(\lambda)$.
The left hand side is zero unless $A = B$, when it is one dimensional.
The same is true for the right hand side,
applying (V1) and the fact from Lemma~\ref{pim} that 
$\mathbb{V}(P(A))$
is the projective cover of
$L(A)$. So if $A \neq B$ both sides are non-zero and the
conclusion is clear. Finally if $A = B$,
 let
$f:P(A) \twoheadrightarrow K(A)$ be any non-zero (hence surjective) map.
Then  $\mathbb{V}(f): \mathbb{V}(P(A))
\twoheadrightarrow \mathbb{V}(K(A))$ is non-zero too as $\mathbb{V}$
is exact  and $\mathbb{V}(K(A)) \cong L(A) \neq 0$.
As both sides are one dimensional this is all we need to 
complete the proof in this case too.
\end{proof}

\begin{cor}\label{pic}
The action of $H_d$ on $P_\bc \otimes V^{\otimes d}$
factors through the quotient $H_d(\Lambda)$, making
$P_\bc \otimes V^{\otimes d}$ into a right $H_d(\Lambda)$-module.
Moreover, letting 
\begin{equation}
\rho_d: H_d(\Lambda) \rightarrow \End_{\C}(P_\bc \otimes V^{\otimes d})^{\op}
\end{equation} 
be the resulting homomorphism, 
the following diagram commutes:
$$
\begin{CD}
&\!\!\!\!\!\!H_d(\Lambda) &\\
\\
\!\!\!\!\!\End_{\mathfrak{g}}(P_\bc \otimes V^{\otimes d})
@>\sim > k_d> &\End_{W(\lambda)}(V_\bc^{\otimes d})
\end{CD}
\begin{picture}(0,0)
\put(-125,0){\makebox(0,0){$\swarrow$}}
\put(-130,10){\makebox(0,0){$\scriptstyle\rho_d$}}
\put(-125,0){\line(1,1){18}}
\put(-56,0){\makebox(0,0){$\searrow$}}
\put(-56,0){\line(-1,1){18}}
\put(-53,11){\makebox(0,0){$\scriptstyle\Psi_d$}}
\end{picture}
$$
where the 
bottom algebra isomorphism $k_d$ is the map $f \mapsto j_d \circ \V(f) \circ j_d^{-1}$.
\end{cor}

\begin{proof}
This follows by applying Lemma~\ref{snow}
to $M = P_\bc \otimes V^{\otimes d}$
and noting that
for each $x \in H_d$,
the endomorphism of $P_\bc \otimes V^{\otimes d}$
defined by right multiplication by $x$
is equal to $x \id_M$.
\end{proof}

Using the commutative diagram from Corollary~\ref{pic},
we can identify the algebra $H_d(\lambda,\bc)
= \Psi_d(H_d(\Lambda)) = \End_{W(\lambda)}(V_\bc^{\otimes d})^{\op}$ 
instead with the image of the homomorphism
$\rho_d: H_d(\Lambda) \rightarrow \End_{\C}(P_\bc \otimes V^{\otimes d})^{\op}$.
The following theorem, which in the special case
$\bc = \biz$
proves Theorem B from the introduction,
follows immediately.

\begin{thm}\label{thmbreal}
$\End_{\mathfrak{g}}(P_\bc \otimes V^{\otimes d})^{\op} = H_d(\lambda,\bc)$.
\end{thm}

\subsection{Unicity of Whittaker functors}
We wish to formulate some
further results characterizing the functor $\V$
and the category $\mathcal R_\bc(\lambda)$.
The following theorem shows that $\V$ is a quotient functor 
in the general sense of \cite[$\S$III.1]{Gab}.

\begin{thm}
The functor $\V: \mathcal O_\bc(\lambda)
\rightarrow \mathcal R_\bc(\lambda)$ 
satisfies the following universal property: given any abelian category $\mathcal C$ and any exact
functor
$F:\mathcal O_\bc(\lambda) \rightarrow \mathcal C$
such that $F(K(A)) = 0$ for all $A \in \Col_\bc(\lambda) \setminus
\Std_\bc(\lambda)$, 
there exists an exact functor
$G:\mathcal R_\bc(\lambda) \rightarrow \mathcal C$
such that $F \cong G \circ \V$.
\end{thm}

\begin{proof}
In view of Corollary~\ref{sum} and the corresponding statement
for $\mathcal O_\bc(\lambda)$ explained after (\ref{expafter}), it suffices to prove that the
restriction $\V:\mathcal O_\bc^d(\lambda) \rightarrow
\mathcal R_\bc^d(\lambda)$ satisfies the analogous universal
property with
$\Col_\bc(\lambda) \setminus
\Std_\bc(\lambda)$
replaced by
$\Col^d_\bc(\lambda) \setminus
\Std^d_\bc(\lambda)$.
Since $P_\bc \otimes V^{\otimes d}$
is projective,
Corollaries~\ref{q} and \ref{only}
together with some well known general theory
imply that
the functor
$\Hom_{\mathfrak{g}}(P_\bc \otimes V^{\otimes d}, ?)$
possesses exactly this universal property.
Hence so does the restriction of $\V$,
thanks to Lemma~\ref{snow}
and Theorem~\ref{cta}.
\end{proof}

\begin{cor}
There is an equivalence of categories between
 $\mathcal R_\bc(\lambda)$
and the quotient category
$\mathcal O_\bc(\lambda) / \mathcal J$, where
$\mathcal J$ is the Serre subcategory of
$\mathcal O_\bc(\lambda)$
generated by the irreducible modules
$\{K(A)\mid A \in \Col_\bc(\lambda) \setminus \Std_\bc(\lambda)\}$.
\end{cor}

Next we turn our attention to the
unicity of the functor $\V$.
There is an obvious embedding
$\iota_d:H_d(\Lambda)\hookrightarrow H_{d+1}(\Lambda)$
under which
the elements $x_i$ and $s_j$ of $H_d(\Lambda)$ map to the elements
of $H_{d+1}(\Lambda)$ with the same name; it 
is easy to see using the basis from 
Lemma~\ref{ze} that this map is indeed injective.
Recalling (\ref{m1}) and Corollary~\ref{pic}, we observe that the following diagram commutes:
$$
\begin{CD}
\End_{\mathfrak{g}}(P_\bc \otimes V^{\otimes d})^{\op}
&@>\lambda_d>>&\End_{\mathfrak{g}}(P_\bc \otimes V^{\otimes (d+1)})^{\op}\\
@A\rho_d AA&&@AA\rho_{d+1} A\\
H_d(\Lambda)&@>\iota_d >>&H_{d+1}(\Lambda)\\
@V\Psi_d VV&&@VV\Psi_{d+1} V\\
\End_{W(\lambda)}(V_\bc^{\otimes d})^{\op}
&@>> \kappa_d>& \End_{W(\lambda)}(V_\bc^{\otimes (d+1)})^{\op}
\end{CD}
\begin{picture}(0,0)
\put(-228,3){\oval(50,80)[l]}
\put(7,3){\oval(50,80)[r]}
\put(-263,4){\makebox(0,0){$\scriptstyle k_d$}}
\put(46,4){\makebox(0,0){$\scriptstyle k_{d+1}$}}
\put(7,-37.65){\makebox(0,0){$\leftarrow$}}
\put(-228,-37.65){\makebox(0,0){$\rightarrow$}}
\end{picture}
$$
Here,
$\lambda_d$ is $f \mapsto f \otimes \id_V$ and
$\kappa_d$ is $f \mapsto  \mu_{V_\bc^{\otimes d}, V} \circ (f \circledast \id_V) \circ \mu_{V_\bc^{\otimes d}, V}^{-1}$.
Obviously, the map $\lambda_d$ is injective, hence so too is
$\kappa_d$. So:
\begin{itemize}
\item 
all three horizontal maps $\lambda_d, \iota_d$ and $\kappa_d$ 
in the diagram are injections;
\item
the outside maps $k_d$ and $k_{d+1}$ are isomorphisms; 
\item 
the remaining vertical maps $\rho_d, \rho_{d+1}, \Psi_d$ and $\Psi_{d+1}$
are surjections.
\end{itemize}
The commutativity of the diagram implies that
the embedding $\iota_d$ factors through the quotients
to induce an embedding
$H_d(\lambda,\bc) \hookrightarrow H_{d+1}(\lambda,\bc)$.
The next lemma connects the resulting induction 
functor $\operatorname{ind}_d^{d+1}  = 
(H_{d+1}(\lambda,\bc) \otimes_{H_d(\lambda,\bc)} ?)$
and the functors
$(? \otimes V)$ and $(? \circledast V)$.

\begin{lem}\label{sun}
There are isomorphisms
\begin{align*}
\sigma: \operatorname{ind}_d^{d+1} \circ \Hom_{W(\lambda)}(V_\bc^{\otimes d}, ?)
\stackrel{\sim}{\longrightarrow}
&\Hom_{W(\lambda)}(V_\bc^{\otimes (d+1)}, ?)
\circ (? \circledast V),\\
\tau: \operatorname{ind}_d^{d+1} \circ \Hom_{\mathfrak{g}}(P_\bc \otimes V^{\otimes d}, ?)
\stackrel{\sim}{\longrightarrow}
&\Hom_{\mathfrak{g}}(P_\bc \otimes V^{\otimes (d+1)}, ?)
\circ (? \otimes V)
\end{align*}
of functors from $\mathcal R_\bc^d(\lambda)$
to $H_{d+1}(\lambda,\bc)\text{-\rm{mod}}$
and from $\mathcal O_\bc^d(\lambda)$
to $H_{d+1}(\lambda,\bc)\text{-\rm{mod}}$, respectively, such that the following
diagram commutes:
$$
\begin{CD}
\operatorname{ind}_d^{d+1}
\circ \Hom_{\mathfrak{g}}(P_\bc \otimes V^{\otimes d}, ?)
&@>\bid\gamma >>&
\operatorname{ind}_d^{d+1}
\circ
\Hom_{W(\lambda)}(V_\bc^{\otimes d}, ?) \circ \V\phantom{.}
\\
&&&&@VV \sigma \bid V\\
&&&&\Hom_{W(\lambda)}(V_\bc^{\otimes (d+1)}, ?) \circ (? \circledast V)
\circ \V\phantom{.}\\
@VVV&&@AA \bid \nu A\\
\Hom_{\mathfrak{g}}(P_\bc \otimes V^{\otimes (d+1)}, ?)
\circ (? \otimes V)&@>>\gamma \bid >&
 \Hom_{W(\lambda)}(V_\bc^{\otimes (d+1)}, ?) \circ \V \circ (? \otimes V).\end{CD}
\begin{picture}(0,0)
\put(-274.28,5){\makebox(0,0){$\scriptstyle \tau$}}
\put(-266.78,-8){\line(0,1){42.5}}
\end{picture}
$$
In particular, the functor $\operatorname{ind}_{d}^{d+1}$ is exact.
Moreover, the definition of $\tau$ does not depend on the choices of
$\V$ or $\nu$.
\end{lem}

\begin{proof}
Let us first define the natural transformation $\sigma$.
For a $W_d(\lambda,\bc)$-module $M$, 
there is a natural $H_d(\lambda,\bc)$-module homomorphism
$$
\Hom_{W(\lambda)}(V_\bc^{\otimes d}, M)
\rightarrow
\Hom_{W(\lambda)}(V_\bc^{\otimes (d+1)}, M \circledast V),
\qquad
f \mapsto (f \circledast \id_V) \circ \mu_{V_\bc^{\otimes d}, V}^{-1}.
$$
By adjointness of tensor and hom, this induces a natural
$H_{d+1}(\lambda,\bc)$-module homomorphism
$$
\sigma_M: H_{d+1}(\lambda,\bc) \otimes_{H_d(\lambda,\bc)}
\Hom_{W(\lambda)}(V_\bc^{\otimes d}, M)
\rightarrow
\Hom_{W(\lambda)}(V_{\bc}^{\otimes (d+1)}, M \circledast V).
$$
In this way, we have defined  a natural transformation $\sigma$ between the 
functors $\operatorname{ind}_d^{d+1}\circ \Hom_{W(\lambda)}(V_\bc^{\otimes d}, ?)$
and $\Hom_{W(\lambda)}(V_\bc^{\otimes (d+1)}, ?) \circ (? \circledast V)$.

Next we prove that $\sigma$ is an isomorphism.
We observe to start with that $\sigma$ gives an isomorphism 
when evaluated on the module $M = V_\bc^{\otimes d}$.
This follows because in that case,
$\mu_{V^{\otimes d}_\bc, V} \circ \sigma_M$
is the canonical isomorphism
between $H_{d+1}(\lambda,\bc) \otimes_{H_d(\lambda,\bc)}
H_d(\lambda,\bc)$ and $H_{d+1}(\lambda,\bc)$.
As $V_{\bc}^{\otimes d}$ is a projective generator
for $\mathcal R_{\bc}^d(\lambda)$ by Theorem~\ref{cta},
we deduce by naturality that $\sigma$ gives an isomorphism when evaluated
on an arbitrary projective module in
$\mathcal R_{\bc}^d(\lambda)$.
Finally let $M$ be any module in $\mathcal R_{\bc}^d(\lambda)$
and let $P_2 \rightarrow P_1 \rightarrow M \rightarrow 0$
be an exact sequence where $P_1$ and $P_2$ are projective.
Since 
 $\operatorname{ind}_d^{d+1}\circ \Hom_{W(\lambda)}(V_\bc^{\otimes d}, ?)$
and $\Hom_{W(\lambda)}(V_\bc^{\otimes (d+1)}, ?) \circ (? \circledast V)$
are both right exact, we get a commutative diagram with exact rows
(we abbreviate $\operatorname{H} = \Hom_{W(\lambda)}$):
$$
\begin{CD}
\operatorname{ind}_d^{d+1} (\operatorname{H}(V_\bc^{\otimes d}, 
P_2)) &\longrightarrow&
\operatorname{ind}_d^{d+1} (\operatorname{H}(V_\bc^{\otimes d}, 
P_1)) &\longrightarrow&
\operatorname{ind}_d^{d+1} (\operatorname{H}(V_\bc^{\otimes d}, 
M)) \longrightarrow 0\\
@V\sigma_{P_2}VV@V\sigma_{P_1}VV@VV\sigma_M V\\
\operatorname{H}(V_\bc^{\otimes (d+1)}, 
P_2\circledast V) &\:\longrightarrow\:&
\operatorname{H}(V_\bc^{\otimes (d+1)}, 
P_1\circledast V) &\:\longrightarrow\:&
\operatorname{H}(V_\bc^{\otimes (d+1)}, 
M\circledast V) \longrightarrow 0\\
\end{CD}
$$
We have already shown that the first two vertical maps are isomorphisms,
hence get that 
the third map $\sigma_M$ is an isomorphism by the five lemma.
This completes the proof that $\sigma$ is an isomorphism of functors.
In particular we deduce at this point that $\operatorname{ind}_d^{d+1}$
is exact, since $? \circledast V$ is exact
and $\Hom_{W(\lambda)}(V_\bc^{\otimes d}, ?)$ and
$\Hom_{W(\lambda)}(V_\bc^{\otimes (d+1)}, ?)$ are equivalences by
Theorem~\ref{cta}.

Finally we define the natural transformation $\tau$.
For each $M \in \mathcal O_\bc^d(\lambda)$, the map
$$
\Hom_{\mathfrak{g}}(P_\bc \otimes V^{\otimes d}, M)
\rightarrow
\Hom_{\mathfrak{g}}(P_\bc \otimes V^{\otimes (d+1)}, M \otimes V),
\qquad
f \mapsto f \otimes \id_V
$$
is a natural $H_d(\lambda,\bc)$-module homomorphism.
So by adjointness of tensor and hom once again, we get induced
a natural $H_{d+1}(\lambda,\bc)$-module homomorphism
$$
\tau_M:H_{d+1}(\lambda,\bc) \otimes_{H_d(\lambda,\bc)}
\Hom_{\mathfrak{g}}(P_\bc \otimes V^{\otimes d}, M)
\rightarrow
\Hom_{\mathfrak{g}}(P_\bc \otimes V^{\otimes (d+1)}, M \otimes V).
$$
This definition does not involve $\V$ or $\nu$ in any way.
Now we check that the given diagram commutes.
In particular, this implies that $\tau$ is an isomorphism,
since we already know that
the other four maps are.
Take any $f \in \Hom_{\mathfrak{g}}(P_\bc \otimes V^{\otimes d}, M)$.
The image of $1 \otimes f$
in 
$\Hom_{W(\lambda)}(V_\bc^{\otimes (d+1)}, \V(M) \circledast V)$
going clockwise around the diagram is
$(\V(f) \circledast \id_V) \circ 
(j_d^{-1} \circledast \id_V) \circ \mu_{V_\bc^{\otimes d}, V}^{-1}$.
The image going counterclockwise is
$\nu_M \circ \V(f \otimes \id_V) \circ j_{d+1}^{-1}$.
By naturality of $\nu$, we have that
$\nu_M \circ \V(f \otimes \id_V)
=
(\V(f) \circledast \id_V) \circ \nu_{P_\bc \otimes V^{\otimes d}}$.
Also, $j_{d+1}^{-1} = i_{d+1}^{-1} \circ \mu_{d+1}^{-1}$,
$i_{d+1}^{-1} = \nu_{P_\bc \otimes V^{\otimes d}}^{-1}
\circ (i_d^{-1} \circledast \id_V)$ and
$\mu_{d+1}^{-1} = (\mu_d^{-1} \circledast \id_V) \circ \mu_{V_\bc^{\otimes d}, V}^{-1}$ by (\ref{bimod}) and (\ref{m1}).
Putting these things together completes the proof. 
\end{proof}

Now we can prove the desired characterization of the 
Whittaker functor $\V$.
This is quite useful: for instance, using this theorem one
can prove that $\V$ commutes with duality in an appropriate sense.

\begin{thm}\label{unicity}
Given another exact functor
$\V':\mathcal O_\bc(\lambda) \rightarrow W(\lambda)\text{-\rm{mod}}$
and an
isomorphism of functors $\nu':\V' \circ (? \otimes V) \rightarrow
(? \circledast V) \circ \V'$
satisfying properties analogous to (V1)--(V2) above,
there is a unique (up to a scalar) isomorphism of functors
$\varphi:\V \rightarrow \V'$ such that 
the natural transformations
$\bid \varphi \circ \nu$ and $\nu' \circ \varphi \bid$
from $\V \circ (? \otimes V)$ to $(? \circledast V) \circ \V'$
are equal, i.e.
the following diagram commutes
$$
\begin{CD}
\V(M \otimes V)&@>\nu_M>>&\V(M) \circledast V\\
@V\varphi_{M \otimes V} VV&&@VV \varphi_M \circledast \id_{V} V\\
\V'(M \otimes V)&@>\nu'_M>>&\V'(M) \circledast V
\end{CD}
$$
for all $M \in \mathcal O_\bc(\lambda)$.
\end{thm}

\begin{proof}
Fix the choice of scalar
by choosing 
an isomorphism $i_0': \V'(P_\bc) \stackrel{\sim}{\rightarrow}
\C_\bc$.
From this, we get a
$(W(\lambda), H_d)$-bimodule isomorphism
$j_d': \V'(P_\bc \otimes V^{\otimes d})
\rightarrow V_\bc^{\otimes d}$
in exactly the same way as (\ref{bimod}).
Then we rerun all the above arguments 
with $\V$ replaced by $\V'$
to get analogues of Lemmas~\ref{snow} and \ref{sun};
we denote the resulting isomorphisms of functors
by $\gamma', \sigma'$ and $\tau'$.

We claim that there exists a unique
isomorphism $\bar\varphi: \V \rightarrow \V'$
between the restrictions of the functors
$\V, \V'$ to the full subcategory
$\bigoplus_{d \geq 0} \mathcal O_\bc^d(\lambda)$
of $\mathcal O_\bc(\lambda)$, such that
$i_0  = i_0'\circ
\bar\varphi_{P_\bc}$ and
$\bid\bar \varphi \circ \nu = \nu' \circ \bar\varphi \bid$.
The theorem exactly as stated can be deduced 
from this claim as follows. 
For any $r \geq 0$, 
the claim with $\bc$ replaced by $\bc - r\bid$ 
implies that there exists a unique isomorphism of functors
$\varphi^{(r)}: \V \rightarrow \V'$
between the restrictions of $\V, \V'$ to the full subcategory
$\bigoplus_{d \geq 0} \mathcal O_{\bc-r\bid}^d(\lambda)$ 
of $\mathcal O_\bc(\lambda)$, such that 
the restriction of $\varphi^{(r)}$ to
$\bigoplus_{d \geq 0} \mathcal O_{\bc}^d(\lambda)$
equals $\bar\varphi$ and
$\bid \varphi^{(r)}\circ \nu = \nu' \circ \varphi^{(r)} \bid$.
Then, to construct $\varphi: \V \rightarrow \V'$ in general,
take any module $M \in \mathcal O_\bc(\lambda)$ and
define $\varphi_M:\V(M) \rightarrow \V'(M)$ to be the map $\varphi^{(r)}_M$,
for any $r\geq 0$ chosen so that
$M$ belongs to the subcategory
$\bigoplus_{d \geq 0} \mathcal O_{\bc -r \bid}^d(\lambda)$.

To prove the claim, let us first
construct $\bar\varphi$.
It suffices by additivity to define it
on $M \in \mathcal O_\bc^d(\lambda)$
for each $d \geq 0$.
The composite map $\gamma_M' \circ \gamma_M^{-1}$ is an
isomorphism between $\Hom_{W(\lambda)}(V_{\bc}^{\otimes d}, 
\V(M))$ and 
$\Hom_{W(\lambda)}(V_{\bc}^{\otimes d}, 
\V'(M))$.
Let $\bar\varphi_M: \V(M) \rightarrow \V'(M)$
be the isomorphism obtained from this 
by first applying the 
functor $V_{\bc}^{\otimes d} \otimes_{H_d(\Lambda)} ?$
and then pushing through the 
isomorphism
$V_\bc^{\otimes d} \otimes_{H_d(\Lambda)} \Hom_{W(\lambda)}(V_\bc^{\otimes d}, ?) \rightarrow \operatorname{Id}$
arising from the canonical adjunction between tensor and hom.
This certainly satisfies the property
that $i_0 = i_0' \circ \bar\varphi_{P_\bc}$.
To check that 
$\bid\bar \varphi \circ \nu = \nu' \circ \bar\varphi \bid$
it suffices to show that the top square in the following diagram
commutes:
$$
\begin{array}{c}
\Hom_{\mathfrak{g}}(P_\bc \otimes V^{\otimes (d+1)}, ?)
\circ (? \otimes V)\\
\\\\
\begin{CD}
\Hom_{W(\lambda)}(V_\bc^{\otimes(d+1)}, ?)
\circ \V \circ (? \otimes V)\!
& @> \bid\bar\varphi\bid >> & \!\Hom_{W(\lambda)}(V_\bc^{\otimes (d+1)}, ?)
\circ \V' \circ (? \otimes V)\\
@V \bid\nu VV&&@VV \bid \nu' V\\
\Hom_{W(\lambda)}(V_\bc^{\otimes(d+1)}, ?) \circ (? \circledast V)
\circ \V \!&@> \bid\bid\bar\varphi  >>& \!\Hom_{W(\lambda)}(V_\bc^{\otimes(d+1)}, ?) \circ
(? \circledast V) \circ \V' \\
@A\sigma \bid AA&&@AA \sigma \bid A \\
\operatorname{ind}_{d}^{d+1}\circ \Hom_{W(\lambda)}(V_\bc^{\otimes d}, ?)
\circ \V\! &@> \bid\bid\bar\varphi >>&\! \operatorname{ind}_{d}^{d+1} \circ
\Hom_{W(\lambda)}(V_\bc^{\otimes d}, ?) \circ \V'
\end{CD}\\\\\\
\operatorname{ind}_d^{d+1} \circ \Hom_{\mathfrak{g}}(P_\bc \otimes V^{\otimes d}, ?)
\end{array}
\begin{picture}(0,0)
\put(-95,-57){\makebox(0,0){$\nearrow$}}
\put(-90,-68){\makebox(0,0){$\scriptstyle \bid \gamma'$}}
\put(-95,-57){\line(-1,-1){15}}
\put(-95,60){\makebox(0,0){$\searrow$}}
\put(-90,68){\makebox(0,0){$\scriptstyle\gamma' \bid$}}
\put(-95,60){\line(-1,1){15}}
\put(-265,60){\makebox(0,0){$\swarrow$}}
\put(-270,68){\makebox(0,0){$\scriptstyle\gamma \bid$}}
\put(-265,60){\line(1,1){15}}
\put(-265,-57){\makebox(0,0){$\nwarrow$}}
\put(-268,-68){\makebox(0,0){$\scriptstyle\bid \gamma$}}
\put(-265,-57){\line(1,-1){15}}
\put(-181,-68){\line(0,1){24}}
\put(-181,-26){\line(0,1){24}}
\put(-181,14){\line(0,1){24}}
\put(-181,55){\line(0,1){12}}
\put(-180.95,71){\makebox(0,0){$\uparrow$}}
\put(-188,64){\makebox(0,0){$\scriptstyle \tau$}}
\end{picture}
$$
This follows because the top and bottom triangles
commute by definition of $\bar\varphi$, 
the left and right hand pentagons commute by Lemma~\ref{sun},
and the bottom square commutes by naturality of $\sigma$.

It just remains to prove that $\bar\varphi$ is unique.
Taking $M = P_\bc \otimes V^{\otimes d}$, the 
definitions force
$\bar\varphi_M: \V(M)
\rightarrow \V'(M)$ to be the composite
$(j_d')^{-1} \circ j_d$.
As we noted in the proof of Corollary~\ref{tt},
any injective module in $\mathcal O_\bc^d(\lambda)$
embeds into a direct sum of copies of $M$, so this is good enough.
\end{proof}

For use in the next section, 
we record one more useful fact about $\V$.

\begin{lem}\label{vd}
For any finite dimensional $\mathfrak{p}$-module $M$,
there is a natural $W(\lambda)$-module isomorphism
between $\V(U(\mathfrak{g}) \otimes_{U(\mathfrak{p})} (\C_{-\rho}  \otimes M))$
and the restriction of $M$ from $U(\mathfrak{p})$
to $W(\lambda)$.
In particular,
for $A \in \Col^d_\bc(\lambda)$,
there
is a $W(\lambda)$-module isomorphism
$\V(N(A)) \cong V(A)$.
\end{lem}

\begin{proof}
This is proved in \cite[Lemmas 8.16--8.17]{BKrep}
for one particular choice of the functor $\V$.
It follows for all other choices by Theorem~\ref{unicity}.
\end{proof}

\subsection{Degenerate analogue of the
Dipper-Mathas equivalence}
We give one application of Theorem~\ref{thmbreal}
to recover the degenerate analogue of a theorem
of Dipper and Mathas \cite{DM}. This theorem
reduces most questions about the
algebras $H_d(\Lambda)$ to the
case that $\Lambda$ is an integral weight as in the introduction.

\begin{thm}\label{dmthm}
Suppose $l = l' + l''$ and $\lambda' = (p_1',\dots,p_n')$ 
resp. $\lambda'' = (p_1'',\dots,p_n'')$
is the partition whose diagram consists of the leftmost $l'$
resp. rightmost $l''$ columns of $\lambda$.
Let $\bc' = (c_1,\dots,c_{l'})$ and $\bc'' = (c_{l'+1},\dots,c_{l})$,
and assume that no entry of $\bc'$ lies in the same coset
of $\C$ modulo $\Z$ as an entry of $\bc''$.
Then there is an algebra isomorphism
$$
H_d(\lambda,\bc) \cong \bigoplus_{d'+d''=d}
M_{\binom{d}{d'}}\left(H_{d'}(\lambda',\bc') \otimes H_{d''}(\lambda'',\bc'')\right),
$$
i.e. $H_d(\lambda,\bc)$
is isomorphic to a direct sum of matrix algebras over tensor
products of degenerate cyclotomic Hecke defined relative
to $\bc'$ and $\bc''$.
\end{thm}

\begin{proof}
Let 
$\mathfrak{g}' = \mathfrak{gl}_{q_1+\cdots+q_{l'}}(\C)$,
$\mathfrak{g}'' = \mathfrak{gl}_{q_{l'+1}+\cdots+q_l}(\C)$,
and take all other notation from Corollary~\ref{equivcat}.
Note that $P_\bc \cong R(P_{\bc'} \boxtimes P_{\bc''})$.
Let $V = V' \oplus V''$ be the obvious decomposition of the natural
$\mathfrak{g}$-module $V$ as the direct sum of the natural
$\mathfrak{g}'$-module $V'$ and the natural $\mathfrak{g}''$-module $V''$.
Observe as 
a $\mathfrak{q}$-module that $V^{\otimes d}$ has a filtration whose
sections are isomorphic to a direct sum of $\binom{d}{d'}$ copies of
${V'}^{\otimes d'} \boxtimes {V''}^{\otimes d''}$, one such
for each decomposition $d = d' + d''$.
So by the tensor identity, we deduce that
$P_\bc \otimes V^{\otimes d}$
has a filtration whose sections
are isomorphic to a direct sum of
$\binom{d}{d'}$ copies of
$$
R\left(P_{\bc'} \otimes {V'}^{\otimes d'} \boxtimes
P_{\bc''} \otimes {V''}^{\otimes d''}\right),
$$
one such for each decomposition $d=d'+d''$.
Recalling the block decomposition from Corollary~\ref{equivcat},
the latter module belongs to the subcategory
$\mathcal O_{\bc}^{d',d''}(\lambda)$.
Hence, the filtration splits and we have proved that
$$
P_\bc \otimes V^{\otimes d}
\cong 
R\left(
\bigoplus_{d'+d''=d} 
\binom{d}{d'}P_{\bc'} \otimes {V'}^{\otimes d'} \boxtimes
P_{\bc''} \otimes {V''}^{\otimes d''}\right).
$$
Computing the endomorphism algebra of the right hand side
using the fact from Corollary~\ref{equivcat} 
that $R$ is an equivalence of categories, we get from this that
\begin{multline*}
\End_{\mathfrak{g}}(P_\bc \otimes V^{\otimes d})
\cong\\
\bigoplus_{d'+d''=d}
M_{\binom{d}{d'}}\left(\End_{\mathfrak{g'}}\left(P_{\bc'} \otimes {V'}^{\otimes d'}\right)
\otimes 
\End_{\mathfrak{g''}}\left(P_{\bc''} \otimes {V''}^{\otimes d''}\right)\right).
\end{multline*}
Combining this with Theorem~\ref{thmbreal} completes the proof.
\end{proof}

\begin{cor}
Suppose $\Lambda = \Lambda'+\Lambda''$,
where $\Lambda' = \Lambda_{a_1}+\cdots+\Lambda_{a_k}$
and $\Lambda'' = \Lambda_{b_1}+\cdots+\Lambda_{b_l}$
for some $a_1,\dots,a_k,b_1,\dots,b_l \in \C$
such that no element from the set
$\{a_1,\dots,a_k\}$ lies in the same coset of $\C$ modulo $\Z$
as an element from the set $\{b_1,\dots,b_l\}$.
Then
$$
H_d(\Lambda)
\cong \bigoplus_{d'+d''=d} M_{\binom{d}{d'}}\left(H_{d'}(\Lambda') \otimes H_{d''}(\Lambda'')\right).
$$
\end{cor}

\section{Permutation modules and Specht modules}

Continue to assume 
that the origin $\bc = (c_1,\dots,c_l) \in \C^l$
satisfies (\ref{cond}).
The immediate goal is to prove Lemma~\ref{pt}.
To do this, we are going to exploit 
some weight idempotents
belonging to the algebra $W_d(\lambda,\bc)$
like in \cite[ch.~3]{Green}.
Then we use similar techniques to 
deduce the relationship between
various natural families of modules in $\mathcal O_\bc^d(\lambda)$
and $H_d(\lambda,\bc)$-mod.

\subsection{\boldmath Proof of Lemma~5.6}
Recall that $\Tab^d(\lambda)$ denotes the set of all
$\lambda$-tableaux with non-negative integer entries
summing to $d$.
For $A \in \Tab^d(\lambda)$
with column reading $\gamma(A) = (a_1,\dots,a_N)$, let
$$
\bi(A) = 
(\underbrace{1,\dots,1}_{a_1\text{ times}},\underbrace{2,\dots,2}_{a_2\text{ times}},\dots, \underbrace{N, \dots,N}_{a_N\text{ times}})
\in I^d.
$$
If all entries of $A$ except the rightmost entries in each row are zero, i.e.
we have that $a_i \neq 0 \Rightarrow R(i) = \varnothing$ for each $i=1,\dots,N$,
then we call $A$ an {\em idempotent tableau}.
Let $\Idem^d(\lambda)$ denote the set of all idempotent tableaux
in $\Tab^d(\lambda)$.
For $A \in \Idem^d(\lambda)$, 
the pair $(\bi(A), \bi(A))$ belongs to $J^d$, so it makes sense to define
$$
e_A = \Xi_{\bi(A), \bi(A)} \in W_d(\lambda,\bc),
$$
recalling Theorem~\ref{basis}. In terms of matrix units, 
Lemma~\ref{rec} implies that
\begin{equation}\label{idem}
e_A =
\sum_{\substack{\bi \in I^d \\ \row(\bi) \sim \row(\bi(A))}}
e_{\bi,\bi}.
\end{equation}
It follows immediately that the elements
$\{e_A\mid A \in \Idem^d(\lambda)\}$
are mutually orthogonal idempotents summing to $1$.
For any left $W_d(\lambda,\bc)$-module $M$, we therefore have a vector 
space decomposition
\begin{equation}\label{dec}
M = \bigoplus_{A \in \Idem^d(\lambda)} e_A M.
\end{equation}
(We remark that this decomposition
is precisely the weight space decomposition of $M$ in the sense of
\cite[(5.1)]{BKrep}.)
Similarly, we have that
\begin{equation}\label{pdec}
W_d(\lambda,\bc) = \bigoplus_{A \in \Idem^d(\lambda)}
W_d(\lambda,\bc) e_A,
\end{equation}
so each $W_d(\lambda,\bc) e_A$ is a projective left
$W_d(\lambda,\bc)$-module. 

For the next lemma, 
recall the $U(\mathfrak{p})$-modules $Z^A_{\bc}(V)$ from (\ref{after}),
for each $A \in \Tab^d(\lambda)$.
We will from now on view these as $W(\lambda)$-modules
by restriction. Since $Z^A_{\bc}(V)$ is a submodule of $V_{\bc}^{\otimes d}$,
the $W(\lambda)$-module $Z^A_{\bc}(V)$ belongs to the category
$\mathcal R^d_\bc(\lambda)$, i.e. it is the inflation of
a $W_d(\lambda,\bc)$-module.

\begin{lem}\label{is}
For $A \in \Idem^d(\lambda)$,
the map
$$
W_d(\lambda,\bc) e_A \rightarrow Z^A_\bc(V),
\quad
x e_A \mapsto x v_{\bi(A)}
$$
is an isomorphism of left $W_d(\lambda,\bc)$-modules.
\end{lem}

\begin{proof}
Note first that $v_{\bi(A)} \in e_A V_\bc^{\otimes d}$
does indeed belong to the submodule $Z^A_\bc(V)$, so this is a well-defined
homomorphism.
Letting $d_j$ be the sum of the entries on the $j$th row of $A$, 
the stabilizer of $\bi(A)$
in $S_d$ 
is the parabolic subgroup $S_{d_1} \times\cdots \times S_{d_n}$.
Pick a set $I(A)$ of representatives for
the orbits of
$S_{d_1} \times\cdots\times S_{d_n}$
on the set $\{\bi \in I^d\mid (\bi,\bi(A)) \in J^d\}$.
Then by Theorem~\ref{basis}, the vectors
$\{\Xi_{\bi,\bi(A)}\mid \bi \in I(A)\}$
form a basis for $W_d(\lambda,\bc) e_A$.
Moreover, the definition of $\Xi_{\bi,\bi(A)}$ gives that
$$
\Xi_{\bi,\bi(A)} v_{\bi(A)}
= 
\sum_{(\bj,\bi(A)) \sim (\bi,\bi(A))}
v_{\bj}.
$$
These vectors for all $\bi \in I(A)$ form a basis for
$Z^A_\bc(V)$.
\end{proof}

\begin{thm}\label{defer}
For any $r \geq 0$ and $d = Nr$,
$\C_\bc \circledast \det^{r}$ is a projective
$W_{d}(\lambda,\bc)$-module.
\end{thm}

\begin{proof}
Note that $\C_{\bc} \circledast 
\det^{r} \cong \C_{\bc} \otimes \det^{r}
\cong \C_{\bc+r\bid}$ 
is the irreducible $W_d(\lambda,\bc)$-module
$L(A_{\bc + r \bid})$.
No other tableau in $\Col^d_\bc(\lambda)$
has the same content as $A_{\bc + r \bid}$, as follows easily on recalling
the explicit definition of $A_{\bc}$ from $\S$4.
This means that $L(A_{\bc + r\bid})$ 
is the only irreducible module
in its block, so we just need to show that
$L(A_{\bc+r\bid})$ appears with multiplicity one
as a composition factor of the regular 
module $W_d(\lambda,\bc)$.
Moreover, it shows that $V(A_{\bc+r\bid}) \cong L(A_{\bc + r \bid})$,
and for $A_{\bc + r \bid} \neq B \in \Col^d_{\bc}(\lambda)$
the standard module $V(B)$ 
does not even have $L(A_{\bc+r\bid})$ as a composition factor.

The first part of Lemma~\ref{s} shows that
$Z_\bc^A(V)$ has a filtration as a $\mathfrak{p}$-module with sections
of the form $V(B)$ for $B \in \Col_\bc^d(\lambda)$,
each $V(B)$ appearing with multiplicity $K_{B,A}$.
Hence it also has such a filtration when viewed as a 
$W_d(\lambda,\bc)$-module.
Therefore, by Lemma~\ref{is} and (\ref{pdec}),
the regular module $W_d(\lambda,\bc)$ has a filtration
with sections of the form $V(B)$ for $B \in \Col_{\bc}^d(\lambda)$,
each $V(B)$ appearing with multiplicity 
$\sum_{A \in \Idem^d(\lambda)} K_{B,A}$.

Combining the conclusions of the previous two paragraphs,
it just remains to show that 
$K_{A_\bc+r\bid, A}$ is zero for all but one $A \in \Idem^d(\lambda)$,
and for that $A$ it is equal to one.
This is an easy combinatorial exercise starting from the definition
of the generalized Kostka numbers before the statement of Lemma~\ref{s};
the only $A \in \Idem^d(\lambda)$ for which $K_{A_\bc+r\bid, A}$ is non-zero
is the one whose $j$th row has rightmost entry $r p_j$
for each $j$ with $p_j > 0$.
\end{proof}

Lemma~\ref{pt} from the previous section follows immediately
from Theorem~\ref{defer} and the definition of the category
$\mathcal R_\bc(\lambda)$, so we have now completed the proofs
of Theorems A and B from the introduction.

\subsection{The row removal trick}
The next goal is to prove that
$H_d(\lambda,\bc)$ is a sum of blocks of $H_d(\Lambda)$,
as remarked in the introduction.
Pick $\bar n \leq n$ and let
\begin{equation}\label{hatlambda}
\bar\lambda = (p_1,\dots,p_{\bar n})
\end{equation}
be the partition
obtained by removing the largest $(n-\bar n)$ parts 
from our fixed partition $\lambda$.
Let $\bar N = N - p_{\bar n+1}- \cdots - p_n$,
$\bar I = \{1,\dots,\bar N\}$ and
$\bar{\mathfrak{g}} = \mathfrak{gl}_{\bar N}(\C)$
with natural module $\bar V$.
Working from the diagram of $\bar \lambda$,
we define a grading on $\bar{\mathfrak{g}}$, hence subalgebras
$\bar{\mathfrak{m}}, \bar{\mathfrak{h}}$ and $\bar{\mathfrak{p}}$,
as usual.
The finite $W$-algebra $W(\bar \lambda)$ is
then a subalgebra of $U(\bar{\mathfrak{p}})$.
We denote the analogue of tensor space for $W(\bar\lambda)$
by $\bar V_{\bc}^{\otimes d}$,
with basis $\{v_\bi\mid \bi \in \bar I^d\}$.
Let $\bar\Phi_{d,\bc}:W(\bar\lambda)
\rightarrow \End_{\C}(\bar V_\bc^{\otimes d})$ be the resulting
representation, with image
$W_d(\bar\lambda,\bc)$.
The same degenerate cyclotomic Hecke
algebra $H_d(\Lambda)$ that acts on the right on
$V_\bc^{\otimes d}$ also acts on the right on
$\bar V_\bc^{\otimes d}$ according to the formula
from Lemma~\ref{action}.
Let 
$\bar\Psi_d:H_d(\Lambda) \rightarrow \End_{\C}(\bar V_\bc^{\otimes d})^{\op}$ be the resulting homomorphism, with image $H_d(\bar\lambda,\bc)$.

Define an embedding $\bar I^d \hookrightarrow I^d$
under which $\bi \in \bar I^d$ maps to the unique 
element $\hat \bi \in I^d$ such that
$\row(\hat \bi)$ and $\col(\hat \bi)$ (computed from the diagram of $\lambda$)
are equal to $\row(\bi)$ and $\col(\bi)$ (computed from the diagram
of $\bar\lambda$). Let $e$ denote the idempotent
\begin{equation}\label{eidem}
e = \sum_{\bi \in \bar I^d} e_{\hat \bi, \hat \bi} \in \End_{\C}(V_\bc^{\otimes d}).
\end{equation}
In the notation of (\ref{idem}), 
we have that $e = \sum_A e_A$
summing over all $A \in \Idem^d(\lambda)$ such that
all entries on rows $\bar n+1,\dots,n$ of $A$ are zero.
This shows that $e$ belongs to the subalgebra $W_d(\lambda,\bc)$,
so it is an $H_d(\Lambda)$-equivariant projection and
$e V_\bc^{\otimes d}$ is an $H_d(\Lambda)$-submodule of $V_\bc^{\otimes d}$.
It is obvious from Lemma~\ref{action} that the map
\begin{equation}
\alpha:\bar V_\bc^{\otimes d} \stackrel{\sim}{\longrightarrow}
eV_{\bc}^{\otimes d}, \qquad
v_{\bi} \mapsto v_{\hat\bi}
\end{equation}
is an isomorphism 
of right $H_d(\Lambda)$-modules.
Hence, there is a unique surjective algebra homomorphism
\begin{equation}\label{pid}
\pi:H_d(\lambda,\bc) \twoheadrightarrow H_d(\bar\lambda,\bc)
\end{equation}
such that $\pi \circ \Psi_d = \bar\Psi_d$.
Let
\begin{equation}
\beta:\End_{\C}(\bar V_{\bc}^{\otimes d})
\stackrel{\sim}{\longrightarrow} 
\End_{\C}(eV_\bc^{\otimes d})=
e\End_{\C}(V_\bc^{\otimes d})e,
\qquad
e_{\bi,\bj} \mapsto e_{\hat\bi,\hat\bj}
\end{equation}
be the natural algebra isomorphism induced by $\alpha$.

\begin{lem}\label{r1}
The restriction of the map $\beta$
to the subalgebra $W_d(\bar\lambda,\bc)$ gives an algebra isomorphism
$\beta:W_d(\bar\lambda,\bc) \rightarrow e W_d(\lambda,\bc) e$
with $\beta(\Xi_{\bi,\bj})= \Xi_{\hat \bi, \hat \bj}$.
\end{lem}

\begin{proof}
Using Theorem~\ref{dcthm}, $\beta$ maps the subalgebra
$W_d(\bar\lambda,\bc) = \End_{H_d(\Lambda)}(\bar V_\bc^{\otimes d})$ 
to
$e \End_{H_d(\Lambda)}(V_\bc^{\otimes d}) e = e W_d(\lambda,\bc) e$.
The fact that $\beta$ maps the basis element $\Xi_{\bi,\bj}$
of $W_d(\bar\lambda,\bc)$ to the element $\Xi_{\hat\bi,\hat\bj}$
is now clear from Theorem~\ref{basis}.
\end{proof}

There are also algebra embeddings
\begin{align}
\gamma:W(\bar\lambda)\hookrightarrow W(\lambda),
\qquad
T_{i,j}^{(r)} &\mapsto T_{i,j}^{(r)},\\
\delta:U(\bar{\mathfrak{p}}) \hookrightarrow U(\mathfrak{p}),
\qquad\quad\!
e_{i,j} &\mapsto e_{\hat i, \hat j}.
\end{align}
For a proof that the first of these is well-defined, 
see the first paragraph of \cite[$\S$6.5]{BKrep}.
We stress that $\gamma$ is definitely 
{\em not} the same as the restriction of $\delta$: that usually 
does not even map $W(\bar\lambda)$ into $W(\lambda)$.
Nevertheless, we still have the following result.

\begin{lem}\label{r2}
Viewing the $\mathfrak{p}$-module
$V_\bc^{\otimes d}$ as a $\bar{\mathfrak{p}}$-module
via the embedding $\delta$, the subspace $e V_\bc^{\otimes d}$
is a $\bar{\mathfrak{p}}$-submodule
and the map
$\alpha:\bar V_\bc^{\otimes d} \rightarrow e V_\bc^{\otimes d}$
is a $\bar{\mathfrak{p}}$-module isomorphism.
Viewing the $W(\lambda)$-module
$V_\bc^{\otimes d}$ as a $W(\bar\lambda)$-module
via the embedding $\gamma$, the subspace
$e V_\bc^{\otimes d}$
is a $W(\bar\lambda)$-submodule and
$\alpha$
is a $W(\bar\lambda)$-module isomorphism.
Moreover, the following diagram commutes:
$$
\begin{CD}
W(\bar\lambda)&@>\bar\Phi_{d,\bc}>>&W_d(\bar\lambda,\bc)
&
\hookrightarrow&\End_{\C}(\bar V_\bc^{\otimes d})&
@<<<&U(\bar{\mathfrak{p}})\phantom{.}\\
&&&&@V\beta VV@VV\beta V&&\\
&&&&W_d(\lambda,\bc)e&\:\:\:\hookrightarrow\:\:\:&\End_{\C}(V_\bc^{\otimes d}) e&&\\
@VVV&&@A\eps AA@AA \eps A&&@VVV\\
W(\lambda)&@>>\Phi_{d,\bc}>&W_d(\lambda,\bc)
&
\hookrightarrow&\End_{\C}(V_\bc^{\otimes d})&@<<<&U({\mathfrak{p}}).
\end{CD}
\begin{picture}(0,0)
\put(-244.96,5){\makebox(0,0){$\scriptstyle \gamma$}}
\put(-237.46,-9){\line(0,1){43.5}}
\put(-5.25,5){\makebox(0,0){$\scriptstyle \delta$}}
\put(-11.75,-9){\line(0,1){43.5}}
\end{picture}
$$
Here, $\eps:\End_{\C}(V_\bc^{\otimes d})
\rightarrow \End_{\C}(V_\bc^{\otimes d})e = 
\Hom_{\C}(e V_\bc^{\otimes d}, V_\bc^{\otimes d})$
denotes the natural restriction map.
\end{lem}

\begin{proof}
The first statement is obvious.
To prove the second statement (which is actually never used elsewhere in the article), one needs to use
the explicit definition of the generators \cite[(3.15)--(3.17)]{BKrep}
of $W(\bar\lambda)$;
the key point is that for $1 \leq i,j \leq \bar n$ and $x < \bar n$, 
all monomials $e_{i_1,j_1}\cdots e_{i_s,j_s}$
in the expansion of the elements $T_{i,j;x}^{(r)}$ from \cite[(3.11)]{BKrep} 
that have $\row(i_t) > \bar n$ for some $t$ necessarily act as zero on 
$e V_\bc^{\otimes d}$.
Now let us verify the commutativity of the left hand pentagon.
Take $x \in W(\bar\lambda)$
and $v \in \bar V_\bc^{\otimes d}$.
We just need to check that $\Phi_{d,\bc}(\gamma(x))$ and
$\beta(\bar\Phi_{d,\bc}(x))$ act in the same way on $\alpha(v) \in e V_\bc^{\otimes d}$:
$$
\beta(\bar\Phi_{d,\bc}(x)) \alpha(v)
=
\alpha(\bar\Phi_{d,\bc}(x) v) =
\alpha(x v) = \gamma(x) \alpha(v)
= \Phi_{d,\bc}(\gamma(x)) \alpha(v).
$$
The proof that the right hand pentagon commutes
is similar, using the first statement.
The middle two squares commute by definition.
\end{proof}

Let us call $A \in \Col^d_\bc(\lambda)$
{\em contractible} if 
the entry in the $i$th row and $j$th column is equal to
$(c_j+1-i)$ (the same entry as in $A_\bc$)
for all $i=\bar n+1,\dots,n$ and $j=1,\dots,p_i$.
In that case, we define $\bar A \in \Col^d_\bc(\bar\lambda)$ 
to be the
tableau obtained from $A$ by removing the rows 
numbered $\bar n+1,\dots,n$.

\begin{lem}\label{sim}
Identifying $W_d(\bar\lambda,\bc)$
with
$e W_d(\lambda,\bc) e$
via the isomorphism $\beta$ from Lemma~\ref{r1}, there are
$W_d(\bar\lambda,\bc)$-module isomorphisms
$$
e V(A) \cong \left\{
\begin{array}{ll}
V(\bar A)&\text{if $A$ is contractible,}\\
0&\text{otherwise}
\end{array}\right.
$$
for all $A \in \Col^d_\bc(\lambda)$ and
$$
e L(A) \cong \left\{
\begin{array}{ll}
L(\bar A)&\text{if $A$ is contractible,}\\
0&\text{otherwise}
\end{array}\right.
$$
for all $A \in \Std^d_\bc(\lambda)$.
Moreover, given $A, B \in \Col^d_\bc(\lambda)$
such that $A$ is contractible and $B$ is not,
we have that $\theta(A) \neq \theta(B)$.
\end{lem}

\begin{proof}
Using the right hand part of the commutative diagram from 
Lemma~\ref{r2},
it suffices to prove
the first isomorphism at the level of $\bar{\mathfrak{p}}$-modules: 
viewing the irreducible $\mathfrak{p}$-module $V(A)$ as a 
$\bar{\mathfrak{p}}$-module via the embedding $\delta$, we need to show that
the $\bar{\mathfrak{p}}$-submodule $e V(A)$ is either isomorphic to
$V(\bar A)$ or is zero according to whether $A$ is contractible or not.
Since $e V(A)$ is just a sum of certain weight spaces of $V(A)$, this
follows from a well known result at the level of $\mathfrak{h}$-modules.

To deduce the second isomorphism, recall that if $A \in \Std^d_\bc(\lambda)$
then $L(A)$ is the unique irreducible quotient of $V(A)$ (viewed now as a 
$W(\lambda)$-module).
Hence $e L(A)$ is zero if $A$ is not contractible.
If it is contractible, then the highest weight vector of $V(A)$
lies in the subspace $e V(A)$, so $e L(A)$ is not zero.
Moreover, it is irreducible by a well known result \cite[(6.2b)]{Green},
and it is a quotient of $e V(A) \cong V(\bar A)$.
This implies that $e L(A) \cong L(\bar A)$ as required.

For the last statement, suppose that $A$ is contractible and
$\theta(A) = \theta(B)$.
The $n$th row of $A$ contains the entries
$c_1+1-n, \dots, c_l + 1 - n$.
Since $\theta(B) = \theta(A)$, $B$ must contain these entries somewhere,
hence, because of (\ref{cond}), $B$ must contain these entries in its
$n$th row too.
Now remove the $n$th row and repeat the argument to deduce that
$B$ is contractible.
\end{proof}

\begin{thm}\label{block}
There is a central idempotent $f \in H_d(\lambda,\bc)$
such that the restriction of the map $\pi$ from (\ref{pid}) 
defines an algebra isomorphism
$$
\pi: f H_d(\lambda,\bc) f \stackrel{\sim}{\longrightarrow} H_d(\bar\lambda,\bc).
$$
Moreover, identifying $H_d(\bar\lambda,\bc)$ with
$f H_d(\lambda,\bc) f$ in this way,
there is
a natural $H_d(\bar\lambda,\bc)$-module isomorphism
$$
\Hom_{W(\bar\lambda)}(\bar V_\bc^{\otimes d}, eM) 
\cong
f \Hom_{W(\lambda)}(V_\bc^{\otimes d}, M)
$$
for any left $W_d(\lambda,\bc)$-module $M$.
\end{thm}

\begin{proof}
By Lemma~\ref{sim}, given $A, B \in \Std_\bc^d(\lambda)$
such that $A$ is contractible and $B$ is not, the corresponding
irreducible $W_d(\lambda,\bc)$-modules
$L(A), L(B)$ have different central characters.
Hence, there is a central idempotent $f\in W_d(\lambda,\bc)$
such that
\begin{equation}\label{sand}
f L(A) = \left\{
\begin{array}{ll}
L(A)&\text{if $A$ is contractible,}\\
0&\text{otherwise,}
\end{array}
\right.
\end{equation}
for every $A \in \Std_\bc^d(\lambda)$.
Moreover, 
identifying $W_d(\bar\lambda,\bc)$ with 
$e W_d(\lambda,\bc) e$ via the isomorphism $\beta$ from Lemma~\ref{r1},
well known general theory implies that 
the functor $M \mapsto e M$ arising from the idempotent $e=ef=fe$
defines an equivalence between the categories
$f W_d(\lambda,\bc)f$-mod
and $W_d(\bar\lambda,\bc)$-mod.
Since the endomorphism $f$ of $V_\bc^{\otimes d}$
centralizes $W_d(\lambda,\bc)$, it 
is also a 
central idempotent in $H_d(\lambda,\bc)
= \End_{W(\lambda)}(V_\bc^{\otimes d})^{\op}$,
and we have that
\begin{align*}
f H_d(\lambda,\bc) f = \End_{W_d(\lambda,\bc)}(f V_\bc^{\otimes d})^{\op}
&\cong \End_{e W_d(\lambda,\bc) e}(e V_\bc^{\otimes d})^{\op}\\
&\cong \End_{W_d(\bar\lambda,\bc)}(\bar V_\bc^{\otimes d})^{\op} = H_d(\bar\lambda,\bc).
\end{align*}
The composite isomorphism is simply the restriction of $\pi$, hence that is an isomorphism as claimed.
Moreover, for any $W_d(\lambda,\bc)$-module $M$, we have that
\begin{align*}
f\Hom_{W_d(\lambda,\bc)}(V_\bc^{\otimes d},M)
&=\Hom_{W_d(\lambda,\bc)}(V_\bc^{\otimes d},fM)
=
\Hom_{W_d(\lambda,\bc)}(f V_\bc^{\otimes d}, fM)\\
&\cong
\Hom_{e W_d(\lambda,\bc)e}(e V_\bc^{\otimes d}, eM)
\cong
\Hom_{W_d(\bar\lambda,\bc)}(\bar V_\bc^{\otimes d}, eM)
\end{align*}
as $H_d(\bar\lambda,\bc)$-modules.
\end{proof}

\begin{cor}\label{centidem}
$H_d(\lambda,\bc) \cong f H_d(\Lambda) f$
for some central idempotent $f \in H_d(\Lambda)$.
\end{cor}

\begin{proof}
If at least $d$ parts of $\lambda$ are equal to $l$,
then 
$H_d(\lambda,\bc) \cong H_d(\Lambda)$ by Lemma~\ref{ze}.
So the theorem  shows in particular that 
$H_d(\bar\lambda,\bc) \cong f H_d(\Lambda) f$
for some central idempotent $f \in H_d(\Lambda)$.
\end{proof}

\subsection{Divided power modules and permutation modules}
Call $S \in \Idem^d(\lambda)$ a {\em special tableau}
if all its entries are $\leq 1$ and all its non-zero entries appear
in column $l$.
Such tableaux exist if and only if
at least $d$ parts of $\lambda$ are equal to $l$.
In that case, the functor $\Hom_{W(\lambda)}(V_\bc^{\otimes d}, ?)$
is particularly easy to understand, because by the following lemma
it just amounts to projecting onto a certain weight space.
Moreover, using Theorem~\ref{block}, it is usually possible to reduce 
computations involving the functor
$\Hom_{W(\lambda)}(V_\bc^{\otimes d}, ?)$ in general to this special case.

\begin{lem}\label{t}
Let $S \in \Idem^d(\lambda)$ be a special tableau.
Then there is an algebra isomorphism
$$
h:H_d(\Lambda)
\rightarrow
e_S W_d(\lambda,\bc) e_S
$$
under which $x \in H_d(\Lambda)$
maps to the unique element 
$h(x) \in e_S W_d(\lambda,\bc)e_S$ with
$h(x) v_{\bi(S)} = v_{\bi(S)} x$.
Moreover, identifying $H_d(\Lambda)$ with $e_S W_d(\lambda,\bc) e_S$ in this way,
there is a natural left
$H_d(\Lambda)$-module isomorphism
$$
\Hom_{W(\lambda)}(V_\bc^{\otimes d}, M) \stackrel{\sim}{\longrightarrow}
e_S M,
\qquad
\theta \mapsto \theta(v_{\bi(S)})
$$
for any left $W_d(\lambda,\bc)$-module $M$.
\end{lem}

\begin{proof}
By Lemma~\ref{is} we have that
$V_\bc^{\otimes d} \cong W_d(\lambda,\bc) e_S$.
Hence, using Theorem~\ref{dcthm} and Lemma~\ref{ze} too,
we get isomorphisms
\begin{equation*}
e_S W_d(\lambda,\bc) e_S
\cong \End_{W_d(\lambda,\bc)}(W_d(\lambda,\bc) e_S)^{\op}
\cong
\End_{W(\lambda)}(V_\bc^{\otimes d})^{\op}
\cong H_d(\Lambda).
\end{equation*}
This is the isomorphism $h$.
Moreover,
we have that
$$
\Hom_{W(\lambda)}(V_\bc^{\otimes d}, M)
\cong
\Hom_{W_d(\lambda,\bc)}(W_d(\lambda,\bc) e_S, M)
\cong
e_S M
$$
which proves the second statement.
\end{proof}

As a first application of this lemma,
we can compute the image of the divided power modules
$Z(A,\bc)$ under the functor $\Hom_{\mathfrak{g}}(P_\bc \otimes V^{\otimes d}, ?)$.
For $A \in \Tab^d(\lambda)$
with $\gamma(A) = (a_1,\dots,a_N)$,
introduce the following elements of $H_d(\lambda,\bc)$:
\begin{align}
w_A &= \sum_{w \in S_{a_1}\times\cdots\times S_{a_N}} w,\\
x_A &= \prod_{i=1}^N \prod_{j=1}^{a_i}
\prod_{k=\col(i)+1}^l
(x_{a_1+\cdots+a_{i-1}+j} -(c_k+q_k- n)).\label{ma1}
\end{align}
Note these elements commute, indeed, $x_A$ centralizes all elements
of the subgroup $S_{a_1}\times\cdots\times S_{a_N}$. 
The corresponding {\em permutation module} is
the following left ideal of $H_d(\lambda,\bc)$:
\begin{align}\label{pm}
M(A, \bc) &= 
H_d(\lambda,\bc) x_Aw_A.
\end{align}
When $\bc = \biz$, we denote this module simply by
$M(A)$.

\begin{thm}\label{perm1}
For $A \in \Tab^d(\lambda)$, we have that
$$
M(A,\bc) \cong
\Hom_{W(\lambda)}(V_\bc^{\otimes d}, Z_\bc^A(V))\cong
\Hom_{\mathfrak{g}}(P_\bc \otimes V^{\otimes d}, Z(A, \bc)).
$$
Let
$\gamma(A) = (a_1,\dots,a_N)$,
$\col(\bi(A)) = (n_1,\dots,n_d)$
and $S_d / S_{a_1} \times\cdots\times S_{a_N}$ be some choice of coset
representatives.
Then, in the special case that at least $d$ parts of $\lambda$ are equal to $l$,
the vectors $$
\left\{w x_1^{r_1}\cdots x_d^{r_d} x_Aw_A\:\big|\:
0 \leq r_i < n_i,
w \in S_d  / S_{a_1} \times \cdots \times S_{a_N}\right\}
$$
give a basis for $M(A,\bc)$ and
$$\dim M(A,\bc) = 
\displaystyle \frac{d!}{a_1!\cdots a_N!} \col(1)^{a_1}\cdots\col(N)^{a_N}.
$$
\end{thm}

\begin{proof}
Note by Lemma~\ref{vd} that
$\V(Z(A,\bc)) \cong Z_\bc^A(V)$.
Given this, the second isomorphism follows by Lemma~\ref{snow}.
Now consider the first isomorphism in the special case
that at least $d$ parts of $\lambda$ are equal to $l$.
Let $S \in \Idem^d(\lambda)$ be a special tableau.
Identifying $H_d(\Lambda)$ with
$e_S W_d(\lambda,\bc) e_S$
and
$W_d(\lambda,\bc) e_S$ with $V_\bc^{\otimes d}$
according to the first part of Lemma~\ref{t} and Lemma~\ref{is},
the map
\begin{equation}\label{bmiso}
H_d(\Lambda) \rightarrow e_S V_\bc^{\otimes d},\qquad
x \mapsto v_{\bi(S)} x = h(x) v_{\bi(S)}
\end{equation}
is a bimodule isomorphism.
It maps the left ideal $H_d(\Lambda) x_Aw_A$ isomorphically
onto the left $e_S W_d(\lambda,\bc) e_S$-submodule
of $e_S V_\bc^{\otimes d}$ generated by the vector
$v_{\bi(S)} x_Aw_A$.
We just need to show this submodule is
equal to
$e_S Z_\bc^A(V)$, since that is isomorphic to
$\Hom_{W(\lambda)}(V_\bc^{\otimes d}, Z_\bc^A(V))$
by the last part of Lemma~\ref{t}.

Define $\bj \in I^d$ 
so that
$\col(\bj) = \col(\bi(A))$ and $\row(\bj) = \row(\bi(S))$.
The key observation is that
\begin{equation}\label{thekey}
v_{\bi(S)} x_A = 
v_{\bj}.
\end{equation}
To see this, use
Lemma~\ref{action} to apply the elements from the product
on the right hand side of (\ref{ma1}) in order
starting with the one involving $x_1$.
In particular, this shows that
$v_{\bi(S)} x_A w_A = v_{\bj} w_A$; 
the latter vector clearly lies in $e_S Z_\bc^A(V)$.
Hence, the $e_S W_d(\lambda,\bc) e_S$-submodule
of $e_S V_\bc^{\otimes d}$ generated by the vector
$v_{\bi(S)} x_Aw_A$ is contained in
$e_S Z_\bc^A(V)$.
Now take $0 \leq r_i < n_i$ and $w \in S_d$.
Define $\bj(w) \in I^d$ so that
$\col(\bj(w)) = \col(\bi(A))$ and $\row(\bj(w)) = \row(\bi(S) \cdot w)$.
Then we have that
$$
v_{\bi(S)} w x_1^{r_1} \cdots x_d^{r_d} x_Aw_A
= 
v_{L_1^{r_1}\circ \cdots \circ L_d^{r_d}(\bj(w))}w_A
+ (*)
$$
where $(*)$ is a sum of vectors of $e_S V_\bc^{\otimes d}$ of strictly
smaller degree in the natural grading.
The vectors 
$$
\left\{v_{L_1^{r_1}\circ \cdots \circ L_d^{r_d}(\bj(w))}w_A
\:\bigg|\:0 \leq r_i < n_i, w \in S_d / S_{a_1}\times\cdots\times S_{a_N}\right\}
$$
give the obvious basis for the weight space
$e_S Z_\bc^A(V)$. Since we know already that
$$
v_{\bi(S)} w x_1^{r_1} \cdots x_d^{r_d} x_Aw_A
= h(w x_1^{r_1}\cdots x_d^{r_d}) v_{\bi(S)} x_Aw_A
$$
lies in $e_S Z_\bc^A(V)$, it follows that the lower terms
$(*)$ all belong to $e_S Z_\bc^A(V)$, and by a triangular
change of basis argument the vectors
$$
\left\{v_{\bi(S)} w x_1^{r_1} \cdots x_d^{r_d} x_Aw_A\:\big|\:
0 \leq r_i < n_i, w \in S_d / S_{a_1}\times\cdots\times S_{a_N}\right\}
$$
form a basis for $e_S Z_\bc^A(V)$.
This proves that the $e_S W_d(\lambda,\bc)e_S$-submodule of
$e_S V_\bc^{\otimes d}$ generated by the vector $v_{\bi(S)} x_Aw_A$
is equal to $e_S Z_\bc^A(V)$, and at the same time 
shows
that the vectors from in the statement of the theorem form a basis
for $M(A,\bc)$.

Finally, 
to deduce the existence of the first isomorphism 
in general from this special case, assume
that $\bar\lambda$ is obtained by removing parts from $\lambda$
like in (\ref{hatlambda}).
Take $A \in \Tab^d(\lambda)$ such that all entries in
rows $\bar n+1,\dots,n$ of $A$ are zero, and let
$\bar A \in \Tab^d(\bar\lambda)$ be the tableau 
obtained by removing these rows.
Let $e$ be the idempotent from (\ref{eidem}).
Identifying $e W_d(\lambda,\bc) e$ with $W_d(\bar\lambda,\bc)$
via the isomorphism $\beta$ from Lemma~\ref{r1}, we 
obviously have that
$e Z^{A}_\bc(V) \cong Z^{\bar A}_\bc(\bar V)$ as $W(\bar\lambda)$-modules.
We have shown already that
$\Hom_{W(\lambda)}(V^{\otimes d}_\bc, Z^{A}_\bc(V))
\cong M(A,\bc)$.
So now Theorem~\ref{block} implies that
$\Hom_{W(\bar\lambda)}(\bar V^{\otimes d}_\bc, Z^{\bar A}_\bc(\bar V))
\cong \Psi_d(M(A,\bc))$ 
as $H_d(\bar\lambda,\bc)$-modules.
It just remains to observe that $\Psi_d(x_A) = x_{\bar A}$, so
$\Psi_d(M(A,\bc)) = M(\bar A,\bc)$.
\end{proof}

\subsection{Cyclotomic Schur algebras and the proof of Theorem C}
Now we let $S_d(\lambda,\bc)$ denote the finite dimensional algebra
\begin{align*}
S_d(\lambda,\bc) &= 
\End_{H_d(\Lambda)}\bigg(\bigoplus_{A \in \Tab^d(\lambda)} M(A,\bc)\bigg)^{\!\op}\cong
\End_{W(\lambda)}\bigg(\bigoplus_{A \in \Tab^d(\lambda)} Z^A_{\bc}(V)\bigg)^{\!\op}
\!\!\!,
\end{align*}
the second isomorphism being a consequence of Theorems~\ref{cta} 
and \ref{perm1}.
As usual, if $\bc = \biz$, we denote this algebra simply by $S_d(\lambda)$
and (for reasons explained shortly) call 
it the {\em cyclotomic Schur algebra}.
The first part of the following theorem is due to Stroppel 
\cite[Theorem 10.1]{St}, 
and can be viewed as a parabolic analogue of Soergel's Struktursatz.

\begin{thm}\label{ucsb}
The functor $\Hom_{\mathfrak{g}}(P_\bc \otimes V^{\otimes d}, ?): \mathcal O_\bc^d(\lambda) \rightarrow H_d(\lambda,\bc)\text{-\rm{mod}}$
is fully faithful on projective modules.
In particular,
it induces an algebra isomorphism
$$
\End_{\mathfrak{g}}\bigg(\bigoplus_{A \in \Tab^d(\lambda)} Z(A,\bc)\bigg)^{\!\op}
\stackrel{\sim}{\longrightarrow}
S_d(\lambda,\bc).
$$
Hence the category $\mathcal O^d_{\bc}(\lambda)$ is equivalent to the
category of all finite dimensional left $S_d(\lambda,\bc)$-modules.
\end{thm}

\begin{proof}
By Theorem~\ref{i} and Corollaries~\ref{only} and \ref{q}, 
$P_\bc \otimes V^{\otimes d}$
is a direct sum of all the self-dual projective indecomposable modules
in $\mathcal O_\bc^d(\lambda)$, each appearing with some non-zero multiplicity.
Applying \cite[Theorem 2.10]{KSX} 
as in the proof of \cite[Theorem 3.2]{KSX}, 
the proof of the first statement reduces to checking that there exists an 
exact sequence
$0 \rightarrow P_0 \rightarrow P_1 \rightarrow P_2$ 
in $\mathcal O^d_{\bc}(\lambda)$ 
such that $P_0$ is a projective generator and 
$P_1, P_2$ are self-dual projectives.
To see this, there is for $A \in \Tab^d(\lambda)$
a short exact sequence of $\mathfrak{p}$-modules
$0 \rightarrow Z_\bc^A(V) \rightarrow V_\bc^{\otimes d} \rightarrow K \rightarrow 0$ for some finite dimensional $\mathfrak{p}$-module $K$.
Inducing, we get a short exact sequence
$0 \rightarrow Z(A,\bc) \rightarrow P_\bc \otimes V^{\otimes d} \rightarrow L
\rightarrow 0$ where $L \in \mathcal O^d_\bc(\lambda)$ has a parabolic Verma flag. 
By Theorem \ref{i}, the injective hull $I$ of $L$
is a self-dual projective module
and we have constructed an exact sequence
$$
0 \rightarrow Z(A,\bc) \rightarrow P_\bc \otimes V^{\otimes d} \rightarrow
I.
$$
By Theorem~\ref{pgen}, the module
$\bigoplus_{A \in \Tab^d(\lambda)} Z(A,\bc)$ is a projective generator
for $\mathcal O^d_{\bc}(\lambda)$, so it just remains
to take the direct sum of these exact sequences 
over all $A \in \Tab^d(\lambda)$.
Moreover, 
the image of the projective generator
$\bigoplus_{A \in \Tab^d(\lambda)} Z(A,\bc)$
under the functor
$\Hom_{\mathfrak{g}}(P_\bc \otimes V^{\otimes d}, ?)$
is isomorphic to $\bigoplus_{A \in \Tab^d(\lambda)} M(A,\bc)$
by Theorem~\ref{perm1}.
So the second statement now follows from the first, and then the final 
statement follows by well known general theory.
\end{proof}

To deduce Theorem C, it just remains to explain the relationship between
the algebra $S_d(\lambda,\bc)$ and the degenerate analogue of the cyclotomic $q$-Schur algebra of Dipper, James and Mathas from \cite{DJM};
see also \cite[$\S$6]{AMR} where the degenerate case is considered 
explicitly.
The key observation 
is that if at least $d$ parts of $\lambda$ are equal to $l$ (when $H_d(\lambda,\bc) = H_d(\Lambda)$) our permutation module
$M(A,\bc)$ is essentially the same as 
the permutation module $M^\bmu$ of Dipper, James and Mathas,
for the parameters $Q_1,\dots,Q_l$ defined from
$Q_i = c_i+q_i-n$ and the $l$-tuple of compositions
$\bmu = (\mu^{(1)},\dots,\mu^{(l)})$ defined so that the parts of
$\mu^{(i)}$ are the entries in the $i$th column of $A$.
More precisely, the permutation module of Dipper, James and Mathas
is exactly the right ideal $w_A x_A H_d(\Lambda)$, which is
isomorphic to our left module $M(A,\bc)$ viewed as a right
module via the antiautomorphism
\begin{equation}\label{star}
*:H_d \rightarrow H_d \qquad 
x_i \mapsto x_i, \:s_j \mapsto s_j.
\end{equation}
Since the degenerate analogue of the cyclotomic $q$-Schur algebra of Dipper, James and Mathas is the endomorphism algebra of a certain direct sum of
their permutation modules, this implies that 
$S_d(\lambda,\bc)$ is at least Morita equivalent to their algebra.
Hence, Theorem C from the introduction is a
special case of Theorem~\ref{ucsb}.
The following lemma further clarifies the relationship 
between permutation modules and the tensor space
$V_{\bc}^{\otimes d}$, recalling (\ref{dec}).

\begin{lem}
Assume that $A \in \Idem^d(\lambda)$.
Then the following are isomorphic as 
right $H_d(\Lambda)$-modules:
\begin{itemize}
\item[(i)] the weight space $e_A V_{\bc}^{\otimes d}$;
\item[(ii)] the right ideal $w_A x_A H_d(\Lambda)$ of $H_d(\Lambda)$;
\item[(iii)] the right ideal $w_A x_A H_d(\lambda,\bc)$ of $H_d(\lambda,\bc)$;
\item[(iv)] the permutation module $M(A,\bc)$ viewed as a right 
module via $*$.
\end{itemize}
\end{lem}

\begin{proof}
Assume to start with that
at least $d$ parts of $\lambda$ are equal to $l$,
so that a special tableau $S \in \Idem^d(\lambda)$ exists
and $H_d(\lambda,\bc) = H_d(\Lambda)$.
Using (\ref{star}), it follows that
the dimension of $w_A x_A H_d(\Lambda)$ is the same
as the dimension of $M(A,\bc)$, which we computed in Theorem~\ref{perm1}.
This agrees with the dimension of the weight space
$e_A V_\bc^{\otimes d}$.
The bimodule isomorphism from (\ref{bmiso}) identifies
$w_A x_A H_d(\Lambda)$ with the right $H_d(\Lambda)$-submodule
of $e_S V_\bc^{\otimes d}$ generated by the vector
$v_{\bi(S)} w_Ax_A$.
Let $\bj \in I^d$ satisfy $\col(\bj) = \col(\bi(A))$
and $\row(\bj) = \row(\bi(S))$.
Left multiplication by $\Xi_{\bj,\bi(A)}$ clearly defines a right
$H_d(\Lambda)$-module homomorphism
$$
\theta:e_A V_\bc^{\otimes d} \rightarrow e_S 
V_\bc^{\otimes d}.
$$ 
By the definition of
$\Xi_{\bj,\bi(A)}$
from Theorem~\ref{basis}, and using (\ref{thekey})
for the penultimate equality, we have that
$$
\theta(v_{\bi(A)}) = 
\Xi_{\bj,\bi(A)} v_{\bi(A)} = \sum_{w \in S_{a_1}\times\cdots\times S_{a_N}} v_{\bj\cdot w}= v_{\bj} w_A = v_{\bi(S)} x_Aw_A = v_{\bi(S)} w_A x_A.
$$
Hence, the image of $\theta$ contains the 
copy
$v_{\bi(S)} w_Ax_A H_d(\Lambda)$ of $w_A x_A H_d(\Lambda)$.
We have already observed that 
$\dim e_A V_\bc^{\otimes d} = \dim w_A x_A H_d(\Lambda)$,
so we now get that $\theta$ is injective and its image is 
isomorphic to
$w_A x_A H_d(\Lambda)$.
This proves the modules (i) and (ii) are isomorphic in this special case.

Now assume that $\bar\lambda$ is obtained by removing parts from $\lambda$,
like in (\ref{hatlambda}).
Take $A \in \Tab^d(\lambda)$ such that all entries in
rows $\bar n+1,\dots,n$ of $A$ are zero, and let
$\bar A \in \Tab^d(\bar\lambda)$ be the tableau 
obtained by removing these rows.
Let $e$ be the idempotent from (\ref{eidem}).
By the previous paragraph, $e_A V_\bc^{\otimes d}$ is isomorphic
to $w_A x_A H_d(\Lambda)$
as a right $H_d(\Lambda)$-module.
On the other hand, recalling Lemma~\ref{r1},
we have that $e_{A} = \beta(e_{\bar A})=e_A e$, so
$$
e_{A} V_\bc^{\otimes d} = e_A e V_\bc^{\otimes d}
= \beta(e_{\bar A}) \alpha(\bar V_\bc^{\otimes d}) = \alpha(e_{\bar A} 
\bar V_\bc^{\otimes d}),
$$
which is isomorphic to $e_{\bar A} \bar 
V_\bc^{\otimes d}$ as a right $H_d(\Lambda)$-module.
This proves that the modules (i) and (ii) are isomorphic in general.

To deduce the remaining statements in the lemma,
note the module (i) is an $H_d(\lambda,\bc)$-module,
so by the isomorphism between (i) and (ii) we see that 
the module (ii) is the inflation of an $H_d(\lambda,\bc)$-module too.
Recalling Theorem~\ref{block}, 
this easily implies that (ii) and (iii) are isomorphic, 
and it is obvious from (\ref{pm}) that (iii) and (iv) are isomorphic.
\end{proof}

\subsection{Parabolic Verma modules and Specht modules}
The final task
is to give an intrinsic description of the 
$H_d(\Lambda)$-modules
which correspond 
to parabolic Verma
modules and standard modules under 
the functors
$\Hom_{\mathfrak{g}}(P_\bc \otimes V^{\otimes d}, ?)$
and $\Hom_{W(\lambda)}(V_\bc^{\otimes d}, ?)$, respectively.
For a partition $\mu$ of $d$,
let $S(\mu)$ denote the usual (irreducible) {\em Specht module}
for the group algebra $\C S_d$ parametrized by the partition $\mu$.
Given $a \in \C$, there is the evaluation homomorphism
\begin{equation}\label{ev}
\ev_a:H_d \twoheadrightarrow \C S_d, \qquad x_1 \mapsto a,\: s_j \mapsto s_j.
\end{equation}
We denote the lift 
of the Specht module $S(\mu)$ to $H_d$ through this map 
 by $\ev_a^*(S(\mu))$ and also refer to it as a Specht module.

Let $A \in \Col^d_{\bc}(\lambda)$.
Define an $l$-tuple of partitions
$\bmu = (\mu^{(1)},\dots,\mu^{(l)})$ by letting
the parts of $\mu^{(i)}$ be the entries in the $i$th column of
$(A - A_{\bc})$.
Say $\mu^{(i)}$ is a partition of $d_i$ for each $i=1,\dots,l$.
We view
$H_{(d_1,\dots,d_l)} = 
H_{d_1}\otimes\cdots\otimes H_{d_l}$ 
as a subalgebra of $H_d$ in the obvious way; see e.g.
\cite[$\S$3.4]{Kbook}. The (outer) tensor product
\begin{equation}
\ev_{c_1+q_1-n}^*(S(\mu^{(1)}))\boxtimes\cdots\boxtimes \ev_{c_l+q_l-n}^*(S(\mu^{(l)}))
\end{equation}
is a left $H_{(d_l,\dots,d_1)}$-module.
Define the {\em Specht module} $S(A)$  by setting
\begin{align}
S(A) &= 
H_d \otimes_{H^{\phantom{A}}_{(d_1,\dots,d_l)}}
\left(\ev_{c_1+q_1-n}^*(S(\mu^{(1)}))\boxtimes\cdots\boxtimes
\ev_{c_l+q_l-n}^*(S(\mu^{(l)}))\right).\label{dsm}
\end{align}
By \cite[3.2.2]{Kbook},
$H_d$ is a free right
$H_{(d_1,\dots,d_l)}$-module with basis given by any set
of $S_d / S_{d_1}\times\cdots\times S_{d_l}$-coset representatives, 
which implies that
\begin{equation}\label{df}
\dim S(A) =
\frac{d!}{d_1! \cdots d_l!} \times \dim S(\mu^{(1)})
\times \cdots \times \dim S(\mu^{(l)}).
\end{equation}
It is quite easy to see directly from the definition (\ref{dsm})
that the action of $H_d$ on $S(A)$
factors through to the quotient
$H_d(\Lambda)$, so $S(A)$ can be viewed as
an $H_d(\Lambda)$-module.
This also follows from the following theorem, which shows moreover that
$S(A)$ is an $H_d(\lambda,\bc)$-module.

\begin{thm}\label{final}
For $A \in \Col^d_{\bc}(\lambda)$, we have that
$$
S(A) \cong
 \Hom_{W(\lambda)}(V_\bc^{\otimes d}, V(A))
\cong
\Hom_{\mathfrak{g}}(P_\bc\otimes V^{\otimes d}, N(A))
$$
as left $H_d$-modules.
\end{thm}

\begin{proof}
The second isomorphism follows from Lemmas~\ref{snow} and \ref{vd}.
To prove the first isomorphism, 
assume to start with that at least $d$ parts of $\lambda$
are equal to $l$, and let $S$ be a special tableau. 
Let $\gamma(A - A_{\bc}) = (a_1,\dots,a_N)$ and
$\bmu=(\mu^{(1)},\dots,\mu^{(l)})$ be the 
$l$-tuple of partitions so that $\mu^{(i)}$ is the partition of $d_i$
whose parts are the
entries in the $i$th column of $A - A_\bc$.
Let $Z = Z^{A - A_\bc}_\bc(V)$. 
The grading on $V_\bc^{\otimes d}$ induces a grading 
$Z = \bigoplus_{i \geq b} Z_i$
on the subspace $Z$; 
the bottom degree is $b = \sum_{i=1}^N a_i (l - \col(i))$.
The subspace $Z_{> b} = \bigoplus_{i > b} Z_i$ is 
stable under the action of $\mathfrak{p}$,
and clearly the quotient $\mathfrak{p}$-module $Z / Z_{> b}$
is isomorphic to the inflation of the 
$\mathfrak{h}$-module $Z_b$.
Letting $V_{(j)}$ denote the span of all
$\{v_i\mid \col(i) = j\}$, we have that
\begin{equation*}
Z_b = 
Z_{c_1}^{\mu^{(1)}}\left(V_{(1)}\right) \boxtimes\cdots 
\boxtimes Z_{c_l}^{\mu^{(l)}}\left(V_{(l)}\right)
\end{equation*}
as an $\mathfrak{h}$-module,
where for a composition $\mu = (\mu_1,\dots,\mu_r)$ and a vector space $E$
we write  
$Z^\mu(E) = Z^{\mu_1}(E) \otimes\cdots\otimes Z^{\mu_r}(E)$.
Hence, by the well known classical theory, there is an
$S_{d_1}\times\cdots\times S_{d_l}$-module isomorphism
$$\Hom_{\mathfrak{h}}(V_{\bc}^{\otimes d}, Z_b)
\cong M(\mu^{(1)}) \boxtimes\cdots\boxtimes M(\mu^{(l)}),$$
where for a composition $\mu = (\mu_1,\dots,\mu_r)$ of $d$ we write
$M(\mu)$ for the usual permutation module for the symmetric group 
$S_d$, i.e. the module induced from the trivial repesentation of 
the parabolic subgroup $S_{\mu_1}\times\cdots\times S_{\mu_r}$.
Composing the natural maps
$$
\Hom_{\mathfrak{h}}(V_{\bc}^{\otimes d},Z_b)
\cong
\Hom_{\mathfrak{p}}(V_{\bc}^{\otimes d}, Z / Z_{> b})
\hookrightarrow
\Hom_{W(\lambda)}(V_{\bc}^{\otimes d}, Z / Z_{> b}),
$$
we have constructed an embedding
$$
f:M(\mu^{(1)}) \boxtimes\cdots\boxtimes
M(\mu^{(l)}) 
\hookrightarrow
\Hom_{W(\lambda)}(V_\bc^{\otimes d}, Z / Z_{> b})
$$
at the level of $S_{d_1}\times\cdots\times S_{d_l}$-modules.
Let us describe the image of $f$ explicitly in coordinates.
Identify
$H_d(\Lambda)$ with $e_S W_d(\lambda,\bc) e_S$
and $\Hom_{W(\lambda)}(V_\bc^{\otimes d}, Z / Z_{> b})$
with $e_S Z / e_S Z_{> b}$ according to Lemma~\ref{t}.
Let $\bj \in I^d$ 
satisfy
$\row(\bj) = \row(\bi(S))$ and $\col(\bj) = \col(\bi(A-A_\bc))$.
Then the image of $f$ is identified with the span of the vectors
$\{v_{\bj} w  w_{A-A_\bc} + e_S Z_{> b}\mid w \in S_{d_1}\times\cdots\times S_{d_l}\}$.
Moreover, by (\ref{thekey}), we know that
$v_{\bj} w w_{A-A_\bc} = v_{\bi(S)} x_{A-A_\bc} w w_{A-A_\bc}$.
Hence, for any $i=1,\dots,l$, we have that
\begin{align*}
h(x_{d_1+\cdots+d_{i-1}+1})
v_{\bj} w w_{A-A_\bc}
&=
v_{\bi(S)} x_{d_1+\cdots+d_{i-1}+1} x_{A-A_\bc} w w_{A-A_\bc}\\
&=
v_{\bi(S)} x_{A-A_\bc} x_{d_1+\cdots+d_{i-1}+1} w w_{A-A_\bc}\\
&= 
v_{\bj} 
x_{d_1+\cdots+d_{i-1}+1} w w_{A-A_\bc}\\
&\equiv (c_i+q_i-n) v_{\bj} w w_{A-A_\bc} \pmod{e_S Z_{> b}},
\end{align*}
the final congruence following by
an explicit calculation using Lemma~\ref{action}.
For a composition $\mu$ of $d$,
let us now write $\ev_a^*(M(\mu))$ for the $H_d$-module
arising by lifting the permutation module
$M(\mu)$ from $\C S_d$ to $H_d$ via the evaluation
homomorphism (\ref{ev}).
Then $\ev_{c_1+q_1-n}^*(M(\mu^{(1)})) \boxtimes\cdots\boxtimes
\ev_{c_l+q_l-n}^*(M(\mu^{(l)}))$ is the unique
lift of $M(\mu^{(1)}) \boxtimes\cdots\boxtimes M(\mu^{(l)})$
to $H_{(d_1,\dots,d_l)}$ such that
$x_{d_1+\cdots+d_{i-1}+1}$ acts as $(c_i+q_i-n)$
for each $i=1,\dots,l$.
So the calculation we just made establishes that 
the map $f$ is actually an embedding
\begin{multline*}
f:\ev_{c_1+q_1-n}^*(M(\mu^{(1)})) \boxtimes\cdots\boxtimes \ev_{c_l+q_l-n}^*(M(\mu^{(l)})) \\\hookrightarrow
\Hom_{W(\lambda)}(V_{\bc}^{\otimes d}, Z / Z_{> b})
\cong e_S Z / e_S Z_{> b}
\end{multline*}
at the level of $H_{(d_1,\dots,d_l)}$-modules.
Moreover, the image of $f$ contains the vector
$v_{\bj} w_{A-A_\bc} + e_S Z_{> b}$ which is a generator
of $e_S Z / e_S Z_{> b}$ as an $H_d$-module.

Now, as an $\mathfrak{h}$-module,
$V(A)$ can be viewed as the (outer) tensor product 
$$
V(A) = V(A_1)\boxtimes\cdots\boxtimes V(A_l),
$$
where $A_i$ is the $i$th column of $A$ and $V(A_i)$ is the
corresponding irreducible $\mathfrak{gl}_{q_i}(\C)$-module.
It is well known that there is a 
surjective $\mathfrak{gl}_{q_i}(\C)$-module homomorphism
$Z_{c_i}^{\mu^{(i)}}\!\!\left(V_{(i)}\right) \twoheadrightarrow V(A_i)$.
So there is a surjective $\mathfrak{h}$-module
homomorphism $Z_b \twoheadrightarrow V(A)$,
hence a surjective $\mathfrak{p}$-module homomorphism
$Z / Z_{> b} \twoheadrightarrow V(A)$.
Moreover, by the classical theory again,
$$
\Hom_{\mathfrak{h}}(V_{\bc}^{\otimes d}, V(A))
\cong S(\mu^{(1)})\boxtimes\cdots\boxtimes S(\mu^{(l)})
$$
as an $S_{d_1}\times\cdots\times S_{d_l}$-module.
Combining this with the conclusion of the previous paragraph, 
we have constructed a commutative square
\begin{align*}
&\ev_{c_1+q_1-n}^*(M(\mu^{(1)}))
\boxtimes\cdots\boxtimes
\ev_{c_l+q_l-n}^*(M(\mu^{(l)}))\\
&\hspace{13mm}\parallel\\
&\begin{CD}
\Hom_{\mathfrak{p}}(V_{\bc}^{\otimes d}, Z / Z_{> b})
&\:\:\hookrightarrow\:\:&\Hom_{W(\lambda)}(V_{\bc}^{\otimes d}, Z / Z_{> b})\\
@VVV@VVV\\
\Hom_{\mathfrak{p}}(V_{\bc}^{\otimes d}, V(A))
&\:\:\hookrightarrow\:\:&\Hom_{W(\lambda)}(V_{\bc}^{\otimes d}, V(A))
\end{CD}\\
&\hspace{13mm}\parallel\\
&S(\mu^{(1)})\boxtimes\cdots\boxtimes S(\mu^{(l)})
\end{align*}
where the vertical maps are surjections arising
from the map $Z / Z_{> b} \twoheadrightarrow V(A)$.
This shows that 
there is an $H_{(d_1,\dots,d_l)}$-module homomorphism
$$
\ev_{c_1+q_1-n}^*(S(\mu^{(1)})) \boxtimes\cdots\boxtimes 
\ev_{c_l+q_l-n}^*(S(\mu^{(l)}))
\hookrightarrow \Hom_{W(\lambda)}(V_{\bc}^{\otimes d}, V(A))
$$
whose image generates $\Hom_{W(\lambda)}(V_{\bc}^{\otimes d}, V(A))$ 
as an $H_d$-module.
Applying adjointness of tensor and hom, we get from this a surjective
$H_d$-module homomorphism 
$S(A)\twoheadrightarrow \Hom_{W(\lambda)}(V_{\bc}^{\otimes d}, V(A))$.
In view of (\ref{df}) and another dimension calculation for the right hand side, this is in fact an isomorphism, and we are done
in the special case.

Now for the general case, let $\bar\lambda$ be obtained from $\lambda$
by removing rows as usual.
Take $A \in \Col^d_\bc(\lambda)$ and assume $A$ is contractible,
so $\bar A \in \Col^d_\bc(\bar\lambda)$ is defined.
By Lemma~\ref{sim}, we know that $e V(A) \cong V(\bar A)$.
Hence,  by the last part of Theorem~\ref{block},
we have that 
$\Hom_{W(\bar\lambda)}(\bar V_\bc^{\otimes d}, V(\bar A))
\cong 
f \Hom_{W(\lambda)}(V_\bc^{\otimes d}, V(A))$.
It remains to observe 
using the last statement from Lemma~\ref{sim} that
all composition factors of $V(A)$ are of the form 
$L(B)$ with $B$ contractible.
Hence, recalling (\ref{sand}), we have that $f V(A) = V(A)$.
So, we deduce that
\begin{align*}
f \Hom_{W(\lambda)}(V_\bc^{\otimes d}, V(A))
&= \Hom_{W(\lambda)}(V_\bc^{\otimes d}, f V(A))\\
&=
\Hom_{W(\lambda)}(V_\bc^{\otimes d}, V(A))
\cong S(A)
\end{align*}
as $H_d(\Lambda)$-modules.
\end{proof}

\begin{cor}\label{filt}
For $A \in \Tab^d(\lambda)$, the permutation module $M(A,\bc)$ 
has a Specht flag in which the
Specht module $S(B)$  appears with multiplicity $K_{B,A}$
for each $B \in \Col^d_\bc(\lambda)$, arranged in any order refining the
Bruhat ordering on $\Col^d_\bc(\lambda)$ (most dominant at the bottom).
\end{cor}

\begin{proof}
Immediate from Theorems~\ref{final}, \ref{perm1} and 
Lemma~\ref{s}.
\end{proof}

\begin{rem}\rm
The Specht module
$S(A)$ from (\ref{dsm}) is isomorphic to the degenerate analogue of the
cell module $S^\bmu$ 
of Dipper, James and Mathas.
An explicit construction of a filtration like in Corollary~\ref{filt}
can be deduced from \cite[Theorem 6.3]{AMR} using the 
formalism of cellular algebras. 
\end{rem}

For $A \in \Std^d_\bc(\lambda)$, we let
\begin{align}
D(A) &= \Hom_{\mathfrak{g}}(P_\bc\otimes V^{\otimes d}, K(A))
\cong \Hom_{W(\lambda)}(V_\bc^{\otimes d}, L(A)).
\end{align}
By Theorem~\ref{cta} and Lemma~\ref{ic}, the modules
$\{D(A)\mid A \in \Std^d_\bc(\lambda)\}$ form a complete set of pairwise
non-isomorphic irreducible $H_d(\lambda,\bc)$-modules.

\begin{thm}
For $A \in \Std_\bc^d(\lambda)$,
the Specht module $S(A)$ has an irreducible 
cosocle isomorphic to $D(A)$. Moreover, for
$A \in \Col^d_\bc(\lambda)$ and $B \in \Std^d_\bc(\lambda)$,
we have that
$$
[S(A):D(B)]=[N(A):K(B)],
$$ 
so the composition multiplicities of Specht modules 
can be computed in terms of Kazhdan-Lusztig polynomials.
\end{thm}

\begin{proof}
Apply Theorem~\ref{final} and 
exactness of the functor
$\Hom_{\mathfrak{g}}(P_\bc \otimes V^{\otimes d}, ?)$.
\end{proof}

\begin{rem}\rm
In the case $\bc = \biz$, 
an explicit formula for the composition multiplicity
$[N(A):K(B)]$ in terms of Kazhdan-Lusztig polynomials
is recorded in the formula \cite[(4.8)]{BKrep}: the expression in parentheses
there is exactly $[N(A):K(B)]$ by \cite[Theorem 4.5]{BKrep}.
We refer the reader to \cite[$\S$8.5]{BKrep}
for further combinatorial results about the representation theory
of $W_d(\lambda,\bc)$, all of which can easily be
translated into analogous (but mostly already known) results 
about the representation theory of $H_d(\lambda,\bc)$ using Theorem~\ref{cta}
and Lemma~\ref{sun}. 
We point out finally by \cite[Theorem 6.11]{AMR} that 
$H_d(\lambda,\bc)$, hence also 
the higher level Schur algebra
$W_d(\lambda,\bc)$, is a semisimple algebra
if and only if 
$d \leq q_i - q_j$
for all $1 \leq i < j \leq l$ with $c_i = c_j$.
\end{rem}

\appendix
\section{Symmetric algebras}

In this appendix we include a short proof of the fact that 
the degenerate cyclotomic Hecke algebra $H_d(\Lambda)$ is a 
symmetric algebra.
The $q$-analogue of this result is due to 
Malle and Mathas \cite{MM}. The degenerate case is actually much easier
and, in view of Theorem~\ref{thmbreal},
is closely related to a general theorem
of Mazorchuk and Stroppel \cite[Theorem 4.6]{MS}.
We will continue to work over the ground field $\C$
as above, but note that the same argument works over an arbitrary commutative 
ground ring $R$.
Recall that a 
finite dimensional algebra $A$
is {\em symmetric} if it possesses a {\em symmetrizing form},
i.e. a linear map
$\tau: A \rightarrow \C$ such that 
 $\tau(ab) = \tau(ba)$ for all $a,b \in A$
and 
whose kernel contains no non-zero left or right ideal of $A$.
The following lemma is an easy exercise.

\begin{lem}\label{a1}
Let $\tau:
\C_l[x_1,\dots,x_d] \,{\scriptstyle{\rtimes\!\!\!\!\!\bigcirc}}\, \C S_d
\rightarrow \C$ be the linear map
sending the monomial $x_1^{r_1} \cdots x_d^{r_d} w$
to $1$ if $r_1=\cdots = r_d = l-1$ and $w = 1$,
or to $0$ otherwise, for all $r_1,\dots,r_d \geq 0$ and $w \in S_d$.
Then $\tau$ is a symmetrizing form,
hence
$\C_l[x_1,\dots,x_d] \,{\scriptstyle{\rtimes\!\!\!\!\!\bigcirc}}\, \C S_d$
is a symmetric algebra.
\end{lem}

Now let $t = (l-1)^d$.
Recall from Lemma~\ref{ze} that $H_d(\Lambda)$ is a filtered algebra
with filtration
$$
\C S_d = \operatorname{F}_0 H_d(\Lambda)
\subseteq
\operatorname{F}_1 H_d(\Lambda)
\subseteq \cdots \subseteq
\operatorname{F}_t H_d(\Lambda) = H_d(\Lambda).
$$
The associated graded algebra
$\gr H_d(\Lambda)$ is identified with 
the twisted tensor product 
$\C_l[x_1,\dots,x_d] \,{\scriptstyle{\rtimes\!\!\!\!\!\bigcirc}}\, \C S_d$.
For any $0 \leq s \leq t$, 
let $$
\gr_s: \operatorname{F}_s H_d(\Lambda)
\rightarrow 
\C_l[x_1,\dots,x_d] \,{\scriptstyle{\rtimes\!\!\!\!\!\bigcirc}}\, \C S_d
$$ 
be the
map sending an element to its degree $s$ graded component.

\begin{thm}
Let $\hat\tau:H_d(\Lambda)
\rightarrow \C$ be the linear map
sending $x_1^{r_1} \cdots x_d^{r_d} w$
to $1$ if $r_1=\cdots = r_d = l-1$ and $w = 1$,
or to $0$ otherwise, for all $0 \leq r_1,\dots,r_d < l$ and $w \in S_d$.
Then $\hat\tau$ is a symmetrizing form,
hence $H_d(\Lambda)$ is a symmetric algebra.
\end{thm}

\begin{proof}
The key observation is that 
$\hat\tau = \tau \circ \gr_t$ where
$\tau$ is symmetrizing form
from Lemma~\ref{a1}.
To prove that $\hat\tau$ is symmetric, observe for 
$0 \leq s \leq t$ and any 
$x \in \operatorname{F}_s H_d(\Lambda)$
and $y \in \operatorname{F}_{t-s} H_d(\Lambda)$
that 
\begin{align*}
\hat\tau(xy) &= \tau(\gr_t(xy))=
\tau(\gr_s(x) \gr_{t-s}(y))\\
&=
\tau(\gr_{t-s}(y)\gr_s(x))
=\tau(\gr_t(yx)) = \hat\tau(yx).
\end{align*}
In particular, this shows that
$\hat\tau(xw) = \hat \tau(wx)$ 
for all $x \in H_d(\Lambda)$ and $w\in S_d$.
It just remains to show that
$\hat\tau(x_i y) = \hat\tau(yx_i)$ for any $i=1,\dots,d$
and $y \in H_d(\Lambda)$.
We are already done if $y \in \operatorname{F}_{t-1} H_d(\Lambda)$,
so we may assume that $y = x_1^{l-1} \cdots x_d^{l-1} w$
for some $w \in S_d$.
But then
\begin{align*}
\hat\tau(x_i y)
&= 
\hat\tau(x_i x_1^{l-1} \cdots x_d^{l-1} w)
=
\hat\tau(w x_i 
x_1^{l-1} \cdots x_d^{l-1})\\
&=
\hat\tau(w 
x_1^{l-1} \cdots x_d^{l-1}x_i )
=
\hat\tau(
x_1^{l-1} \cdots x_d^{l-1} w x_i )
= \hat\tau(yx_i),
\end{align*}
using the fact that $x_1^{l-1} \cdots x_d^{l-1}$ is central in $H_d(\Lambda)$.
Now suppose that $I$ is a left or right ideal of $H_d(\Lambda)$
such that $\hat\tau(I) = 0$.
Define a filtration on $I$ by setting $\operatorname{F}_s I= 
I \cap \operatorname{F}_s H_d(\Lambda)$.
Then $\gr I$ is a left or right ideal in
$\C_l[x_1,\dots,x_d] \,{\scriptstyle{\rtimes\!\!\!\!\!\bigcirc}}\, \C S_d$ 
such that $\tau(\gr I) = 0$.
Hence $\gr I = 0$, so $I = 0$ too.
\end{proof}

\end{document}